\documentclass[b5paper,twoside]{irmaems}
\pdfoutput=1
\usepackage{amssymb} 
\usepackage{amsmath}
\usepackage{makeidx}
\usepackage{latexsym}
\usepackage{subfigure}
\usepackage{graphicx}
\usepackage{xy}
\xyoption{all}

\newcommand{\B}{\mathcal{B}}
\newcommand{\CP}{\mathbb{CP}}
\newcommand{\T}{\mathcal{T}}
\renewcommand{\P}{\mathcal{P}}

\newcommand{\ML}{\mathcal{ML}}
\newcommand{\PML}{\mathbb{P}\mathcal{ML}}

\newcommand{\V}{\mathcal{V}}
\newcommand{\Z}{\mathbb{Z}}
\newcommand{\R}{\mathbb{R}}
\newcommand{\C}{\mathbb{C}}
\newcommand{\N}{\mathbb{N}}
\newcommand{\E}{\mathcal{E}}
\newcommand{\Q}{\mathcal{Q}}
\newcommand{\QF}{\mathcal{QF}}
\newcommand{\PQ}{\mathbb{P}^+Q}
\newcommand{\D}{\mathcal{D}}
\newcommand{\F}{\mathcal{F}}
\newcommand{\scc}{\mathcal{S}}
\renewcommand{\H}{\mathbb{H}}
\newcommand{\oo}{G}
\newcommand{\planes}{\mathcal{H}^{2,1}}
\newcommand{\rep}{\mathcal{R}}
\newcommand{\X}{\mathcal{X}}
\newcommand{\sslash}{\slash\!\!\slash}
\newcommand{\free}{\ast}

\newcommand{\param}{{\mathchoice{\mkern1mu\mbox{\raise2.2pt\hbox{$\centerdot$}}\mkern1mu}{\mkern1mu\mbox{\raise2.2pt\hbox{$\centerdot$}}\mkern1mu}{\mkern1.5mu\centerdot\mkern1.5mu}{\mkern1.5mu\centerdot\mkern1.5mu}}}

\newcommand{\tensor}{\otimes}
\renewcommand{\bar}[1]{\overline{#1}}

\renewcommand{\Im}{\imag}

\renewcommand{\L}{\mathrm{L}}
\newcommand{\PL}{\mathbb{P}\L}

\DeclareMathOperator{\Hess}{\mathrm{Hess}}

\DeclareMathOperator{\imag}{\mathrm{Im}}

\newcommand{\Pl}{\mathrm{Pl}}
\newcommand{\PSL}{\mathrm{PSL}}
\newcommand{\SL}{\mathrm{SL}}
\newcommand{\SU}{\mathrm{SU}}
\newcommand{\SO}{\mathrm{SO}}

\renewcommand{\sl}{\mathfrak{sl}}
\newcommand{\osc}{\mathrm{osc}}
\newcommand{\gr}{\mathrm{gr}}
\newcommand{\Gr}{\mathrm{Gr}}
\newcommand{\Hom}{\mathrm{Hom}}
\newcommand{\tr}{\mathrm{tr}}
\newcommand{\hol}{\mathrm{hol}}
\newcommand{\wrap}{\mathrm{wr}}
\newcommand{\interior}{\mathrm{int}}
\renewcommand{\Tilde}{\widetilde}

\newcommand{\npc}{\textsc{NPC}}
\newcommand{\catzero}{$\text{\sc CAT}(0)$}
\newcommand{\catk}{$\text{\sc CAT}(\kappa)$}

\newcommand{\boldpoint}[1]{\medskip\par\noindent\textbf{#1}}

\renewcommand\theequation{\arabic{section}.\arabic{equation}}

\theoremstyle{definition} %%% for statements in roman typeface

 \newtheorem{definition}{Definition}[section]

\theoremstyle{plain}      %%% for statements in italic typeface

 \newtheorem{theorem}[definition]{Theorem}
 \newtheorem{corollary}[definition]{Corollary}
 \newtheorem{lemma}[definition]{Lemma}

\begin{document}

\title{Complex Projective Structures}

\author{David Dumas\thanks{
Work partially supported by a NSF postdoctoral research
fellowship.}}

\address{
Department of Mathematics, Statistics, and Computer Science\\
University of Illinois at Chicago\\[6mm]
January 24, 2009
}

\maketitle

\vspace{-16mm}

\tableofcontents

\section{Introduction}\label{sec:introduction}

In this chapter we discuss the theory of complex projective structures
on compact surfaces and its connections with Teichm\"uller theory,
$2$- and $3$-dimensional hyperbolic geometry, and representations of
surface groups into $\PSL_2(\C)$.  Roughly speaking, a complex
projective structure is a type of $2$-dimensional geometry in which
M\"obius transformations play the role of geometric congruences (this
is made precise below).  Such structures are abundant---hyperbolic,
spherical, and Euclidean metrics on surfaces all provide examples of
projective structures, since each of these constant-curvature
$2$-dimensional geometries has a model in which its isometries are
M\"obius maps.  However, these examples are not representative of the
general situation, since most projective structures are not induced by
locally homogeneous Riemannian metrics.

Developing a more accurate picture of a general projective structure is
the goal of the first half of the chapter (\textbf{\S\S2--4}).  After
some definitions and preliminary discussion (in \textbf{\S2}), we
present the complex-analytic theory of projective structures in
\textbf{\S3}.  This theory has its roots in the study of automorphic
functions and differential equations by Klein \cite[Part~1]{klein:ausgewahlte},
Poincar\'e \cite{poincare:groups-linear-equations}, Riemann
\cite{riemann:hypergeometric}, and others in the late
nineteenth century (see \cite{hejhal:monodromy-and-poincare}
\cite[\S1]{hejhal:monodromy} for further historical discussion and
references), while its more recent history is closely linked to
developments in Teichm\"uller theory and deformations of Fuchsian and
Kleinian groups (e.g.~\cite{earle:variation}
\cite{gunning:special-coordinate} \cite{gunning:riemann-surfaces} \cite{hejhal:monodromy}
\cite{hubbard:monodromy} \cite{kra:deformations}
\cite{kra:deformations2} \cite{kra:generalization} \cite{kra-maskit:remarks}).

In this analytic approach, a projective structure is represented by a
holomorphic quadratic differential on a Riemann surface, which is
extracted from the geometric data using a M\"obius-invariant
differential operator, the \emph{Schwarzian derivative}.  The inverse
of this construction describes every projective structure in terms of
holomorphic solutions to a linear ordinary differential equation (the
\emph{Schwarzian equation}).  In this way, many properties of
projective structures and their moduli can be established using tools
from complex function theory.  However, in spite of the success of
these techniques, the analytic theory is somewhat detached from the
underlying geometry.
In particular, the analytic parameterization of projective structures
does not involve an explicit geometric construction, such as one has
in the description of hyperbolic surfaces by gluing
polygons.

In \textbf{\S4} we describe a more direct and geometric construction
of complex projective structures using \emph{grafting}, a gluing
operation on surfaces which is also suggested by the work of the
nineteenth-century geometers (e.g.~\cite{klein}), but whose
significance in complex projective geometry has only recently been fully
appreciated.  Grafting was used by Maskit \cite{maskit:grafting},
Hejhal \cite{hejhal:monodromy}, and Sullivan-Thurston
\cite{sullivan-thurston} to construct certain deformations of Fuchsian
groups, and in later work of Thurston (unpublished, see
\cite{kamishima-tan:grafting}) it was generalized to give a universal
construction of complex projective surfaces starting from basic
hyperbolic and Euclidean pieces.

This construction provides another coordinate system for the moduli
space of projective structures, and it reveals an important connection
between these structures and convex geometry in $3$-dimensional
hyperbolic space.  However, the explicit geometric nature of complex
projective grafting comes at the price of a more complicated parameter
space, namely, the piecewise linear manifold of \emph{measured
geodesic laminations} on hyperbolic surfaces.  In particular, the lack of a
differentiable structure in this coordinate system complicates the
study of variations of complex projective structures, though there has
been some progress in this direction using a weak notion of
differentiability due to Thurston \cite{thurston:minimal-stretch} and
Bonahon \cite{bonahon:variations}.

After developing the analytic and geometric coordinates for the moduli
space of projective structures, the second half of the chapter is
divided into two major topics: In \textbf{\S5}, we describe the
relation between projective structures and the
$\PSL_2(\C)$-representations of surface groups, their deformations,
and associated problems in hyperbolic geometry and Kleinian groups.
The key to these connections is the \emph{holonomy representation} of
a projective structure, which records the topological obstruction to
analytically continuing its local coordinate charts over the entire
surface.  After constructing a parameter space for such
representations and the holonomy map for projective structures, we
survey various developments that center around two basic questions:
\begin{itemize}
\item Given a projective structure, described in either analytic or
  geometric terms, what can be said about its holonomy representation?
\item Given a $\PSL_2(\C)$-representation of a surface group, what
  projective structures have this as their holonomy representation, if
  any?
\end{itemize}
We discuss partial answers to these general questions, along with much
more detailed information about certain classes of holonomy
representations (e.g.~Fuchsian groups).

Finally, in \textbf{\S6} we take up the question of relating the
analytic and geometric coordinate systems for the space of projective
structures, or equivalently, studying the interaction between the
Schwarzian derivative and complex projective grafting.  We describe
asymptotic results that relate compactifications of the analytic and
geometric parameter spaces using the geometry of measured foliations
on Riemann surfaces.  Here a key tool is the theory of harmonic
maps\index{harmonic map} between Riemann surfaces and from Riemann
surfaces to $\R$-trees, and the observation that two geometrically
natural constructions in complex projective geometry (the
\emph{collapsing} and \emph{co-collapsing} maps) are closely
approximated by harmonic maps.  We close with some remarks concerning
infinitesimal compatibility between the geometric and analytic
coordinate systems, once again using the limited kind of differential
calculus that applies to the grafting parameter space.

\boldpoint{Scope and approach.}  Although this chapter covers a range
of topics in complex projective geometry, is not intended to be a
comprehensive guide to the subject.  Rather, we have selected several
important aspects of the theory (the Schwarzian derivative, grafting,
and holonomy) and concentrated on describing their interrelationships
while providing references for further reading and exploration.  As a
result, some major areas of research in complex projective structures
are not mentioned at all (circle packings \cite{kojima-mizushima-tan}
\cite{kojima:survey}, the algebraic-geometric aspects of the theory
\cite[\S11]{gkm}, and generalizations to punctured or open Riemann
surfaces \cite{kra:deformations2} \cite{luo}, to name a few) and
others are only discussed in brief.

We have also included some detail on the basic analytic and geometric
constructions in an attempt to make this chapter a more useful
``invitation'' to the theory.  However, where we discuss more advanced
topics and results of recent research, it has been necessary to refer
to many concepts and results that are not thoroughly developed here.

Finally, while we have attempted to provide thorough and accurate
references to the literature, the subject of complex projective
structures is broad enough (and connected to so many other areas of
research) that we do not expect these references to cover every
relevant source of additional information.  We hope that the
references included below are useful, and regret any inadvertent
omissions.

\boldpoint{Acknowledgments.}  The author thanks Richard Canary, George
Daskalopoulos, William Goldman, Brice Loustau, Albert Marden, Athanase
Papadopoulos, Richard Wentworth, and Michael Wolf for helpful
discussions and suggestions related to this work, and Curt McMullen
for introducing him to the theory of complex projective structures.

\section{Basic definitions}\label{sec:basics}

\boldpoint{Projective structures.}  Let $S$ be an oriented surface.  A
\emph{complex projective structure}\index{complex projective
  structure}\index{projective structure!complex} $Z$ on $S$ is a
maximal atlas of charts mapping open sets in $S$ into $\CP^1$ such
that the transition functions are restrictions of M\"obius
transformations.  For brevity we also call these \emph{projective
  structures} or \emph{$\CP^1$-structures}.

We often treat a projective structure $Z$ on $S$ as a surface in its
own right---a \emph{complex projective surface}.
Differentiably, $Z$ is the same
as $S$, but $Z$ has the additional data of a restricted atlas of
projective charts.

Two projective structures $Z_1$ and $Z_2$ on $S$ are
\emph{isomorphic}\index{isomorphism!of projective structures} if there is
an orientation-preserving diffeomorphism $\iota : Z_1 \to Z_2$ that
pulls back the projective charts of $Z_2$ to projective charts of
$Z_1$, and \emph{marked isomorphic} if furthermore $\iota$ is
homotopic to the identity.

Our main object of study is the space $\P(S)$ of marked isomorphism
classes of projective structures on a compact surface $S$.  Thus far,
we have only defined $\P(S)$ as a set, but later we will equip it with
the structure of a complex manifold.

\boldpoint{Non-hyperbolic cases.}
Projective structures on compact surfaces are most interesting when
$S$ has genus $g \geq 2$: The sphere has a unique projective structure
(by $S^2 \simeq \CP^1$) up to isotopy, while a projective structure on
a torus is always induced by an affine structure \cite[\S9, pp.~189-191]{gunning:riemann-surfaces}.  We therefore make the
assumption that $S$ has genus $g \geq 2$ unless stated otherwise.

\boldpoint{First examples.}
The projective structure of $\CP^1$ itself (using the identity for
chart maps) also gives a natural projective structure on any open set
$U \subset \CP^1$.  If $U$ is preserved by a group $\Gamma$
of M\"obius transformations acting freely and properly
discontinuously, then the quotient surface $X = U / \Gamma$ has a
natural projective structure in which the charts are local inverses of
the covering $U \to X$.

In particular any Fuchsian group $\Gamma \subset \PSL_2(\R)$ gives
rise to a projective structure on the quotient surface $\H / \Gamma$
and a Kleinian group $\Gamma \subset \PSL_2(\C)$ gives a projective
structure on the quotient of its domain of discontinuity
$\Omega(\Gamma) / \Gamma$.  Rephrasing the latter example, the ideal
boundary of a hyperbolic $3$-manifold has a natural projective
structure.

\boldpoint{Locally M\"obius maps.}
A map $f : Z \to W$ between complex projective surfaces is
\emph{locally M\"obius} if for every sufficiently small open set $U
\subset Z$, the restriction $\left . f \right|_U$ is a M\"obius
transformation with respect to projective coordinates on $U$ and
$f(U)$.
Examples of such maps include isomorphisms and covering maps of
projective surfaces (where the cover is given the pullback projective
structure) and inclusions of open subsets of surfaces.

\boldpoint{Developing maps.}
A projective structure $Z$ on a surface $S$ lifts to a projective
structure $\Tilde{Z}$ on the universal cover $\Tilde{S}$.
A \emph{developing map}\index{developing map} for $Z$ is an immersion $f :
\Tilde{S} \to \CP^1$ such that the restriction of $f$ to any
sufficiently small open set in $\Tilde{S}$ is a projective chart for
$\Tilde{Z}$.  Such a map is also called a \emph{geometric realization}
of $Z$ (e.g.~\cite[\S6]{gunning:special-coordinate}) or a \emph{fundamental
  membrane} \cite{hejhal:monodromy}.

Developing maps always exist, and are essentially unique---two
developing maps for a given structure differ by post-composition with
a M\"obius transformation.  Concretely, a developing map can be
constructed by analytic continuation starting from any basepoint $z_0
\in \Tilde{Z}$ and any chart defined on a neighborhood $U$ of $z_0$.
Another chart $V \to \CP^1$ that overlaps $U$ can be adjusted by a
M\"obius transformation so as to agree on the overlap, gluing to give
a map $(U \cup V) \to \CP^1$.  Continuing in this way one defines a
map on successively larger subsets of $\Tilde{Z}$, and the limit is a
developing map\index{developing map} $\Tilde{Z} \to \CP^1$.  The simple connectivity of
$\Tilde{Z}$ is essential here, as nontrivial homotopy classes of loops
in the surface create obstructions to unique analytic continuation of
a projective chart.

For a fixed projective structure, we will speak of \emph{the}
developing map when the particular choice is unimportant or implied.

\boldpoint{Holonomy representation.}
The developing map $f : \Tilde{S} \to \CP^1$ of a projective structure
$Z$ on $S$ has an equivariance property with respect to the action of
$\pi_1(S)$ on $\Tilde{S}$: For any $\gamma \in \pi_1(S)$, the composition
$f \circ \gamma$ is another developing map for $Z$.  Thus there exists
$A_\gamma \in \PSL_2(\C)$ such that
\begin{equation}
\label{eqn:equivariance}
f \circ \gamma = A_\gamma \circ f
\end{equation}
The map $\gamma \mapsto A_\gamma$ is a homomorphism $\rho : \pi_1(S)
\to \PSL_2(\C)$, the \emph{holonomy representation}%
\index{holonomy representation} (or \emph{monodromy representation}%
\index{monodromy representation}) of the projective structure.

\boldpoint{Development-holonomy pairs.}  The developing
map\index{developing map} and holonomy representation form the
\emph{development-holonomy pair}\index{development-holonomy pair}
$(f,\rho)$ associated to the projective structure $Z$.

This pair determines $Z$ uniquely, since restriction of $f$ determines
a covering of $S$ by projective charts.  Post-composition of the
developing map with $A \in \PSL_2(\C)$ conjugates $\rho$, and
therefore the pair $(f,\rho)$ is uniquely determined by $Z$ up to the
action of $\PSL_2(\C)$ by
\begin{equation*}
(f,\rho) \mapsto (A \circ f, \rho^A) \; \;\text{where} \;\;
  \rho^A(\gamma) = A \rho(\gamma) A^{-1}
\end{equation*}

Conversely, any pair $(f,\rho)$ consisting of an immersion $f :
\Tilde{S} \to \CP^1$ and a homomorphism $\rho : \pi_1(S) \to \PSL_2(\C)$
that satisfy \eqref{eqn:equivariance} defines a projective structure
on $S$ in which lifting $U \subset S$ to $\Tilde{S}$ and applying $f$
gives a projective chart (for all sufficiently small open sets $U$).

Thus we have an alternate definition of $\P(S)$ as the quotient of the
set of development-holonomy pairs by the $\PSL_2(\C)$ action and
by precomposition of developing maps with orientation-preserving
diffeomorphisms of $S$ homotopic to the identity.  We give the set of
pairs of maps $(f,\rho)$ the compact-open topology, and
$\P(S)$ inherits a quotient topology.
We will later see that $\P(S)$ is homeomorphic to $\R^{12g-12}$.

\boldpoint{Relation to $(G,X)$-structures.}  There is a very general
notion of a geometric structure defined by a Lie group $G$ acting by
diffeomorphisms on a manifold $X$.  A
\emph{$(G,X)$-structure}\index{G-X-structure@$(G,X)$-structure} on a
manifold $M$ is an atlas of charts mapping open subsets of $M$ into $X$
such that the transition maps are restrictions of elements of $G$.

In this language, complex projective structures are $(\PSL_2(\C),
\CP^1)$-structures.  Some of the properties of projective structures
we develop, such as developing maps\index{developing map}, holonomy
representations\index{holonomy representation}, deformation spaces,
etc., can be applied in the more general setting of
$(G,X)$-structures.  See \cite{goldman:geometric-structures} for a
survey of $(G,X)$-structures and analysis of several low-dimensional
examples.

\boldpoint{Circles.}  Because M\"obius transformations map circles to
circles, there is a natural notion of a circle on a surface with a
projective structure $Z$: A smooth embedded curve $\alpha \subset Z$
is a \emph{circular arc} if the projective charts
map (subsets of) $\alpha$ to circular arcs in $\CP^1$.  Equivalently,
the embedded curve $\alpha$ is a circular arc if the developing map
sends any connected component of the preimage of $\alpha$ in
$\Tilde{Z}$ to a circular arc in $\CP^1$.  A closed circular arc on
$Z$ is a \emph{circle}\index{circle}.

Small circles are ubiquitous in any projective structure: For any $z
\in Z$ there is a projective chart mapping a contractible neighborhood
of $z$ to an open set $V \in \CP^1$.  The preimage of any circle
contained in $V$ is a homotopically trivial circle for the projective
structure $Z$.  Circles that bound disks on a projective surface have
an important role in Thurston's projective grafting construction (see
\S\ref{sec:projective-grafting}).

Circles on a projective surface can also be homotopically nontrivial.
For example, any simple closed geodesic on a hyperbolic surface is a
circle, because its lifts to $\H$ are half-circles or vertical lines
in the upper half-plane.  The analysis of circles on more general
projective surfaces would be a natural starting point for the
development of synthetic complex projective geometry; Wright's study
of circle chains and Schottky-type dynamics in the Maskit slice of
punctured tori is an example of work in this direction
\cite{wright:maskit-boundary}.

\boldpoint{Forgetful map.}
Since M\"obius transformations are holomorphic, a projective structure
$Z \in \P(S)$ also determines a complex structure, making $S$ into a
compact Riemann surface.  In this way, marked isomorphism
of projective structures corresponds to marked isomorphism of Riemann
surfaces, and so there is a natural (and continuous) \emph{forgetful
  map}\index{forgetful map}
\begin{equation*}
\pi : \P(S) \to \T(S)
\end{equation*}
where $\T(S)$ is the Teichm\"uller space of marked isomorphism classes
of complex structures on $S$.  (See e.g.~\cite{lehto:univalent},
\cite{imayoshi-taniguchi}, \cite{hubbard:book} for background on
Teichm\"uller spaces.)  As a matter of terminology, if $Z$ is a
projective structure with $\pi(Z) = X$, we say $Z$ is a projective
structure \emph{on the Riemann surface} $X$.

The forgetful map is surjective: By the uniformization theorem, every
complex structure $X \in \T(S)$ arises as the quotient of $\H$ by a
Fuchsian group $\Gamma_X$, and the natural projective structure on $\H
/ \Gamma_X$ is a preimage of $X$ by $\pi$.  We call this the
\emph{standard Fuchsian structure} on $X$%
\index{projective structure!standard Fuchsian}.  The standard Fuchsian
structures determine a continuous section
\begin{equation*}
\sigma_0 : \T(S) \to \P(S).
\end{equation*}

One might expect the fibers of $\pi$ to be large, since isomorphism of
projective structures is a much stronger condition than isomorphism of
complex structures.  Our next task is to describe the fibers
explicitly.

\section{The Schwarzian Parameterization}\label{sec:schwarzian}\index{Schwarzian parameterization}

\subsection{The Schwarzian derivative}

Let $\Omega \subset \C$ be a connected open set.  The \emph{Schwarzian
  derivative}\index{Schwarzian derivative} of a locally injective
holomorphic map $f : \Omega \to \CP^1$ is the holomorphic quadratic
differential
\begin{equation*}
S(f) = \left [ \left(\frac{f''(z)}{f'(z)}\right)' -
\frac{1}{2} \left(\frac{f''(z)}{f'(z)}\right)^2 \right ] \: dz^2.
\end{equation*}
Two key properties make the Schwarzian derivative\index{Schwarzian derivative} useful in the theory
of projective structures:
\begin{enumerate}
\item \emph{Cocycle property.}  If $f$ and $g$ are locally injective
  holomorphic maps such that the composition $f
  \circ g$ is defined, then
\begin{equation*}
S(f \circ g) = g^* S(f) + S(g)
\end{equation*}
\item \emph{M\"obius invariance.} For any $A \in \PSL_2(\C)$, we have
\begin{equation*}
S(A) \equiv 0,
\end{equation*}
and conversely, if $S(f) \equiv 0$, then $f$ is the restriction of a
M\"obius transformation.
\end{enumerate}

\smallskip

Note that the pullback $g^* S(f)$ uses the definition of the
Schwarzian as a quadratic differential.  In classical complex
analysis, the Schwarzian was regarded as a complex-valued
function, with $g^* S(f)$ replaced by $g'(z)^2 S(f)(g(z))$.

An elementary consequence of these properties is that the map $f$ is
almost determined by its Schwarzian derivative%
\index{Schwarzian derivative}; if $S(f) = S(g)$, then the locally
defined map $f \circ g^{-1}$ satisfies $S(f \circ g^{-1}) \equiv 0$,
and so we have $f = A \circ g$ for some $A \in \PSL_2(\C)$.

Further discussion of the Schwarzian derivative can be found in
e.g.~\cite[Ch.~2]{lehto:univalent} \cite[\S6.3]{hubbard:book}.

\boldpoint{Osculation.}
\label{sec:osculation}
Intuitively, the Schwarzian derivative measures the failure of a
holomorphic map to be the restriction of a M\"obius transformation.
Thurston made this intuition precise as follows (see
\cite[\S2]{thurston:zippers}, \cite[\S2.1]{anderson:projective}):  For
each $z \in \Omega$, there is a unique M\"obius transformation that
has the same $2$-jet as $f$ at $z$, called the \emph{osculating
  M\"obius transformation}\index{osculation} $\osc_z f$.

The \emph{osculation map} $\oo : \Omega \to \PSL_2(\C)$ given by $\oo(z) = \osc_zf$ is
holomorphic, and its Darboux derivative (see \cite{sharpe}) is the holomorphic $\sl_2(\C)$-valued
$1$-form
\begin{equation*}
\omega(z) = \oo^{-1}(z) \: d \oo(z).
\end{equation*}
An explicit computation shows that $\omega$ only depends on $f$ through its
Schwarzian derivative; if $S(f) = \phi(z) dz^2$, then
\begin{equation*}
\omega(z) = -\frac{1}{2} \phi(z)
\begin{pmatrix}
z & -z^2\\
1 & -z
\end{pmatrix}
\: dz.
\end{equation*}

\subsection{Schwarzian parameterization of a fiber}

\boldpoint{Fibers over Teichm\"uller space}
For any marked complex structure $X \in \T(S)$, let $P(X) =
\pi^{-1}(X) \subset \P(S)$ denote the set of marked complex projective
structures with underlying complex structure $X$.  The Schwarzian
derivative can be used to parameterize the fiber $P(X)$ as follows:

Fix a conformal identification $\Tilde{X} \simeq \H$, whereby $\pi_1(S)$
acts on $\H$ as a Fuchsian group.  Abusing notation, we use the same
symbol for $\gamma \in \pi_1(S)$ and for its action on $\H$ by a real
M\"obius transformation.

Given $Z \in P(X)$, we regard the developing map\index{developing map} as a meromorphic
function $f$ on $\H$.  The Schwarzian derivative $\Tilde{\phi} = S(f)$
is therefore a holomorphic quadratic differential on $\H$.  Combining
the equivariance property \eqref{eqn:equivariance} of $f$ and the
properties of the Schwarzian derivative, we find
$$ \Tilde{\phi} = S(A_\gamma \circ f) = S(f \circ \gamma) = \gamma^*
\Tilde{\phi},$$ 
Thus we have $\Tilde{\phi} = \gamma^* \Tilde{\phi}$ for all $\gamma
\in \pi_1(S)$, and $\Tilde{\phi}$ descends to a holomorphic quadratic
differential $\phi$ on $X$.  We call $\phi$ the \emph{Schwarzian of
the projective structure}\index{Schwarzian derivative} $Z$.

Let $Q(X)$ denote the vector space of holomorphic quadratic
differentials on the marked Riemann surface $X \in \T(S)$.  By the
Riemann-Roch theorem, we have $Q(X) \simeq \C^{3g-3}$ (see
\cite{jost:compact-riemann-surfaces}).  The Schwarzian defines a map
$P(X) \to Q(X)$.  We will now show that this map is bijective by
constructing its inverse.

\boldpoint{Inverting the Schwarzian.}  Let $\phi(z)$ be a holomorphic
function defined on a contractible open set $\Omega \subset \C$.  Then
the linear ODE (the \emph{Schwarzian equation}\index{Schwarzian
equation})
\begin{equation}
\label{eqn:schwarzian-ode}
 u''(z) + \frac{1}{2} \phi(z) u(z) = 0
\end{equation}
 has a two-dimensional vector space $V$ of holomorphic solutions on
$\Omega$.  Let $u_1(z)$ and $u_2(z)$ be a basis of solutions.  The
Wronskian $W(z)$ of $u_1$ and $u_2$ satisfies $W'(z) = 0$, so it is a
nonzero constant function, and $u_1$ and $u_2$ cannot vanish
simultaneously.

This ODE construction inverts the Schwarzian derivative in the sense
that the meromorphic function $f(z) = u_1(z) / u_2(z)$ satisfies $S(f)
= \phi(z) \: dz^2$ (see \cite{nehari:schwarzian}).  Note that changing
the basis for $V$ will alter $f$ by composition with a M\"obius
transformation (and leave $S(f)$ unchanged).  Furthermore, since
$$ f'(z) = \frac{u_1'(z) u_2(z) - u_1(z) u_2'(z)}{u_2(z)^2} =
\frac{-W(z)}{u_2(z)^2},$$ 
it follows that the holomorphic map $f : \Omega \to \CP^1$ is
locally injective except possibly on $\{ u_2(z) = 0 \} =
f^{-1}(\infty)$.  Applying similar considerations to $1/f(z)$, we find
that $f$ is locally injective away from $\{ u_1(z) = 0 \}$, and thus
everywhere.

The existence of a holomorphic map with a given Schwarzian derivative
can also be understood in terms of maps to the Lie group $\PSL_2(\C)$
and the definition of the Schwarzian in terms of osculation (described
in \S\ref{sec:osculation})\index{osculation}.  Here the quadratic
differential $\phi$ is interpreted as a $\sl_2(\C)$-valued $1$-form,
which satisfies the integrability condition $d \phi + \frac{1}{2}
[\phi,\phi] = 0$ because there are no holomorphic $2$-forms on a
Riemann surface.  The integrating map to $\PSL_2(\C)$ is the
osculation map of a holomorphic function $f$ satisfying $S(f) = \phi$.
See \cite[\S2.2.3, Cor~2.20]{anderson:projective} for details.

\boldpoint{Parameterization of a fiber.}
Given a quadratic differential $\phi \in Q(X)$, lift to the universal
cover $\Tilde{X} \simeq \H$ to obtain $\Tilde{\phi} =
\Tilde{\phi}(z)\: dz^2$.  Applying the ODE construction to
$\Tilde{\phi}(z)$ yields a holomorphic immersion $f_\phi : \H \to
\CP^1$.

For any $\gamma \in \pi_1(S)$ we have $S(f_\phi \circ \gamma) =
\gamma^* \Tilde{\phi} = \Tilde{\phi} = S(f_\phi)$, and thus $f_\phi
\circ \gamma = A_\gamma \circ f_\phi$ for some $A_\gamma \in
\PSL_2(\C)$.  We set $\rho_\phi(\gamma) = A_\gamma$.  Then
$(f_\phi,\rho_\phi)$ determine a development-holonomy
pair\index{development-holonomy pair}\index{developing map}, and thus
a projective structure $X_\phi$ on $S$.  Since $f$ is holomorphic, we
also have $\pi(X_\phi) = X$.

The map $Q(X) \to P(X)$ given by $\phi \mapsto X_\phi$ is inverse to
the Schwarzian map $P(X) \to Q(X)$ because the ODE construction is
inverse to the Schwarzian derivative.  In particular, each fiber of
$\pi : \P(S) \to \T(S)$ is naturally parameterized by a complex vector
space.

\boldpoint{Affine naturality.}
\label{sec:affine-naturality}
The identification $Q(X) \simeq P(X)$ defined above depends on a
choice of coordinates on the universal cover of $X$.  Specifically, we
computed the Schwarzian using the coordinate $z$ of the upper
half-plane.

A coordinate-independent statement is that the Schwarzian derivative is
a measure of the \emph{difference} between a pair of projective
structures on $X$, which we can see as follows: Given $Z_1, Z_2 \in
P(X)$, let $U$ be a sufficiently small open set on $S$ so that there
are projective coordinate charts $z_i : U \to \CP^1$ of $Z_i$ for
$i=1,2$.  We can assume that $\infty \notin z_i(U)$.

The quadratic differential $z_1^*S(z_2 \circ z_1^{-1})$ on $U$ is
holomorphic with respect to the Riemann surface structure $X$.
Covering $S$ by such sets, it follows from the cocycle property that
these quadratic differentials agree on overlaps and define an element
$\phi \in Q(X)$, which is the \emph{Schwarzian of $Z_2$ relative to
$Z_1$}\index{Schwarzian of a projective structure}.  Abusing notation, we write $Z_2 - Z_1 = \phi$.

Thus $P(X)$ has a natural structure of an \emph{affine space} modeled
on the vector space $Q(X)$.  The choice of a basepoint $Z_0 \in P(X)$
gives an isomorphism $P(X) \to Q(X)$, namely $Z \mapsto (Z - Z_0)$.  See
\cite[\S2]{hubbard:monodromy} for details.

From this perspective, the previous identification $P(X) \to Q(X)$
using the Schwarzian of the developing map\index{developing map} on $\H$ is simply $Z
\mapsto (Z - \sigma_0(X))$, that is, it is the Schwarzian relative to
the standard Fuchsian structure.  Complex-analytically, this is not
the most natural way to choose a basepoint in each fiber, though this
will be remedied below (\S\ref{sec:quasi-fuchsian-sections}).

The realization of $P(X)$ as an affine space modeled on a vector space
of differential forms can also be understood in terms of \v{C}ech
cochains on $X$ with a fixed coboundary
\cite[\S3]{gunning:special-coordinate}, or in terms of connections on
a principal $\PSL_2(\C)$-bundle of \emph{projective frames}
\cite[\S2.2]{anderson:projective} (and the related notions of the
\emph{graph} of a projective structure
\cite[\S2]{goldman:geometric-structures} and of \emph{$\sl_2$-opers}
\cite[\S8.2]{frenkel-ben-zvi}).

\subsection{Schwarzian parameterization of $\P(S)$}

\boldpoint{Identification of bundles.}
There is a complex vector bundle $\Q(S) \to \T(S)$ over Teichm\"uller
space whose total space consists of pairs $(X,\phi)$, where $X \in
\T(S)$ and $\phi \in Q(X)$.  In Teichm\"uller theory, this bundle is
identified with the holomorphic cotangent bundle of Teichm\"uller
space (see e.g.~\cite{imayoshi-taniguchi} \cite{hubbard:book}).  Since
Teichm\"uller space is diffeomorphic to $\R^{6g-6}$, the bundle
$\Q(S)$ is diffeomorphic to $\R^{12g-12}$.

Using a section $\sigma : \T(S) \to \P(S)$ to provide basepoints for
the fibers, we can form a bijective \emph{Schwarzian parameterization}
$$
\xymatrix@R-8mm@C-5mm{
\P(S) \ar[r] & \Q(S)\\
Z \ar@{|->}[r] & (\pi(Z), Z - \sigma(\pi(Z)))
}
$$ which is compatible with the maps of these spaces to $\T(S)$.  This
correspondence identifies the zero section of $\Q(S)$ with the section
$\sigma$ of $\P(S)$.  A different section $\sigma$ will result in a
parameterization that differs by a translation in each fiber.

\boldpoint{Compatibility.}  The topology on $\P(S)$ defined using
development-holonomy pairs\index{development-holonomy pair} is
compatible with the topology of $\Q(S)$, in that the bijection induced
by any continuous section $\sigma : \T(S) \to \P(S)$ is a
homeomorphism.  Continuity in one direction is elementary complex
analysis---uniformly close holomorphic developing
maps\index{developing map} have uniformly
close derivatives (on a smaller compact set), and therefore uniformly
close Schwarzian derivatives, making $\P(S) \to \Q(S)$ continuous.  On
the other hand, continuity of $\Q(S) \to \P(S)$ follows from
continuous dependence of solutions to the ODE
\eqref{eqn:schwarzian-ode} on its parameter $\phi$.

\boldpoint{Holomorphic structure.}
\label{sec:holomorphic}
The bundle $\Q(S)$ is a complex manifold, and a holomorphic vector
bundle over $\T(S)$.  The Schwarzian parameterization given by a
section $\sigma : \T(S) \to \P(S)$ transports these structures to
$\P(S)$.  However, two sections $\sigma_1$ and $\sigma_2$ induce the
same complex structure on $\P(S)$ if and only if $(\sigma_1 -
\sigma_2)$ is a holomorphic section of $\Q(S)$.

There is also a natural complex structure on $\P(S)$ that is defined
without reference to its parameterization by $\Q(S)$: The tangent
space $T_Z\P(S)$ can be identified with the cohomology group $H^1(Z,
\V_{proj})$, where $\V_{proj}$ is the sheaf of projective vector
fields over $Z$, i.e.~vector fields that in a local projective
coordinate are restrictions of infinitesimal M\"obius transformations.
This cohomology group is complex vector space,
which gives an integrable almost complex structure $J : T_Z\P(S) \to
T_Z \P(S)$.  (Compare the construction of \cite[Prop.~1,2]{hubbard:monodromy}.)

\boldpoint{Quasi-Fuchsian sections.}
\label{sec:quasi-fuchsian-sections}
Using deformations of Kleinian surface groups, we can construct a
class of sections of $\P(S)$ that transport the complex structure of
$\Q(S)$ to the natural complex structure on $\P(S)$.  Given $X,Y \in
\T(S)$, let $Q(X,Y)$ denote the quasi-Fuchsian group (equipped with an
isomorphism $\pi_1(S) \simeq Q(X,Y)$) that simultaneously uniformizes
$X$ and $Y$ (see e.g. \cite[Ch.~6]{imayoshi-taniguchi}).  This means that
$Q(X,Y)$ has domain of discontinuity $\Omega_+ \sqcup \Omega_-$ with
marked quotient Riemann surfaces
\begin{equation*}
\Omega_+ / Q(X,Y) \simeq X \;\;\;\; \Omega_- / Q(X,Y) \simeq \bar{Y}
\end{equation*}
where $\bar{Y}$ is the complex conjugate Riemann surface of $Y$,
which appears in the quotient because the induced orientation on the
marked surface $\Omega_- / Q(X,Y)$ is opposite that of $S$.

As a quotient of a domain by a Kleinian group, the surface $\Omega_+ /
Q(X,Y)$ also has a natural projective structure, which we denote by
$\Sigma_Y(X)$.  By definition, the underlying Riemann surface of
$\Sigma_Y(X)$ is $X$, so for any fixed $Y \in T(S)$ this defines a
\emph{quasi-Fuchsian section}\index{quasi-Fuchsian
  section}\index{projective structure!quasi-Fuchsian}
\begin{equation*}
\Sigma_Y : \T(S) \to \P(S).
\end{equation*}

These quasi-Fuchsian sections induce the natural complex structure on
$\P(S)$.  We sketch two ways to see this: First, Hubbard uses a
cohomology computation to show that a section induces the canonical
complex structure if and only if it can be represented by a
\emph{relative projective structure}\index{relative projective
  structure}\index{projective structure!relative} on the universal curve over
$\T(S)$ \cite[Prop.~1,2]{hubbard:monodromy}.  The quasi-Fuchsian
groups provide such a structure due to the analytic dependence of the
solution of the Beltrami equation on its parameters
\cite{ahlfors-bers}, and the associated construction of the \emph{Bers
fiber space}\index{Bers fiber space} \cite{bers:fiber-spaces}.

Alternatively, one can show (as in the respective computations of
Hubbard \cite{hubbard:monodromy} and Earle \cite{earle:variation})
that both the canonical complex structure on $\P(S)$ and the complex
structure coming from a quasi-Fuchsian section make the holonomy map
(discussed in \S5) a local biholomorphism, and therefore they are
holomorphically equivalent.

\boldpoint{Norms.}\index{norm!of a quadratic differential}
A norm on the vector space $Q(X)$ induces a natural measure of the
``complexity'' of a projective structure on $X$ (relative to the
standard Fuchsian structure), or of the difference between two
projective structures.  There are several natural choices for such a
norm.

The hyperbolic $L^\infty$ norm $\|\phi\|_\infty$ is the supremum of
the function $|\phi| / \rho^2$, where $\rho^2$ is the area element of
the hyperbolic metric on $X$.  Lifting $\phi$ to the universal cover
and identifying $\Tilde{X} \to \Delta$, we have
$$ \| \phi \|_\infty = \| \Tilde{\phi} \|_\infty = \frac{1}{4} \sup_{z
  \in \Delta} |\Tilde{\phi}(z)| (1 - |z|^2)^2. $$

By Nehari's theorem, a holomorphic immersion $f : \Delta \to \CP^1$
satisfying $\|S(f)\|_\infty \leq \frac{1}{2}$ is injective, while any
injective map satisfies $\| S(f) \|_\infty \leq \frac{3}{2}$ (see
\cite{nehari:schwarzian}, also \cite{pommerenke:univalent}
\cite{lehto:univalent}).  More generally, the norm $\|S(f)\|_\infty$
gives a coarse estimate of the size of hyperbolic balls in $\Delta$ on
which $f$ is univalent \cite[\S3]{kra:deformations}
\cite[Lem.~5.1]{kra-maskit:remarks}.  Thus, when applied to projective
structures, the $L^\infty$ norm reflects the geometry and valence of
the developing map\index{developing map}.

In Teichm\"uller theory, it is more common to use the $L^1$ norm
$\|\phi\|_1$, which is the area of the surface $X$ with respect to the
singular Euclidean metric $|\phi|$.  This norm is conformally natural,
since it does not depend on the choice of a Riemannian metric on $X$.
However, the intrinsic meaning of the $L^1$ norm of the Schwarzian
derivative is less clear.

More generally, given any background Riemannian metric on $X$
compatible with its conformal structure, there is an associated $L^p$
norm on $Q(X)$.  These norms, with $p \in (1,\infty)$ and especially
$p=2$, can be used to apply PDE estimates to the study of projective
structures, as discussed in \S\ref{sec:schwarzian-proof} below.

Note that while any two norms on the finite-dimensional vector space
$Q(X)$ are bilipschitz equivalent, the bilipschitz constants between
the $L^\infty$, $L^1$, and hyperbolic $L^p$ norms on $Q(X)$ diverge as
$X \to \infty$ in Teichm\"uller space.

\section{The Grafting Parameterization}\label{sec:grafting}

\subsection{Definition of grafting}

Grafting is a geometric operation that can be used to build an
arbitrary projective structure by gluing together simple pieces.  We
start by defining grafting in a restricted setting, and then work toward
the general definition.

\begin{figure}
\begin{center}
\includegraphics[width=\textwidth]{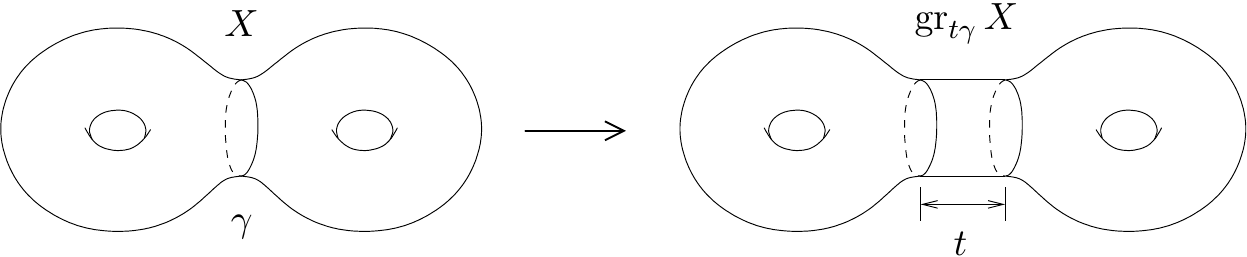}
\end{center}
\caption{Grafting along a simple closed curve.\label{fig:basicgraft}}
\end{figure}

\boldpoint{Grafting simple geodesics.}
Equip a Riemann surface $X \in \T(S)$ with its hyperbolic metric.  Let
$\gamma$ be a simple closed hyperbolic geodesic on $X$.  The basic
grafting construction replaces $\gamma$ with the cylinder $\gamma
\times [0,t]$ to obtain a new surface $\gr_{t \gamma} X$, the
\emph{grafting of $X$ by $t \gamma$}\index{grafting}, as shown in Figure
\ref{fig:basicgraft}.  The natural metric on this surface is partially
hyperbolic (on $X-\gamma$) and partially Euclidean (on the cylinder),
and underlying this metric is a well-defined conformal structure on
$\gr_{t \gamma} X$.

Let $\scc$ denote the set of free homotopy classes of homotopically
nontrivial simple closed curves on $S$.  Then $\scc$ is canonically
identified with the set of simple closed geodesics for any hyperbolic
structure on $S$, and we can regard grafting as a map
\begin{equation*}
\gr : \scc \times \R^+ \times \T(S) \to \T(S).
\end{equation*}
When it is important to distinguish this construction from the
projective version defined below, we will call this \emph{conformal
  grafting}\index{grafting!conformal}\index{grafting}, since the result is a
conformal structure.

\boldpoint{Projective grafting.}
\label{sec:projective-grafting}
The Riemann surface $X$ has a standard Fuchsian projective structure
in which the holonomy of a simple closed geodesic $\gamma$ is
conjugate to $z \mapsto e^\ell z$, where $\ell = \ell(\gamma,X)$ is
the hyperbolic length of $\gamma$.  

For any $t < 2\pi$, let $\Tilde{A_t}$ denote a sector of angle $t$ in
the complex plane, with its vertex at $0$.  The quotient $A_t =
\Tilde{A_t} / \langle z \mapsto e^\ell z \rangle$ is an annulus
equipped with a projective structure, which as a Riemann surface is
isomorphic to the Euclidean product $\gamma \times [0,t]$.

\begin{figure}
\begin{center}
\includegraphics[width=\textwidth]{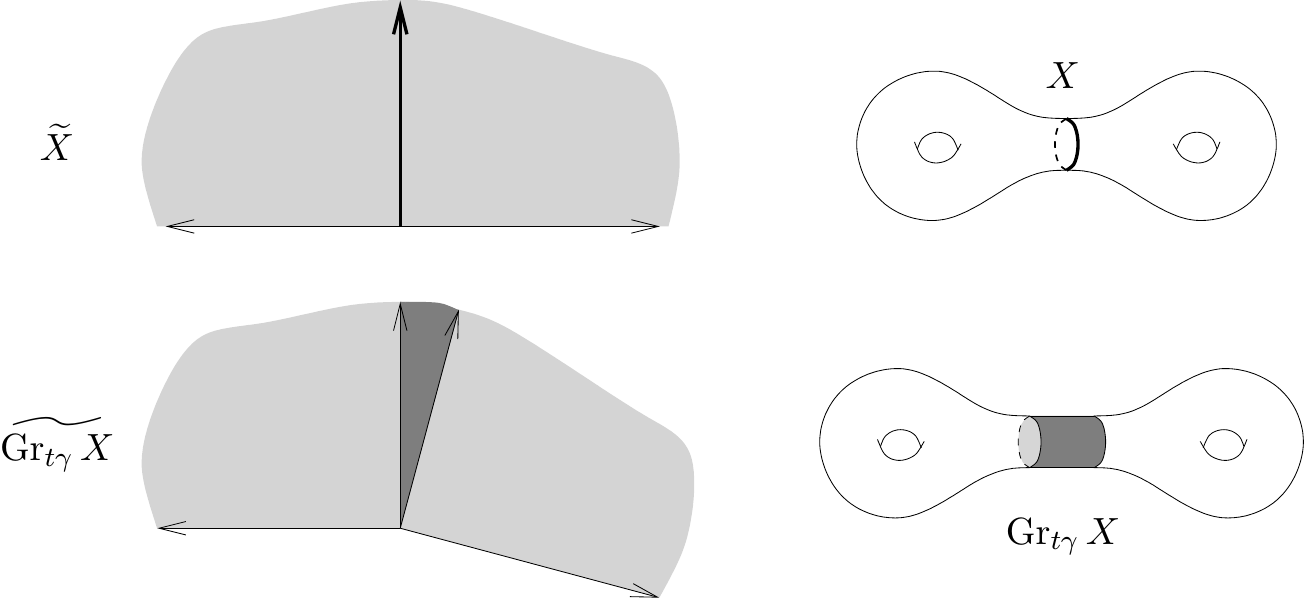}
\end{center}
\caption{Projective grafting\index{grafting!projective}\index{grafting}: Gluing a cylinder into the surface along
  a geodesic corresponds to inserting a sector or lune into each lift
  of the geodesic.  Only one lift is shown here, but the gluing
  construction is repeated equivariantly in $\widetilde{\Gr_{t
  \gamma}X}$.}
\label{fig:sector}
\end{figure}

There is a natural projective structure on the grafted surface $\gr_{t
  \gamma} X$ that is obtained by gluing the standard Fuchsian
projective structure of $X$ to $A_t$; these structures are compatible
due to the matching holonomy around the gluing curves.  In the
universal cover of $X$, this corresponds to inserting a copy of
$\Tilde{A_t}$ in place of each lift of $\gamma$ (see Figure
\ref{fig:sector}), applying M\"obius transformations to $\Tilde{A_t}$
and the complementary regions of $\gamma$ in $\Tilde{X}$ (which are
bounded by circular arcs) so that they fit together.  For sufficiently
small $t$, this produces a Jordan domain in $\CP^1$ that is the image
of the developing map\index{developing map}, while for large $t$ the
developing image is all of $\CP^1$.  We denote the resulting
projective structure by $\Gr_{t \gamma} X$.

Applying a generic M\"obius transformation to the sector $\Tilde{A_t}$
will map it to a $t$-\emph{lune}, the intersection of two round disks
with interior angle $t$.  Thus the projective structure $\Gr_{t
\gamma} X$ corresponds to a decomposition of its universal cover into
$t$-lunes and regions bounded by circular arcs.

The restriction to small values of $t$ in this construction is
not necessary; for $t > 2 \pi$ we simply interpret $\Tilde{A_t}$ as a
``sector'' that wraps around the punctured plane $\C^*$ some number of
times.  Alternatively, we could define $A_t$ for $t \geq 2\pi$ by
gluing $n$ copies of $A_{t/n}$ end-to-end, for a sufficiently large $n
\in \N$.

Therefore we have a \emph{projective grafting
  map}\index{grafting!projective}\index{grafting},
\begin{equation*}
\Gr : \scc \times \R^+ \times \T(S) \to \P(S)
\end{equation*}
which is a lift of grafting through the forgetful map $\pi: \P(S)
\to \T(S)$, i.e.~$\pi \circ \Gr = \gr$.

\boldpoint{Variations on simple grafting.}  Grafting along a simple
geodesic with weight $t = 2 \pi$ was originally used by Maskit
\cite{maskit:grafting}, Hejhal \cite{hejhal:monodromy}, and
Sullivan-Thurston \cite{sullivan-thurston} to construct examples of
exotic Fuchsian projective structures (discussed in
\S\ref{sec:fuchsian-holonomy} below).  Grafting\index{grafting} with weight $2\pi$ is
special because it does not change the holonomy representation of the
Fuchsian projective structure (see \S\ref{sec:holonomy}).

It is possible to extend this holonomy-preserving grafting operation
to certain simple closed curves which are not geodesic, and to projective
structures that are not standard Fuchsian (see
\cite[Ch.~7]{kapovich:book}); this generalization has been important
to some applications in Kleinian groups and hyperbolic geometry
(e.g.~\cite{bromberg:degenerate}\cite[\S5]{brock-bromberg}), and it
will appear again in our description of quasi-Fuchsian projective
structures (\S\ref{sec:quasi-fuchsian-holonomy}).  However, our main
focus in this chapter is a different extension of grafting, defined by Thurston,
which leads to a geometric model for the entire moduli space $\P(S)$.

\boldpoint{Extension to laminations.}
Projective grafting is compatible with the natural completion of $\R^+
\times \scc$ to the space $\ML(S)$ of \emph{measured
  laminations}\index{measured lamination}.  An
element $\lambda \in \ML(S)$ is realized on a hyperbolic surface $X
\in \T(S)$ as a foliation of a closed subset of $X$ by complete,
simple hyperbolic geodesics (some of which may be closed), equipped
with a transverse measure of full support.  A piecewise linear
coordinate atlas for $\ML(S)$ is obtained by integrating transverse
measures over closed curves, making $\ML(S)$ into a $PL$-manifold
homeomorphic to $\R^{6g-6}$. See \cite[Ch.~8-9]{thurston:notes}
\cite{casson-bleiler} \cite[Ch.~3]{penner-harer} for detailed
discussion of measured laminations.

There is continuous extension $\Gr : \ML(S) \times \T(S) \to \P(S)$ of
projective grafting\index{grafting}, which is uniquely determined by the simple
grafting construction because weighted simple curves are dense in
$\ML(S)$.  Similarly, there is an extension of the grafting map $\gr :
\ML(S) \times \T(S) \to \T(S)$ defined by $\gr = \pi \circ \Gr$.
These extensions were defined by Thurston [unpublished], and are
discussed in detail in \cite{kamishima-tan:grafting}.

For a lamination $\lambda \in \ML(S)$ that is supported on a finite
set of disjoint simple closed curves, i.e.~$\lambda = \sum_{i=1}^n t_i
\gamma_i$, the grafting $\gr_\lambda X$ defined by this extension
procedure agrees with the obvious generalization of grafting along
simple closed curves, wherein the geodesics $\gamma_1, \ldots,
\gamma_n$ are simultaneously replaced with cylinders.

For a general measured lamination $\lambda \in \ML(S)$, one can think
of $\gr_{\lambda} X$ as a Riemann surface obtained from $X$ by
thickening the leaves of the lamination $\lambda$ in a manner dictated
by the transverse measure.  This intuition is made precise by the
definition of a canonical stratification of $\gr_{\lambda}X$ in the
next section.

\subsection{Thurston's Theorem}
\label{sec:thurston-theorem}
Projective grafting\index{grafting} is a universal
construction---every projective structure can be obtained from it, and
in exactly one way:

\begin{theorem}[{Thurston [unpublished]}]
\label{thm:thurston}
The projective grafting map $\Gr : \ML(S) \times \T(S) \to \P(S)$ is a
homeomorphism.
\end{theorem}

The proof of Theorem \ref{thm:thurston} proceeds by explicitly
constructing the inverse map $\Gr^{-1}$ using complex projective and
hyperbolic geometry.  We will now sketch this construction; details
can be found in \cite{kamishima-tan:grafting}.

\boldpoint{The embedded case.}  First suppose that $Z \in \P(S)$ is a
projective surface whose developing map\index{developing map} is an
embedding (an \emph{embedded projective structure}%
\index{embedded projective structure}%
\index{projective structure!embedded}).  The image of the developing
map is a domain $\Omega \subset \CP^1$ invariant under the action of
$\pi_1(S)$ by the holonomy representation $\rho$.  In this case, we
will describe the inverse of projective grafting in terms of convex
hulls in hyperbolic space.  See \cite{epstein-marden:convex-hulls} for
details on these hyperbolic constructions.

Considering $\CP^1$ as the ideal boundary of hyperbolic space $\H^3$,
let $\Pl(Z)$ denote the boundary of the hyperbolic convex
hull\index{convex hull} of $(\CP^1 - \Omega)$.  Then $\Pl(Z)$ is a
\emph{convex pleated plane}\index{convex pleated plane}\index{pleated
plane!convex} in $\H^3$ invariant under the action of $\pi_1(S)$ by
isometries.

When equipped with the path metric, the pleated plane $\Pl(Z)$ is
isometric to $\H^2$, and by this isometry, the action of $\pi_1(S)$ on
$\H^3$ corresponds to a discontinuous cocompact action on $\H^2$.  Let
$Y \in \T(S)$ denote the marked quotient surface.

The pleated plane $\Pl(Z)$ consists of totally geodesic pieces
(\emph{plaques}\index{plaque} or \emph{facets}\index{facet}) meeting
along geodesic \emph{bending lines}\index{bending line}.  Applying the
isometry $\Pl(Z) \simeq \H^2$ to the union of the bending lines yields
a geodesic lamination, which has a natural transverse measure
recording the amount of bending of $\Pl(Z)$.  The lamination and
measure are $\pi_1(S)$-invariant, and therefore descend to the
quotient, defining an element $\lambda \in \ML(S)$.

Thus, starting from an embedded projective structure $Z$, we obtain a
hyperbolic structure $Y$ and a measured lamination $\lambda$.  To show
that we have inverted the projective grafting\index{grafting} map, we must check that
$\Gr_\lambda Y = Z$.

\boldpoint{Nearest-point projection.}
\label{sec:nearest-point}
There is a \emph{nearest-point projection map}\index{nearest-point
  projection} $\kappa : \Omega \to \Pl(Z)$ that sends $z \in \Omega$
to the first point on $\Pl(Z)$ that is touched by an expanding family
of horoballs in $\H^3$ based at $z$.  Convexity of $\Pl(Z)$ ensures
that this point is well-defined.  In fact, from each $z \in \Omega$ we
obtain not just a nearest point on $\Pl(Z)$, but also a \emph{support
  plane}\index{support plane} $H_z$ which contains $\kappa(z)$ and whose normal vector at
that point defines a geodesic ray with ideal endpoint $z$.  This gives
a map $\hat{\kappa} : \Omega \to \planes$, where $\planes$ is the
space of planes in $\H^3$ (the \emph{de Sitter space}\index{de Sitter space}).

The \emph{canonical stratification}\index{canonical stratification}%
\index{stratification!canonical}\index{strata} of $\Omega$ is the
decomposition into fibers of the map $\hat{\kappa}$.  Strata are of
two types:
\begin{itemize}
\item $1$-dimensional strata---circular arcs that map homeomorphically
  by $\kappa$ onto bending lines of $\Pl(Z)$, and
\item $2$-dimensional strata---regions with nonempty interior bounded
by circular arcs which map homeomorphically by $\kappa$ to the totally
geodesic pieces of $\Pl(Z)$.
\end{itemize}

\begin{figure}
\begin{center}
\subfigure[Developed image]{\begin{minipage}{0.48\textwidth}\centering\hspace{-0.8mm}\includegraphics[width=5.3cm]{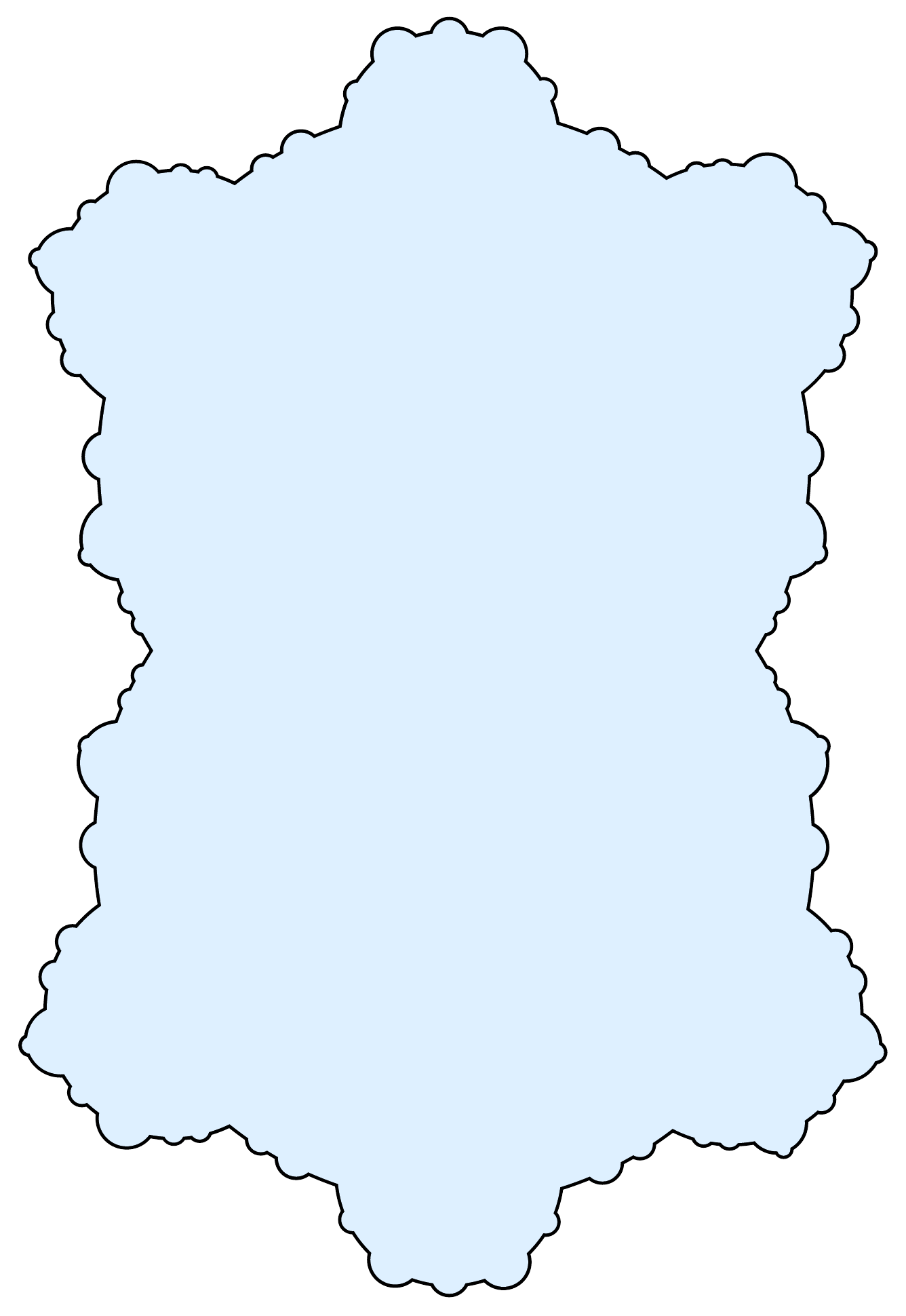}\vspace{2mm}\end{minipage}}
\subfigure[Maximal disks]{\begin{minipage}{0.48\textwidth}\centering\hspace{-0.8mm}\includegraphics[width=5.3cm]{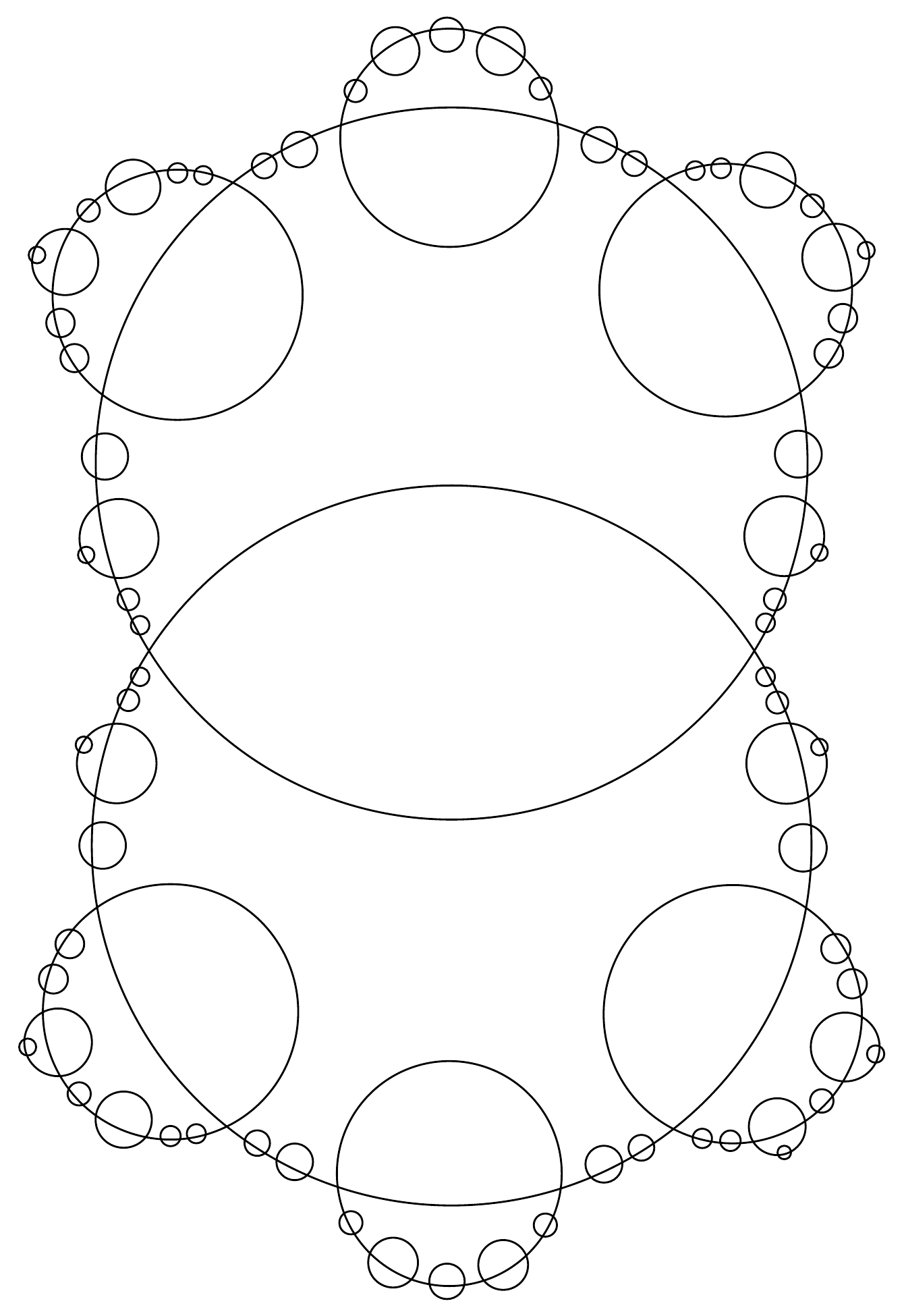}\vspace{2mm}\end{minipage}}\\
\subfigure[Lunes and ideal polygons]{\begin{minipage}[b]{0.48\textwidth}\centering\hspace{-0.8mm}\includegraphics[width=5.3cm]{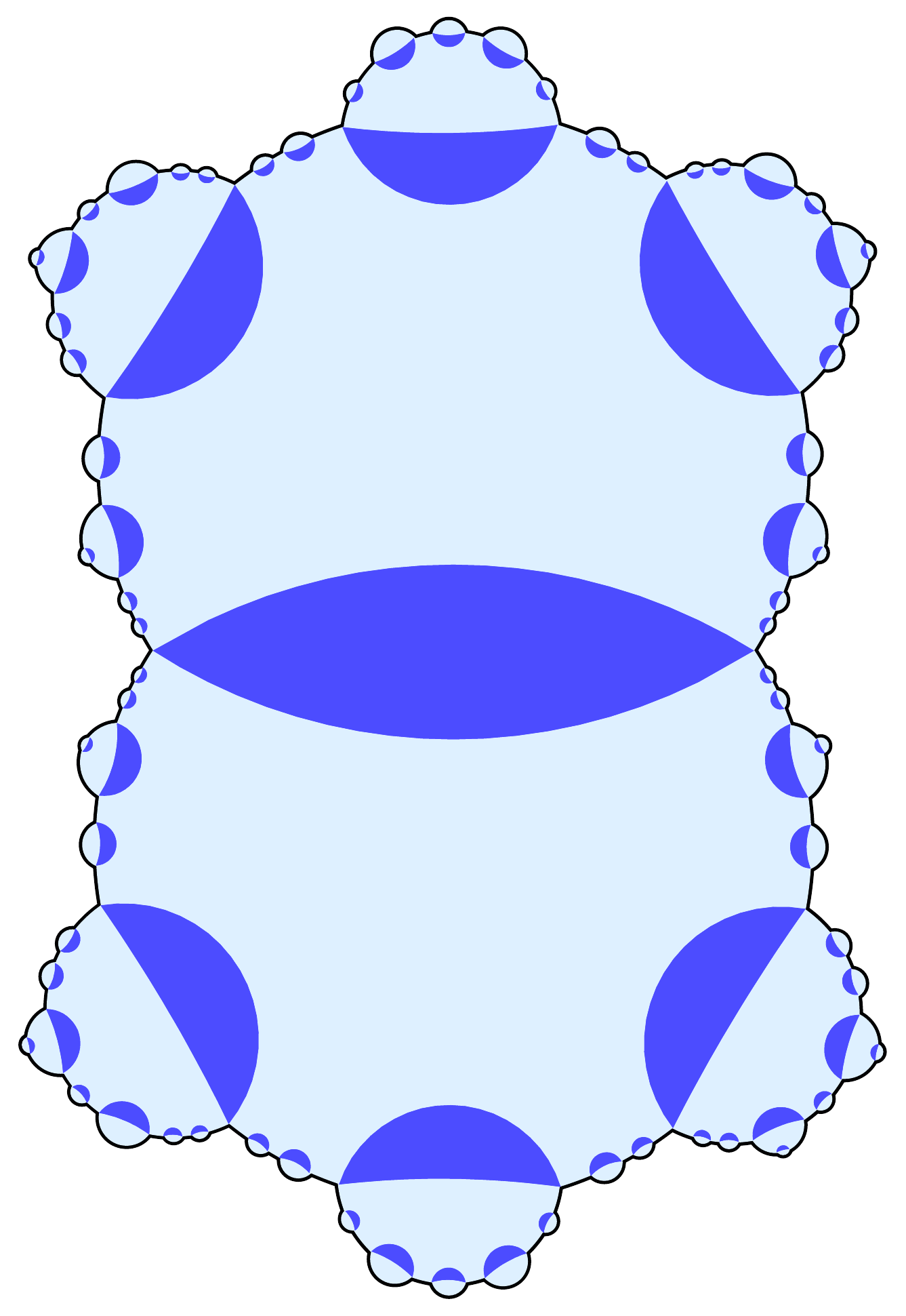}\vspace{2mm}\end{minipage}}
\subfigure[Pleated plane in $\H^3$ (Klein model)]{\begin{minipage}[b]{0.48\textwidth}\centering\hspace{-0.8mm}\includegraphics[width=5.3cm]{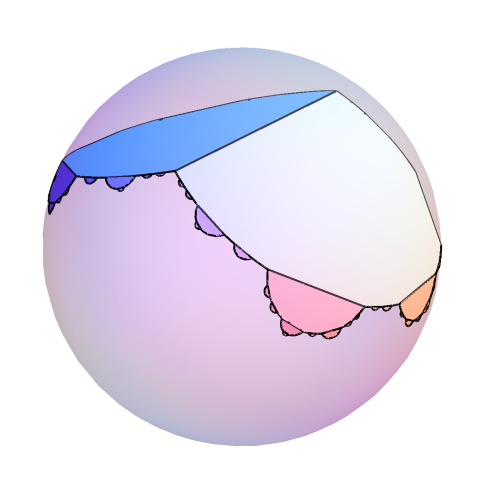}\vspace{2mm}\end{minipage}}
\end{center}
\caption{Four views of a projective structure lifted to the universal
  cover of a surface.  The example shown here is an approximation of
  an embedded structure on a surface of genus $2$, obtained by
  grafting along a separating simple closed curve.  The approximation
  includes only a few of the maximal disks.\label{fig:bend}}
\end{figure}

If $\lambda$ is supported on a single closed geodesic (or on a finite
union of them), the $1$-dimensional strata and the boundary geodesics
of the $2$-dimensional strata in $\Omega \simeq \Tilde{Z}$ fill out a
collection of lunes, and the interiors of the $2$-dimensional strata
correspond by $\kappa$ to the complementary regions of the lift of
$\lambda$, realized geodesically on $Y$, to $\Tilde{Y} \simeq \H^2$.
See Figure \ref{fig:bend} for an example of this type.  This is the
arrangement of lunes and circular polygons giving the projective
structure of $\Gr_\lambda Y$, an so $Z = \Gr_\lambda Y$. A limiting
argument shows the same holds for general $\lambda$.

\boldpoint{The general case.}  The key to inverting the projective
grafting\index{grafting} in the embedded case is the construction of the convex
pleated plane $\Pl(Z)$.  For general $Z \in \P(S)$, this is replaced
by a \emph{locally convex pleated plane}\index{locally convex pleated
  plane}\index{pleated plane!locally convex} defined using the
projective geometry of $Z$ itself, rather than its developed image.

Let $f : \Tilde{Z} \to \CP^1$ be the developing map\index{developing map} of $Z$.  A
\emph{round disk}\index{round disk} in $\Tilde{Z}$ is an open subset $U$ such that $f$
is injective on $U$ and $f(U)$ is an open disk in $\CP^1$.  The round
disks in $\Tilde{Z}$ are partially ordered with respect to inclusion.
A maximal element for this ordering is a \emph{maximal round
  disk}\index{maximal round disk}.

Each maximal round disk $U$ in $\Tilde{Z}$ corresponds to a disk in
$\CP^1$, and thus to an oriented plane $H_U$ in $\H^3$.  Allowing $U$
to vary over all maximal round disks in $\Tilde{Z}$ gives a family of
oriented planes, and the envelope of this family is a locally convex
pleated plane $\Pl(Z)$.  

The rest of the convex hull construction generalizes as follows: The
intrinsic geometry of $\Pl(Z)$ is hyperbolic, with quotient $Y$, and
the bending of $\Pl(Z)$ is recorded by a measured lamination
$\lambda$.  In place of the nearest-point projection and support
planes, we have a \emph{collapsing map}\index{collapsing map} $\kappa
: \Tilde{Z} \to \Pl(Z)$ and a \emph{co-collapsing
  map}\index{co-collapsing map} $\hat{\kappa} : \Tilde{Z} \to \planes$
(see also \cite[\S2,\S7]{dumas:antipodal}).  The fibers of
$\hat{\kappa}$ induce a canonical stratification\index{canonical stratification}%
 of $\Tilde{Z}$, and separating the $1$- and
$2$-dimensional strata\index{strata} describes $Z$ as the projective grafting
$\Gr_\lambda Y$.

Note that the canonical stratification of $\Tilde{Z}$ is
$\pi_1(S)$-invariant, and therefore we have a corresponding
decomposition of $Z$ into $1$- and $2$-dimensional pieces.  We will
also refer to this as the canonical stratification.  Similarly, the
collapsing map descends to a map $\kappa : Z \to Y$ between quotient
surfaces, which sends the union of $1$-dimensional strata and boundary
geodesics of $2$-dimensional strata onto the bending lamination
$\lambda \subset Y$.

The canonical stratification for complex projective structures is
discussed further in \cite[\S1.2]{kamishima-tan:grafting}, where it is
also generalized to $n$-manifolds equipped with \emph{flat
  conformal structure}\index{flat conformal structure} (see also
\cite{kulkarni-pinkall:metric} \cite{scannell:desitter}).

\begin{figure}
\begin{center}
\includegraphics{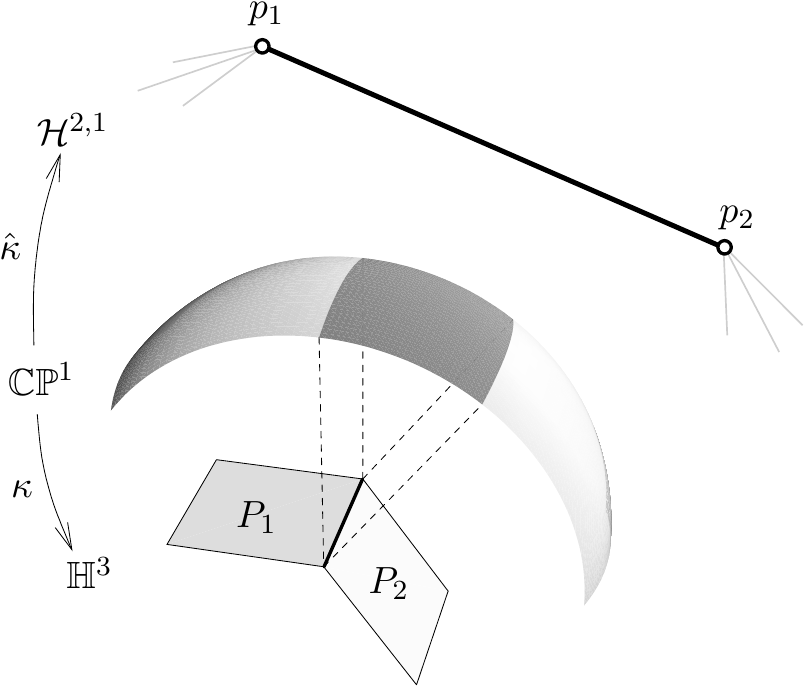}
\end{center}
\caption{A lune between two maximal disks collapses to a
  bending line between two planes ($P_1$,$P_2$), and co-collapses to
  an interval between two points ($p_1$,$p_2$).\label{fig:cocollapse}}
\end{figure}

\boldpoint{Dual trees.}
\label{sec:dual-tree}
When grafting\index{grafting} along a simple closed curve $\gamma$ with weight $t$,
each bending line of the associated pleated plane in $\H^3$ has a
one-parameter family of support planes (see Figure
\ref{fig:cocollapse}).  These give an interval in the image of
$\hat{\kappa}$, and the angle between support planes gives a metric on
this interval, making it isometric to $[0,t] \subset \R$.
Alternatively, this metric could be defined as the restriction of the
Lorentzian metric of $\planes$, where the restriction is positive
definite because any pair of support planes of a given bending
line\index{bending line} intersect (see \cite[\S5]{scannell:desitter}
\cite[\S3,\S6.5]{kulkarni-pinkall:metric}).

The intervals corresponding to different bending lines meet at
vertices corresponding to support planes of flat pieces.  This gives
$\hat{\kappa}$ the structure of a metric tree, the \emph{dual
  tree}\index{dual tree} of
the weighted curve $t \gamma$, denoted $T_{t\gamma}$.  As this
notation suggests, this tree depends only on $t\gamma$ (through the
bending lines, their bending angles, and the adjacency relationship
between bending lines and flat pieces) and not on the quotient
hyperbolic structure of the pleated plane.  The equivariance of the
pleated plane with respect to $\pi_1(S)$ determines an isometric
action of $\pi_1(S)$ on $T_{t\gamma}$.

For a general grafting lamination $\lambda \in \ML(S)$, the image of
$\hat{\kappa}$ has the structure of a \emph{$\R$-tree}\index{R-tree@$\R$-tree} (see
\cite[Ch.~9]{scannell:thesis}
\cite[\S6,\S11]{kulkarni-pinkall:metric}), a geodesic metric space in
which each pair of points is joined by a unique geodesic which is
isometric to an interval in $\R$
\cite[Ch.~2]{morgan-shalen:valuations-trees}.  This \emph{dual
$\R$-tree of $\lambda$}, denoted $T_\lambda$, is also equipped with an
isometric action of $\pi_1(S)$.

\subsection{The Thurston metric}
\label{sec:thurston-metric}
We have seen that when grafting\index{grafting} along a simple closed curve,
the resulting projective surface $\Gr_{t \gamma} X$ has a natural
conformal metric that combines the hyperbolic structure of $X$ and the
Euclidean structure of the cylinder.  This is the \emph{Thurston
  metric}\index{Thurston metric} (or \emph{projective
  metric}\index{projective metric}) on the projective surface.

This definition can be extended to arbitrary projective surfaces by
taking limits of the metrics obtained from an appropriate sequence of
simple closed curves; however, we will prefer an intrinsic description
of the metric based on complex projective geometry.

\boldpoint{Kobayashi construction.}  The Kobayashi
metric\index{Kobayashi metric} on a complex manifold is defined by a
norm on each tangent space, where the length of a vector $v$ is the
infimum of lengths given to it by holomorphically immersed disks (each
of which is equipped with its hyperbolic metric).  For a surface $Z$
with a projective structure, there is a variant of the Kobayashi
metric in which one minimizes length over the smaller class of
\emph{projectively immersed disks}\index{projectively immersed disk},
that is, immersions $\Delta \to \Tilde{Z}$ that are locally M\"obius with
respect to the projective structure on $\Delta$ as a subset of
$\CP^1$.  The resulting ``projective Kobayashi metric'' is the
\emph{Thurston metric} of $Z$ \cite[\S2.1]{tanigawa:grafting}.

\boldpoint{Relation to grafting.}  This intrinsic definition of the
Thurston metric is related to grafting as follows: for each $z \in
\Tilde{Z}$, there is a unique maximal round disk $U \subset \Tilde{Z}$
such that the (lifted) Thurston metric at $z$ agrees with the
hyperbolic metric on $U$.  Furthermore, the set of points in
$\Tilde{Z}$ that correspond to a given maximal disk $U$ is a stratum
in the canonical stratification of $Z$.  Thus the Thurston metric is
built from $2$-dimensional hyperbolic regions with geodesic boundary
and $1$-dimensional geodesic strata\index{strata}.  For $Z =
\Gr_\lambda X$, the union of the hyperbolic strata covers a subset of
$Z$ isometric to $(X - \lambda)$, where $\lambda$ is realized
geodesically on the hyperbolic surface $X$.  When $\lambda = t \gamma$
is supported on a simple closed curve, the $1$-dimensional strata
sweep out Euclidean strips in $\Tilde{\Gr_{t \gamma} X}$, which cover
a Euclidean cylinder in $\Gr_{t \gamma} X$, recovering the synthetic
description of the Thurston metric in this case.

\boldpoint{Conformal metrics and regularity.}  The Thurston metric on
$Z$ is a nondegenerate Riemannian metric compatible with the
underlying complex structure $\pi(Z)$, i.e.~it is a \emph{conformal
  metric}\index{conformal metric} on the Riemann surface.  In local
complex coordinates, the line element of such a metric has the form
$\rho(z) |dz|$, where $\rho(z)$ is the real-valued \emph{density
  function}%
\index{density function}.
In the case of simple grafting\index{grafting}, the density function
of the Thurston metric is smooth on the hyperbolic and Euclidean
pieces, but it is only $C^1$ on the interface between them.  (The
discontinuity in its second derivative is necessary since the
curvature changes along the interface.)  In general, the Thurston
metric of a projective surface is $C^{1,1}$, meaning that its density
function has Lipschitz derivatives, with Lipschitz constant locally
bounded on $\P(S)$ \cite{kulkarni-pinkall:metric}.

\boldpoint{Variation of metrics.}  The Thurston metric is a continuous
function of the projective structure $Z \in \P(S)$ with respect to the
topology of locally uniform convergence of density functions: For a
sequence $Z_n \to Z \in \P(S)$, the Lipschitz bound on the derivatives
of the Thurston metrics shows that uniform convergence follows from
pointwise convergence, which in turn follows from the locally uniform
convergence of the developing maps\index{developing map} $f_n : \Delta
\to \CP^1$ (or from the continuous variation of the associated locally
convex pleated surfaces).

\boldpoint{Area.}
A conformal metric on a Riemann surface with density
function\index{density function} $\rho$
induces an area measure by integration of $\rho^2 = \rho(z)^2 |dz|^2$.

The total area of $\Gr_\lambda X$ with respect to the Thurston metric
is $4 \pi (g-1) + \ell(\lambda,X)$, where $\ell(\lambda,X)$ is the
length of the measured lamination $\lambda$ with respect to the
hyperbolic metric of $X$.  The two terms correspond to the two types
of strata\index{strata}: The union of the $2$-dimensional strata has area $4 \pi
(g-1)$, because it is isometric to the complement of a geodesic
lamination (a null set) in the hyperbolic surface $X$.  The union of
the $1$-dimensional strata has area $\ell(\lambda,X)$, which is the
continuous extension to $\ML(S)$ of the function $t \ell(\gamma,X)$
giving the area of the Euclidean cylinder $\gamma \times [0,t]$ in the
case of simple grafting.

\boldpoint{Curvature.}
The Gaussian curvature $K$ and curvature $2$-form $\Omega$ of a smooth
conformal metric are related to its density function\index{density function} $\rho$ by
\begin{equation}
\label{eqn:curvature}
\begin{split}
K &= - \frac{1}{\rho^2} \Delta \log \rho\\
\Omega &= K \rho^2 = - \Delta \log \rho.
\end{split}
\end{equation}
In particular, such a metric has nonpositive Gaussian curvature if and
only if $\log \rho$ is a subharmonic function.

The Thurston metric is not smooth everywhere, but it is nonpositively
curved (\npc)\index{nonpositively curved (NPC)}, meaning that its
geodesic triangles are thinner than triangles in Euclidean space with
the same edge lengths.  As in the smooth case, this implies that $\log
\rho$ is subharmonic, so we have a nonpositive measure $\Omega = -
\Delta \log(\rho)$ that generalizes the curvature $2$-form
\cite{reshetnyak:conformal} (see also \cite{reshetnyak:survey}
\cite{huber:zum} \cite{huber:subharmonic}
\cite{mese:between-surfaces}).  For the Thurston metric, $\Omega$ is
absolutely continuous, $\Omega = K \rho^2$ where $K$ is the
(a.e.~defined) Gaussian curvature function.  By a generalization of
the Gauss-Bonnet theorem, the total mass of $\Omega$ (which is the
integral of $K$) is $-4 \pi (g-1)$ \cite{huber:subharmonic}.

\boldpoint{Hyperbolic and Euclidean.}
Since the Gaussian curvature of the Thurston metric is $-1$ in the
interior of each $2$-dimensional stratum, and these have total area $4
\pi (g-1)$, the curvature of the Thurston metric is almost everywhere
$0$ in the union of the $1$-dimensional strata\index{strata}.  In this sense,
grafting along a general lamination can be seen as the operation of
inserting a Euclidean ``surface'' in place of a geodesic lamination,
generalizing the case of closed leaves.

\subsection{Conformal grafting maps}
\label{sec:conformal-grafting}
Having discussed the projective grafting construction and its inverse,
we turn our attention to properties of the conformal grafting
map\index{grafting!conformal}\index{grafting} $\gr : \ML(S) \times
\T(S) \to \T(S)$.

Using techniques from the theory of harmonic maps\index{harmonic map}
between surfaces (see \S\ref{sec:harmonic-maps}), Tanigawa showed that
this map is proper when either one of the coordinates is fixed:

\begin{theorem}[{Tanigawa \cite{tanigawa:grafting}}]
\label{thm:properness}
For each $\lambda \in \ML(S)$, the $\lambda$-grafting map $\gr_\lambda
: \T(S) \to \T(S)$ is a proper smooth map. For each $X \in \T(S)$, the
$X$-grafting map $\gr_{\param}X : \ML(S) \to \T(S)$ is a proper
continuous map.
\end{theorem}

Properness allows global properties of these maps to be derived from
local considerations.  For example, Scannell and Wolf showed that the
$\lambda$-grafting map is an immersion, and therefore it is a local
diffeomorphism.  Since a proper local diffeomorphism is a covering
map, this result and Theorem \ref{thm:properness} give:

\begin{theorem}[{Scannell-Wolf \cite{scannell-wolf:grafting}}]
\label{thm:scannell-wolf}
For each $\lambda \in \ML(S)$, the $\lambda$-grafting map $\gr_\lambda
: \T(S) \to \T(S)$ is a diffeomorphism.
\end{theorem}

Earlier, Tanigawa had shown that $\gr_\lambda$ is a diffeomorphism
when $\lambda \in \ML(S)$ is supported on a finite set of simple
closed curves with weights that are integral multiples of $2 \pi$
\cite{tanigawa:grafting}.  This follows from Theorem
\ref{thm:properness} because holonomy considerations (see
\S\ref{sec:holonomy}) imply that $\gr_\lambda$ is a local
diffeomorphism in this case.

In the general case, Scannell and Wolf analyze the Thurston metric and
conformal grafting\index{grafting!conformal}\index{grafting} map
through the interaction of two differential equations: The Liouville
equation, which relates a Riemannian metric to its curvature, and the
Jacobi equation, which determines the variation of a geodesic with
respect to a family of Riemannian metrics.  Analytic estimates for
these equations are used to show that a $1$-parameter family of
graftings $t \mapsto \gr_\lambda X_t$ cannot be conformally equivalent
to first order unless $\left.(d/dt)X_t\right|_{t=0} = 0$, which gives
injectivity of the derivative of $\gr_\lambda$.

As a consequence of Theorem \ref{thm:scannell-wolf}, for any $\lambda
\in \ML(S)$ the set of projective structures with grafting lamination
$\lambda$ projects homeomorphically to $\T(S)$ by the forgetful
map.  That is, the set of such projective structures forms a
smooth section $\sigma_\lambda : \T(S) \to \P(S)$ of $\pi$, which is given by
\begin{equation}
\label{eqn:sigma-lambda}
\sigma_\lambda(X) = \Gr_{\lambda}(\gr_\lambda^{-1}(X)).
\end{equation}
Note that this is compatible with our previous definition of the
standard Fuchsian structure $\sigma_0(X)$, since this is the unique
projective structure on $X$ with zero grafting lamination.  As with
Theorem \ref{thm:scannell-wolf}, in the special case of
$2\pi$-integral weighted multicurves, the existence of these smooth
sections follows from the earlier work of Tanigawa.

Fixing $X$ and varying $\lambda \in \ML(S)$, we can also use Theorem
\ref{thm:scannell-wolf} to parameterize the fiber $P(X)$; that is,
$$ \lambda \mapsto \sigma_\lambda(X)$$ gives a homeomorphism $\ML(S)
\to P(X)$ (compare \cite[\S4]{dumas:schwarzian}).  It is the inverse of
the map which sends $\Gr_\lambda Y \in P(X)$ to $\lambda$.

Building on the Scannell-Wolf result, the author and Wolf showed
that the $X$-grafting map is also a local homeomorphism, leading to:

\begin{theorem}[{Dumas and Wolf \cite{dumas-wolf:grafting}}]
\label{thm:dumas-wolf}
For each $X \in \T(S)$, the $X$-grafting map $\gr_{\param}X : \ML(S)
\to \T(S)$ is a homeomorphism.  Furthermore, this homeomorphism is
bitangentiable.
\end{theorem}

The last claim in this theorem involves the regularity of the
grafting map as $\lambda$ is varied.  Let $f : U \to V$ be a
continuous map, where $U \subset \R^n$ and $V \subset \R^m$ are open
sets.  The \emph{tangent map}\index{tangent map} of $f$ at $x$, denoted $T_xf : \R^n \to
\R^m$, is defined by
$$ T_xf(v) = \lim_{\epsilon \to 0^+} \frac{f(x + \epsilon v) -
  f(x)}{\epsilon}. $$ The map $f$ is \emph{tangentiable}\index{tangentiable} if this limit
exists for all $(x,v) \in U \times \R^n$, and if the convergence is locally uniform in
$v$ when $x$ is fixed.  Intuitively, a tangentiable map is one which
has one-sided derivatives everywhere.  These notions generalize
naturally to maps between smooth manifolds (e.g.~$\T(S)$), piecewise
linear manifolds (e.g.~$\ML(S)$), or manifolds defined by an atlas of
charts with tangentiable transition functions.  Tangentiable maps and
  manifolds are discussed in \cite{bonahon:variations}.

A homeomorphism $f$ is called \emph{bitangentiable}\index{bitangentiable} if both $f$ and
$f^{-1}$ are tangentiable, and if every tangent map of $f$ or $f^{-1}$
is a homeomorphism.  Thus a bitangentiable homeomorphism is the
analogue of a diffeomorphism in the tangentiable category.

The connection between grafting, projective structures, and
tangentiability was studied by Bonahon, following work of Thurston on the
infinitesimal structure of the space $\ML(S)$
\cite{thurston:minimal-stretch}; the fundamental result, which
strengthens Thurston's theorem, is

\begin{theorem}[Bonahon \cite{bonahon:variations}]
\label{thm:tangentiability}
The projective grafting\index{grafting!projective}\index{grafting} map $\Gr : \ML(S)
\times \T(S) \to \P(S)$ is a bitangentiable homeomorphism.
\end{theorem}

The proof of Theorem \ref{thm:dumas-wolf} uses Theorem
\ref{thm:scannell-wolf}, the above result of Bonahon, and a further
complex linearity property of the tangent map of
projective\index{grafting!projective} and conformal
grafting\index{grafting!conformal}\index{grafting} (see \cite{bonahon:variations}
\cite[\S10]{bonahon:shearing}, also \cite[\S3]{dumas-wolf:grafting}).
This complex linearity provides a ``duality'' between variation of
$\gr_\lambda X$ under changes in $X$ and $\lambda$; in a certain
sense, grafting behaves like a holomorphic function, where $X$ and
$\lambda$ are the real and imaginary parts of its parameter,
respectively.  This allows infinitesimal injectivity of $\gr_\param
X_0$ near $\lambda_0$ to be derived from the infinitesimal injectivity
of $\gr_{\lambda_0}$ near $X_0$.

After applying some additional tangentiable calculus, this
infinitesimal injectivity is converted to local injectivity of
$\gr_{\param}X$, from which Theorem \ref{thm:dumas-wolf} follows by
properness (Theorem \ref{thm:properness}).

\section{Holonomy}\label{sec:holonomy}

We now turn our attention to the holonomy representations of
projective structures in relation to the grafting and
Schwarzian coordinate systems for $\P(S)$.  General references for
these matters include \cite{hejhal:monodromy}
\cite{gunning:affine-projective} \cite{earle:variation}
\cite{hubbard:monodromy} \cite{gkm}.

\subsection{Representations and characters}

Let $\rep(S) = \Hom(\pi_1(S),\PSL_2(\C))$ denote the set of homomorphisms
(representations) from $\pi_1(S)$ to $\PSL_2(\C)$, which is
of an affine $\C$-algebraic variety (as a subset of $(\PSL_2(\C))^N \simeq
(\SO_3(\C))^N)$).  The group $\PSL_2(\C)$ acts algebraically on
$\rep(S)$ by conjugation, and there is a quotient \emph{character
  variety}\index{character variety}
\begin{equation*}
\X(S) = \rep(S) \sslash \PSL_2(\C)
\end{equation*}
in the sense of geometric invariant theory.  Concretely, the points of
$\X(S)$ are in one-to-one correspondence with the set of
\emph{characters}, i.e.~$\C$-valued functions on $\pi_1(S)$ of the form
\begin{equation*}
\gamma \mapsto \tr^2(\rho(\gamma))
\end{equation*}
where $\rho \in \rep(S)$.  Mapping a character to its values on an
appropriate finite subset of $\pi_1(S)$ gives an embedding of $\X(S)$
as an affine variety in $\C^n$.  See \cite{heusener-porti} for a
discussion of $\PSL_2(\C)$ character varieties, building on the work
of Culler-Shalen in the $\SL_2(\C)$ case
\cite{culler-shalen:varieties}.  Algebraic and topological properties
of character varieties are also studied in
\cite[\S9]{gunning:vector-bundles} \cite{goldman:topological}
\cite{bcr}.

\boldpoint{Liftability.}
The variety $\X(S)$ splits into two irreducible components according
to whether or not the associated representations lift from
$\PSL_2(\C)$ to $\SL_2(\C)$ (see \cite{goldman:topological} \cite{bcr}).
Denote these by $\X_0(S)$ and $\X_1(S)$, where the former consists of
liftable characters.  Each of these components has complex dimension
$6g-6$, which agrees with the ``expected dimension'', i.e.~$6g-6 =
(\dim \PSL_2(\C))(N_{\mathrm{gens}} - N_{\mathrm{relators}} - 1)$.

\boldpoint{Elementary and non-elementary.}  When working with the
character variety\index{character variety}, complications may arise
because the invariant-theoretic quotient $\X(S)$ is singular, or
because it is not the same as the quotient set $\rep(S) / \PSL_2(\C)$.
However we can avoid most of these difficulties by restricting
attention to a subset of characters (which contains those that arise
from projective structures).

A representation $\rho \in \rep(S)$ is
\emph{elementary}\index{elementary (character)} if its action on
$\H^3$ by isometries fixes a point or an ideal point, or if it
preserves an unoriented geodesic, otherwise it is
\emph{non-elementary}\index{non-elementary (character)}.

A non-elementary representation is determined up to conjugacy by its
character, so there is a one-to-one correspondence between set set of
conjugacy classes of non-elementary representations and the set $\X'(S)
\subset \X(S)$ of characters of non-elementary representations.

The subset $\X'(S)$ is open and lies in the the smooth locus of the
character variety\index{character variety} \cite{gunning:vector-bundles}
\cite{gunning:affine-projective} \cite{goldman:symplectic-nature}.
Thus $\X'(S)$ is a complex manifold of dimension $6g-6$, and is the
union of the open and closed subsets $\X'_i(S) = \X'(S) \cap \X_i(S)$,
$i=1,2$.

\boldpoint{Fuchsian and quasi-Fuchsian spaces.}  The character variety
$\X(S)$ contains the space $\QF(S)$\index{quasi-Fuchsian space} of
conjugacy classes of quasi-Fuchsian representations of $\pi_1(S)$ as
an open subset of $\X'_0(S)$.  The parameterization of $\QF(S)$ by the
pair of quotient conformal structures gives a holomorphic embedding
$$ \T(S) \times \T(\bar{S}) \to \X(S), $$ where $\bar{S}$ represents
the surface $S$ with the opposite orientation (see
\cite[\S4.3]{matsuzaki-taniguchi}).  In this embedding, the diagonal
$\{ (X, \bar{X}) \: | \: X \in \T(S) \}$ corresponds to the set
$\F(S)$\index{Fuchsian space} of Fuchsian representations, giving an
identification $\F(S) \simeq \T(S)$.  Note that this is \emph{not} a
holomorphic embedding of Teichm\"uller space into the character
variety; the image is a totally real submanifold.

\subsection{The holonomy map}
\label{sec:holonomy-map}

Since the holonomy representation $\rho \in \rep(S)$ of a projective
structure $Z$ is determined up to conjugacy, the associated character
$[\rho] \in \X(S)$ is uniquely determined.  Considering $[\rho]$ as a
function of $Z$ gives the \emph{holonomy map}\index{holonomy map}
\begin{equation*}
\hol : \P(S) \to \X(S).
\end{equation*}
In fact, the image of $\hol$ lies in $\X_0'(S)$: A lift to $\SL_2(\C)$
is given by the linear monodromy of the Schwarzian ODE
\eqref{eqn:schwarzian-ode}.  The holonomy representation is
non-elementary because $S$ does not admit an affine or spherical
structure; for details, see \cite[pp.~297-304]{appell-goursat}
\cite[Thm.~3.6]{kamishima:not-surjective}
\cite[\S2]{gunning:affine-projective}
\cite[Thm.~19,Cor.~3]{gunning:riemann-surfaces}.

\boldpoint{Holonomy theorem.}
For hyperbolic structures on compact manifolds, the holonomy
representation determines the geometric structure.  For projective
structures on surfaces, the same is true \emph{locally}:

\begin{theorem}[Hejhal \cite{hejhal:monodromy}, Earle
    \cite{earle:variation}, Hubbard \cite{hubbard:monodromy}]
\label{thm:hejhal}
The holonomy map $\hol : \P(S) \to \X(S)$ is a local biholomorphism.
\end{theorem}

Originally, Hejhal showed that the holonomy map is a local
homeomorphism using a cut-and-paste argument.  Earle and Hubbard gave
alternate proofs of this result, along with differential calculations
showing that the map is locally biholomorphic.  Recall that when
considering $\P(S)$ as a complex manifold, we are using the complex
structure induced by the quasi-Fuchsian sections.

A more general holonomy theorem for $(G,X)$ structures is discussed in
\cite{goldman:geometric-structures}.

\boldpoint{Negative results.}
Despite the simple local behavior described by Theorem
\ref{thm:hejhal}, the global behavior of the holonomy map is quite
complicated:

\begin{theorem}\mbox{}
\begin{enumerate}
\item \label{item:maskit} The holonomy map is not injective.  In fact,
  all of the fibers of the holonomy map are infinite.
\item \label{item:hejhal} The holonomy map is not a covering of
  its image.
\end{enumerate}
\end{theorem}

The non-injectivity in \eqref{item:maskit} follows from the discussion
of $2\pi$-grafting in \S\ref{sec:fuchsian-holonomy} below.  Hejhal
established \eqref{item:hejhal} by showing that the path lifting
property of coverings fails for the holonomy map
\cite{hejhal:monodromy}.  The infinite fibers of the holonomy map
arise from the existence of \emph{admissible curves}\index{admissible
  curve} that can be used to alter a projective structure while
preserving its holonomy \cite{gkm} \cite[Ch.~7]{kapovich:book}; this
is similar to the ``constructive approach'' discussed in
\S\ref{sec:quasi-fuchsian-holonomy} below.

Further pathological behavior of the holonomy map is discussed in
\cite[\S5]{ito:survey}.

\boldpoint{Surjectivity.}
Of course one would like to know which representations arise from the
holonomy of projective structures.  We have seen that in order to
arise from a projective structure, a character must be non-elementary
and liftable (i.e.~$\hol(\P(S)) \subset \X_0'(S)$).  These necessary
conditions are also sufficient:

\begin{theorem}[{Gallo, Kapovich, and Marden \cite{gkm}}]
\label{thm:gkm}
Every non-elementary liftable $\PSL_2(\C)$-representation of
$\pi_1(S)$ arises from the holonomy of a projective structure on $S$.
Equivalently, we have $\hol(P(S)) = \X_0'(S)$.
\end{theorem}
In the same paper it is also shown that the non-elementary
non-liftable representations arise from \emph{branched} projective
structures\index{branched projective structure}%
\index{projective structure!branched}.  In both cases, the developing
map of a projective structure with holonomy representation $\rho$ is
constructed by gluing together simpler projective surfaces that can be
analyzed directly.  A key technical result that enables this
construction is:
\begin{theorem}[{\cite{gkm}}]
\label{thm:pants}
Let $\rho : \pi_1(S) \to \PSL_2(\C)$ be a homomorphism with
non-elementary image.  Then there exists a pants
decomposition\index{pants decomposition} of $S$
such that the restriction of $\rho$ to any component of the
decomposition is a marked rank-$2$ classical Schottky group.  In
particular, the image of every curve in the decomposition is
loxodromic.
\end{theorem}
Projective structures on pairs of pants with loxodromic boundary
holonomy are analyzed in \cite[\S\S6-7]{gkm}. 

\boldpoint{Holonomy deformations.}  We have seen that projective
structures on a Riemann surface $X$ form an affine space modeled on
$Q(X)$ (\S\ref{sec:affine-naturality}).  Thus, given a non-elementary
representation $\rho \in \X_0'(S)$, projective structures provide
deformations of $\rho$ as follows: Find $Z \in \P(S)$ with
$\hol(Z) = \rho$, which is possible by Theorem \ref{thm:gkm}, and
consider the family of holonomy representations $\{ \hol(Z + \phi) \:
| \: \phi \in Q(X)\}$.  This gives a holomorphic embedding of
$\C^{3g-3}$ into $\X(S)$, a family of \emph{projective
  deformations}\index{projective deformation}
of $\rho$.  (Compare \cite{kra:deformations} \cite{kra:deformations2},
where Kra refers to a projective structure on $X$ as a
\emph{deformation} of the Fuchsian group uniformizing $X$.)

These deformations could be compared with the classical
quasi-conformal deformation theory of Kleinian groups.  Projective
deformations are especially interesting because they are insensitive
to the discreteness of the image of a representation, and because they
apply to quasiconformally rigid Kleinian groups.  On the other hand,
it is difficult to describe the global behavior of a projective
deformation explicitly, and there is often no canonical choice for the
preimage of $\rho$ under the holonomy map.

\subsection{Holonomy and bending}

The holonomy map for projective structures is related to the grafting
coordinate system through the notion of \emph{bending deformations}.
We now describe these deformations, mostly following Epstein and
Marden \cite{epstein-marden:convex-hulls}).  In doing so, we are
essentially re-creating the projective grafting construction of
\S\ref{sec:projective-grafting} while working entirely in hyperbolic
$3$-space, and starting with a Fuchsian representation rather than a
hyperbolic surface.

\boldpoint{Bending Fuchsian groups.}
We begin with an algebraic description of bending.  A primitive
element $\gamma \in \pi_1(S)$ representing a simple closed curve that
separates the surface $S$ determines a $\Z$-amalgamated free product
decomposition
\begin{equation*}
\pi_1(S) = \pi_1(S_1) \free_{\langle \gamma \rangle}
\pi_1(S_2)
\end{equation*}
where $(S - \gamma) = S_1 \sqcup S_2$.  Note that the representative
$\gamma$ determines an orientation of the closed geodesic, and using
this orientation, we make the convention that $S_2$ lies to the right
of the curve.  Given a homomorphism $\rho : \pi_1(S) \to \PSL_2(\C)$
and an element $A \in \PSL_2(\C)$ centralizing $\rho(\gamma)$, there
is a deformed homomorphism $\rho'$ uniquely determined by
\begin{equation}
\label{eqn:quakebend}
 \rho'(x) = 
\begin{cases}
\rho(x) &\text{ if } x \in \pi_1(S_1)\\
A \rho(x) A^{-1} &\text { if } x \in \pi_1(S_2).
\end{cases}
\end{equation}
Similarly, a nonseparating curve $\gamma$ corresponds to a
presentation of $\pi_1(S)$ as an HNN extension, and again each
centralizing element gives a deformation of $\rho$.  See
\cite[\S3]{goldman:ergodic-theory} for further discussion of this
deformation procedure.

When $\rho$ is a Fuchsian representation and $A$ is an elliptic
element having the same axis as $\rho(\gamma)$, the homomorphism
$\rho'$ is a \emph{bending deformation}\index{bending deformation} of
$\rho$.  When $A$ rotates by angle $t$ about the axis of
$\rho(\gamma)$, clockwise with respect to the orientation, we denote
the deformed representation by $\beta_{t \gamma}(\rho) = \rho'$.  Up
to conjugacy, this deformation depends only on the angle $t$ and the
curve $\gamma$, not on the representative in $\pi_1(S)$ or the induced
orientation.

The ``bending'' terminology refers to the geometry of the action of
$\pi_1(S)$ on $\H^3$ by $\beta_{t \gamma}(\rho)$. The Fuchsian
representation $\rho$ preserves a plane $\H^2 \subset \H^3$, whereas
we will see that the bending deformation $\beta_{t \gamma}(\rho)$
preserves a locally convex pleated (or \emph{bent}) plane.

In terms of characters, the Fuchsian representation
$\rho_0$ is a point in $\F(S) \simeq \T(S)$ and bending defines a map
\begin{equation*}
\beta : \scc \times \R^+ \times \T(S) \to \X(S).
\end{equation*}
Like grafting, this map extends continuously to measured laminations
\cite[Thm.~3.11.5]{epstein-marden:convex-hulls}, giving
\begin{equation*}
\beta : \ML(S) \times \T(S) \to \X(S).
\end{equation*}
Note that while the bending path $t \mapsto \beta_{t \gamma}(X)$ is
$2\pi$-periodic, there is no apparent periodicity when bending along a
general measured lamination.

\boldpoint{Earthquakes and quakebends.}  The centralizer of a
hyperbolic M\"obius transformation $\gamma \in \PSL_2(\C)$ contains
all of the elliptic and hyperbolic transformations with the same axis
as $\gamma$, but in defining bending we have only considered the
elliptic transformations.  The deformation corresponding (by formula
\eqref{eqn:quakebend}) to a pure translation is known as an
\emph{earthquake}\index{earthquake}, and the common generalization of
a bending or earthquake deformation (corresponding to the full
centralizer) is a \emph{quakebend} \index{quakebend} or \emph{complex
  earthquake}\index{complex earthquake}.  For further discussion of
these deformations, see \cite{epstein-marden:convex-hulls}
\cite{thurston:earthquakes} \cite{mcmullen:complex-earthquakes}.

\boldpoint{Bending cocycles.}
An alternate definition of the bending deformation makes the geometric
content of the construction more apparent.  Realize the simple closed
curve $\gamma$ as a hyperbolic geodesic on the surface $X \in \T(S)$,
and consider the full preimage $\Tilde{\gamma} \subset \H^2$ of
$\gamma$ in the universal cover; thus $\Tilde{\gamma}$ consists of
infinitely many complete geodesics, the \emph{lifts} of $\gamma$.  By
analogy with the terminology for a pleated plane in $\H^3$, the
connected components of $\H^2 - \Tilde{\gamma}$ will be called
\emph{plaques}\index{plaque}.  For the purposes of this discussion we regard $\H^2$
as a plane in $\H^3$, stabilized by $\PSL_2(\R) \subset \PSL_2(\C)$.

Given $x,y \in (\H^2 - \Tilde{\gamma})$, let $(g_1, \ldots, g_n)$ be
the set of lifts of $\gamma$ that separate $x$ from $y$, ordered
according to the way they intersect the oriented geodesic segment from
$x$ to $y$, with $g_1$ closest to $x$.  Orient each geodesic $g_i$ so
that $y$ lies to the right.  For any $t \in \R$, define the
\emph{bending cocycle}\index{bending cocycle} $B(x,y) \in \PSL_2(\C)$
by
$$ B(x,y) = E(g_1, t) E(g_2, t) \cdots E(g_n,t),$$
where $E(g,t)$ is an elliptic M\"obius transformation with fixed axis
$g$ and clockwise rotation angle $t$.

In case $x$ and $y$ lie in a facet\index{facet}, this empty product is
understood to be the identity.  This construction defines a map $B :
(\H^2 - \Tilde{\gamma}) \times (\H^2 - \Tilde{\gamma}) \to
\PSL_2(\C)$.  Clearly we have $B(x,x) = I$ and $B(x,y)$ only depends
on the plaques\index{plaque} containing $x$ and $y$.  Furthermore, the map $B$
satisfies the \emph{cocycle relation}
\begin{equation}
\label{eqn:bending-cocycle}
B(x,y) B(y,z) = B(x,z) \;\;\text{ for all }\;\; x,y,z \in \H^2 -
\Tilde{\gamma},
\end{equation}
and the \emph{equivariance relation}
\begin{equation}B(\gamma x, \gamma y) = \rho_0(\gamma) B(x,y) \rho_0(\gamma)^{-1}\;\; \text{ for all }
\;\; \gamma \in \pi_1(S)
\label{eqn:bending-equivariance}
\end{equation}
where $\rho_0 \in \F(S)$ represents $Y$.

The connection between the bending cocycle and the bending deformation
described above is as follows (compare
\cite[Lem.~3.7.1]{epstein-marden:convex-hulls}).

\begin{lemma}
\label{lem:bending-deformation-from-cocycle}
Given $Y \in \T(S)$, a simple closed curve $\gamma$, and $t \in
\R$, choose a basepoint $O \in (\H^2 - \Tilde{\gamma})$ and define
$$ \rho(\gamma) = B(O, \gamma O) \rho_0(\gamma), $$ where $\rho_0 \in
\F(S)$ represents $Y$ and $B$ is the bending cocycle associated to
$Y$, $\gamma$, and $t$.  Then $\rho$ is a homomorphism, and it lies in the same
conjugacy class as the bending deformation $\beta_{t \gamma}(Y)$.
\end{lemma}

In other words, the bending cocycle records the ``difference'' between
a Fuchsian character $\rho_0$ and the deformed character $\beta_{t
  \gamma}(\rho_0)$.  The bending cocycle and this lemma extend naturally
to measured laminations \cite[\S3.5.3]{epstein-marden:convex-hulls}.  

\boldpoint{Bending and grafting.}
The key observation relating bending and grafting is that the bending
deformation $\beta_\lambda(Y) : \pi_1(S) \to \PSL_2(\C)$ preserves the
locally convex pleated plane\index{locally convex pleated
  plane}\index{pleated plane!locally convex} in $\H^3$ with intrinsic hyperbolic
structure $Y$ and bending lamination $\lambda$.  In exploring this
connection, let us suppose that $\lambda = t \gamma$ is supported on a
simple closed curve.  The pleating map $\Pl : \H^2 \to \H^3$ can be
defined in terms of the bending cocycle as
$$ \Pl(x) = B(O,x) x$$ where as before $O \in (\H^2 - \Tilde{\gamma})$
is a base point.  Equivariance of this map with respect to $\pi_1(S)$
then follows from Lemma \ref{lem:bending-deformation-from-cocycle} and
the properties \eqref{eqn:bending-cocycle}-\eqref{eqn:bending-equivariance} of the
bending cocycle.  As written, this pleating map is only defined on
$\H^2 - \Tilde{\gamma}$, however it extends continuously to $\H^2$
because on the two sides of a lift $g \subset \Tilde{\gamma}$, the
values of $B(O,\param)$ differ by an elliptic M\"obius transformation
that fixes $g$ pointwise.

The same reasoning shows that the image of $\Pl$ is a locally convex
pleated plane: Since $B$ is locally constant away from
$\Tilde{\gamma}$, the plaques\index{plaque} map into planes in $\H^3$, and when two
such plaques share a boundary geodesic $g$, the images of the plaques
in $\H^3$ meet along a geodesic $\Pl(g)$ with bending angle $t$ (which
is to say, their enveloping planes are related by an elliptic M\"obius
transformation fixing their line of intersection, with rotation angle
$t$). 

We have seen that the holonomy of the projective structure $Z =
\Gr_\lambda Y$ also preserves the equivariant pleated plane in $\H^3$
constructed by bending $\Tilde{Y} \simeq \H^2$ along $\lambda$.
This leads to the fundamental relationship between grafting, bending and the
holonomy map (see \cite[\S2]{mcmullen:complex-earthquakes}):
\begin{equation}
\label{eqn:holonomy-grafting}
\hol(\Gr_\lambda Y) = \beta_\lambda(Y).
\end{equation}
For laminations supported on simple closed curves, this is simply the
observation that the processes of inserting lunes into $\H \subset
\CP^1$ (which gives projective grafting\index{grafting!projective}\index{grafting})
and bending $\H^2 \subset \H^3$ along geodesics (which gives the
bending deformation) are related to one another by the convex hull
construction of \S\ref{sec:thurston-theorem}.  The general equality
follows from this case by continuity of $\hol$, $\Gr$, and $\beta$.

Using \eqref{eqn:holonomy-grafting} we can think of projective
grafting as a ``lift'' of the bending map $\beta : \ML(S) \times \T(S)
\to \X(S)$ through the locally diffeomorphic holonomy map $\hol :
\P(S) \to \X(S)$ (which is not a covering).

\subsection{Fuchsian holonomy}
\label{sec:fuchsian-holonomy}

Let $\P_{\F}(S) = \hol^{-1}(\F(S))$ denote the set of all projective
structures with Fuchsian holonomy.

We can construct examples of projective structures in $\P_{\F}(S)$
using grafting.  Because of the $2\pi$-periodicity of bending along a
simple closed geodesic $\gamma$, the projective structures $\{ \Gr_{2
\pi n \gamma} Y \: | \: n \in \N \}$\index{Fuchsian projective
structure}\index{projective structure!Fuchsian} all have the same Fuchsian
holonomy representation $\rho_0$ (up to conjugacy), which is the
representation uniformizing $Y$.  Of course $n=0$ gives the standard
Fuchsian structure on $Y$.

For $n > 0$ these projective structures have underlying Riemann
surfaces of the form $\gr_{2 \pi n \gamma} Y$, and due to the
$2\pi$-lunes inserted in the projective
grafting\index{grafting!projective}\index{grafting} construction, their developing
maps are surjective.  This construction of ``exotic'' Fuchsian
projective\index{exotic Fuchsian projective
  structure}\index{projective structure!exotic Fuchsian} structures is
due independently to Maskit \cite{maskit:grafting}, Hejhal
\cite[Thm.~4]{hejhal:monodromy}, and Sullivan-Thurston
\cite{sullivan-thurston}.

\boldpoint{Goldman's classification.}
Let $\ML_{\Z}(S)$ denote the countable subset of $\ML(S)$ consisting
of disjoint collections of simple closed geodesics with positive
integral weights.  Generalizing the case of a single geodesic, every
projective structure of the form $\Gr_{2 \pi \lambda} Y$ with $\lambda
\in \ML_{\Z}(S)$ has Fuchsian holonomy.  Goldman showed that all
Fuchsian projective structures arise in this way:

\begin{theorem}[Goldman \cite{goldman:fuchsian-holonomy}]
\label{thm:fuchsian-holonomy}
Let $Z \in \P_\F(S)$ and let $Y = \H^2 / \hol(Z)(\pi_1(S))$ be the
hyperbolic surface associated to the holonomy representation.  Then $
Z = \Gr_{2 \pi \lambda} Y$ for some $\lambda \in \ML_{\Z}$.
\end{theorem}

In terms of the holonomy map $\hol : \P(S) \to \T(S)$, this result
shows that we can identify $\P_{\F}(S)$ with countably many
copies of Teichm\"uller space,
$$ \xymatrix{ \Gr^{-1} : \P_{\F}(S) \ar[r]^-{\simeq} & (2 \pi
\ML_{\Z}(S)) \times \T(S),}$$ and the restriction of the holonomy map
to any one of these spaces $\{2 \pi \lambda \} \times \T(S)$ gives the
natural isomorphism $\T(S) \simeq \F(S)$.  

Alternatively, using Theorem \ref{thm:scannell-wolf} in combination
with Theorem \ref{thm:fuchsian-holonomy}, we can characterize
$\P_{\F}(S)$ as the union of countably many sections of $\pi$,
$$ \P_{\F}(S) = \bigcup_{\lambda \in \ML_{\Z}} \sigma_{2 \pi \lambda}(\T(S)).$$ Note
the difference between these two descriptions of $\P_\F(S)$: In
describing it as a union of sections, we see that the intersection of
$\P_{\F}(S)$ with a fiber $P(X) = \pi^{-1}(X)$ consists of a countable
discrete set naturally identified with $\ML_{\Z}(S)$, whereas in the
holonomy picture we describe the intersection of $\P_{\F}(S)$ with
$\hol^{-1}(Y)$ in similar terms.

Describing $\P_\F(S)$ as a union of the smooth sections $\sigma_{2 \pi
  \lambda}(\T(S))$ of $\pi$ also allows us to conclude that each
  intersection between $\P_\F(S)$ and a fiber $P(X)$ is transverse.
  Previously, Faltings established this transversality result in the
  greater generality of \emph{real holonomy}, that is, the projective
  structures in $\hol^{-1}(\X_\R(S))$ where $\X_\R(S) \subset \X(S)$
  consists of real-valued characters of homomorphisms of $\pi_1(S)$
  into $\PSL_2(\C)$.
\begin{theorem}[{Faltings \cite{faltings:real-projective}}]
Let $Z \in P(X)$ be a projective structure with real holonomy.  Then
$\hol(P(X))$ is transverse to $\X_\R(S)$ at $\hol(Z)$.
\end{theorem}
The characters in $\X_\R(S)$ correspond to homomorphisms that are
conjugate into $\SU(2)$ or $\PSL_2(\R)$.  Both cases include many
non-Fuchsian characters, as a homomorphism $\rho : \pi_1(S) \to
\PSL_2(\R)$ is Fuchsian if and only if its Euler class is maximal,
$e(\rho) = 2g-2$.  Goldman describes the projective structures with
real holonomy in terms of grafting in
\cite[\S2.14]{goldman:fuchsian-holonomy}
\cite[pp.~14-15]{goldman:geometric-structures}.

\subsection{Quasi-Fuchsian holonomy} 
\label{sec:quasi-fuchsian-holonomy}
Let $\P_{\QF}(S) = \hol^{-1}(\F(S))$ denote the set of all projective
structures with quasi-Fuchsian holonomy, which is an open subset of
$\P(S)$.

Goldman's proof of Theorem \ref{thm:fuchsian-holonomy} involves a
study of the topology and geometry of developing maps of Fuchsian
projective structures.  The topological arguments apply equally well
to projective structures with quasi-Fuchsian holonomy, and the
information they provide can be summarized as follows:

\begin{theorem}[{Goldman \cite{goldman:fuchsian-holonomy}}]
\label{thm:quasi-fuchsian-holonomy}
Let $Z \in \P_\QF(S)$ have developing map\index{developing map} $f :
\Tilde{Z} \to \CP^1$, and let $\Lambda \subset \CP^1$ be the limit set
of the holonomy group, a Jordan curve with complementary regions
$\Omega_{\pm}$.  Then:
\begin{enumerate}
\item The quotient of the developing preimage of the limit set, denoted
  $\Lambda(Z) = f^{-1}(\Lambda) / \pi_1(S)$, consists of a finite
  collection of disjoint simple closed curves.

\item The quotient of the developing preimage of $\Omega_-$, denoted 
$Z_- = f^{-1}(\Omega_-) / \pi_1(S)$ consists of a
  finite collection of disjoint homotopically essential annuli bounded
  by the curves in $\Lambda(Z)$.  In particular, the curves in
  $\Lambda(Z)$ are naturally grouped into isotopic pairs.
\end{enumerate}
\end{theorem}

Recall that among the two domains of discontinuity, $\Omega_+$ is
distinguished by the fact that the orientation of its quotient marked
Riemann surface agrees with that of $S$, while that of the quotient of
$\Omega_-$ is opposite.

\begin{figure}
\begin{center}
\includegraphics{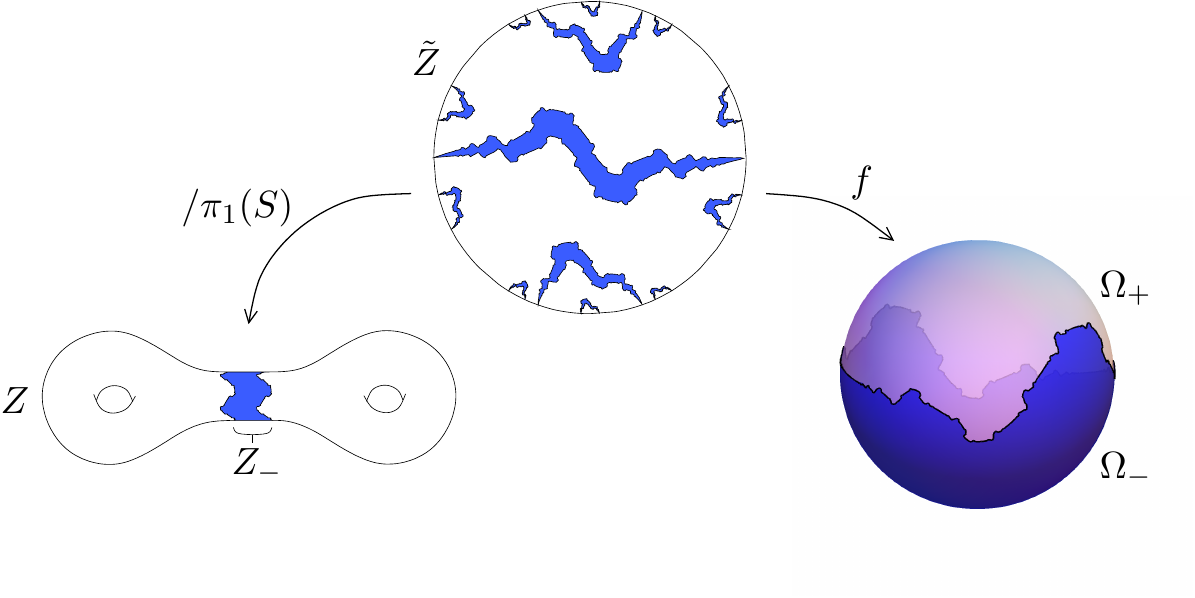}
\end{center}
\caption{The relationship between the developing map $f$ and the
  domains of discontinuity $\Omega_{\pm}$ for a projective structure
  with quasi-Fuchsian holonomy.  In this example, the open set
  $Z_-$ is an annulus, so the wrapping invariant is a simple
  closed curve with unit weight.\label{fig:qfdevel}}
\end{figure}

The topology of a typical (surjective) quasi-Fuchsian developing map
is represented schematically in Figure \ref{fig:qfdevel}.

\boldpoint{Wrapping invariant.}  Given this description of the
preimage of the limit set, there is a natural $\Z$-weighted multicurve
associated to a quasi-Fuchsian projective structure $Z$: Suppose the
collection of annuli $Z_-$ represents homotopy classes $\gamma_1,
\ldots, \gamma_n$, and that there are $n_i$ parallel annuli homotopic
to $\gamma_i$.  Define the \emph{wrapping invariant}\index{wrapping invariant}
$$ \wrap(Z) = \sum_i n_i \gamma_i \in \ML_{\Z}(S).$$ Note that we
could have also defined this using the family of curves $\Lambda(Z)$,
since $2 n_i$ is the number of parallel curves homotopic to
$\gamma_i$. 

Theorem \ref{thm:fuchsian-holonomy} is derived from Theorem
 \ref{thm:quasi-fuchsian-holonomy} by showing that for a Fuchsian
 projective structure $Z$, we have
$$ Z = \Gr_{2 \pi \wrap(Z)} Y,$$ where $Y$ is the quotient of $\H^2$
 by the holonomy group, as above.  In other words, for Fuchsian
 projective structures, the wrapping invariant is the grafting
 lamination (up to a multiple of $2 \pi$).

\boldpoint{Quasi-Fuchsian components.}
Because limit sets vary continuously in $\QF(S)$, the wrapping
invariant is a locally constant function on $\P_{\QF}(S)$.  Thus
$\P_{\QF}(S)$ breaks into countably many subsets
$$ \P_{\QF}(S) = \bigcup_{\lambda \in \ML_\Z(S)} \P_\lambda(S) \text{,
  \; where \;\;} \P_\lambda(S) = \wrap^{-1}(\lambda).$$ We will refer
to these as \emph{components}\index{quasi-Fuchsian
  component} of $\P_{\QF}(S)$.

The quasi-Fuchsian component with zero wrapping invariant,
$\P_{0}(S)$, consists of \emph{standard} quasi-Fuchsian
structures\index{projective structure!standard quasi-Fuchsian}.  The
holonomy map gives a diffeomorphism
$$ \hol : \P_{0}(S) \to \QF(S),$$ where the inverse map associates to
$\rho \in \QF(S)$ the induced projective structure on the quotient
$\Omega_+ / \rho(\pi_1(S))$ of one domain of discontinuity.  The
developing map\index{developing map} of a standard quasi-Fuchsian
projective structure is a Riemann map $f : \H \xrightarrow{\simeq}
\Omega_+$.

The other components $\P_\lambda(S)$, with $\lambda \neq 0$, consist
of \emph{exotic} quasi-Fuchsian projective structures\index{exotic
  quasi-Fuchsian projective structure}%
\index{projective structure!exotic quasi-Fuchsian}; as in the Fuchsian
case, these have surjective developing maps.  Unlike the Fuchsian
case, however, the components $\P_{\lambda}(S)$ do not have a simple
description in terms of the grafting coordinates on $\P(S)$.
Nevertheless, when restricted to one of these components, the holonomy
map
$$ \hol : \P_\lambda(S) \to \QF(S), $$ is again a diffeomorphism.  The
inverse $\QF(S) \to \P_\lambda(S)$ can be constructed by either of two
methods:
\begin{enumerate}
\item \textbf{Constructive approach.} In a generalization of
  $2\pi$-integral projective grafting, one starts with a standard
  quasi-Fuchsian projective structure $Z$ and glues annuli into the
  surface to produce a new projective structure which has the same
  holonomy but which has wrapping invariant $\lambda$.  Allowing the
  starting structure to vary gives a map $\QF(S) \simeq \P_0(S) \to
  \P_{\lambda}(S)$ that is inverse to $\hol$.  See
  \cite[\S1.2]{goldman:fuchsian-holonomy} \cite[\S2.4]{ito:exotic2}
  \cite[Ch.~7]{kapovich:book} for details.

\item \textbf{Deformation approach.}  Starting with a fixed Fuchsian
  representation $\rho_0 \in \F(S)$, any quasi-Fuchsian representation
  $\rho$ can be obtained by a $\rho_0$-equivariant quasiconformal
  deformation.  By pulling back the quasiconformal deformation through
  a developing map\index{developing map}, one can simultaneously deform a Fuchsian
  projective structure $Z_0$ with holonomy $\rho_0$ to obtain a
  quasi-Fuchsian structure $Z$ with holonomy $\rho$.  This deformation
  does not change the wrapping invariant, so starting with $Z_0 =
  \Gr_{2 \pi \lambda} X$ and considering all quasiconformal
  deformations gives the desired map $\QF(S) \to \P_\lambda(S)$.  See \cite[\S3]{shiga-tanigawa}\cite[\S2.5]{ito:exotic}.

\end{enumerate}
Thus the structure of $\P_{\QF}(S)$ is similar to that of $\P_\F(S)$
described above: It consists of countably many connected components
$\P_\lambda(S)$, each of which is diffeomorphic to $\QF(S)$ by the
holonomy map (compare \cite[\S7.2]{kapovich:book} \cite[\S\S2.5-2.6]{ito:exotic}).

\boldpoint{Bumping of quasi-Fuchsian components.}
We say that two components $\P_\lambda(S)$ and $\P_\mu(S)$ \emph{bump}\index{bumping}
if their closures intersect, i.e. if $\bar{\P_\lambda(S)} \cap
\bar{\P_\mu(S)} \neq \emptyset$; an element of the intersection
is called a \emph{bumping point}.  A component $\P_\lambda(S)$
\emph{self-bumps}\index{self-bumping} at $Z \in \P(S)$ if $U \cap \P_\lambda(S)$ is
disconnected for all sufficiently small neighborhoods $U$ of $Z$.
These terms are adapted from similar phenomena in the theory of
deformation spaces of Kleinian groups (surveyed in
\cite{canary:bumponomics}, see also \cite{anderson-canary:bumping}
\cite{anderson-canary-mccullough:bumping} \cite{bromberg-holt:bumping}
\cite{holt:bumping}).

The bumping of quasi-Fuchsian components has been studied by McMullen
\cite{mcmullen:complex-earthquakes}, Bromberg-Holt
\cite{bromberg-holt:projective}, and Ito \cite{ito:exotic}
\cite{ito:exotic2}.  The basic problem of determining which component
pairs bump is resolved by:
\begin{theorem}[Ito \cite{ito:exotic2}]\mbox{}
\begin{enumerate}
\item For any $\lambda, \mu \in \ML_\Z(S)$, the components
  $\P_\lambda(S)$ and $\P_\mu(S)$ bump.
\item For any $\lambda \in \ML_\Z(S)$, the component $\P_\lambda(S)$
  self-bumps at a point in $\bar{\P_0(S)}$.
\end{enumerate}
\end{theorem}
The bumping points constructed in the proof of this theorem are all
derived from a construction of Anderson-Canary that illustrates the
difference between algebraic and geometric convergence for Kleinian
groups \cite{anderson-canary:bumping}.  This construction was first
applied to projective structures by McMullen
 to give an example of bumping
between $\P_\lambda(S)$ and $\P_0(S)$ \cite{mcmullen:complex-earthquakes}.  The holonomy representations
for these bumping examples have accidental parabolics but are not
quasiconformally rigid; recently, Brock, Bromberg, Canary, and Minsky
have shown that these conditions are necessary for bumping
\cite{bbcm:bumping} (compare \cite{ohshika:bumping}).

\subsection{Discrete holonomy}
Let $\D(S) \subset \X(S)$ denote the set of characters of discrete
representations, and let $\P_\D(S)$ denote the set of projective
structures with discrete holonomy\index{discrete holonomy}.  Since
$\F(S) \subset \QF(S) \subset \D(S)$, we have corresponding inclusions
$$ \P_\F(S) \subset \P_\QF(S) \subset \P_\D(S).$$

Because $\hol$ is a local diffeomorphism, topological properties of
$\D(S)$ correspond to those of $\P_\D(S)$.  For example, $\D(S)$ is
closed (see \cite{jorgensen:discrete} \cite{chuckrow}), and its
interior is the set $\QF(S)$ of quasi-Fuchsian representations
\cite{sullivan:qc2} \cite{bers:holomorphic-families}.  Thus
$\P_\D(S)$ is a closed subset of $\P(S)$ with interior $\P_\QF(S)$.

If $Z \in \P_\D(S)$ has holonomy $\rho$, then the associated pleated
plane $\Pl(Z) : \H^2 \to \H^3$ is invariant under the holonomy group
$\Gamma = \rho(\pi_1(S))$ and descends to a \emph{locally convex pleated
surface}\index{locally convex pleated surface}\index{pleated
surface!locally convex} in the quotient hyperbolic manifold $M = \H^3 / \Gamma$:
$$ \xymatrix{
  \H^2 \ar[d]^-{/\pi_1(S)} \ar[r]^-{\Pl(Z)} & \H^3 \ar[d]^-{/\Gamma}\\
  Y \ar[r] & M } $$ Here $Y \in \T(S)$ is the hyperbolic surface such
that $Z = \Gr_\lambda Y$ for some $\lambda \in \ML(S)$.

The pleated surface arising from a projective structure $Z$ with
 discrete holonomy may be one of the connected components of the
 boundary of the \emph{convex core}\index{convex core} of the associated hyperbolic
 manifold $M$.  If so, the projective surface $Z$ is the component of
 the ideal boundary of $M$ on the ``exterior'' side of the pleated
 surface.  Conversely, the ideal boundary and convex core boundary
 surfaces in a complete hyperbolic manifold are related by grafting
 (see \cite[\S5.1]{scannell-wolf:grafting}
 \cite[\S2.8]{mcmullen:complex-earthquakes}).

For more general projective structures with discrete holonomy, the
pleated surface need not be embedded in the quotient manifold, however
it must lie within the convex core (see \cite[\S5.3.11]{ceg:notes-on-notes}).

In addition to the Fuchsian and quasi-Fuchsian cases described above,
projective structures with other classes of discrete holonomy
representations have found application in Kleinian groups and
hyperbolic geometry.  For example, projective structures with
degenerate holonomy are used in Bromberg's approach to the Bers
density conjecture \cite{bromberg:degenerate}, and those with Schottky
holonomy are used in Ito's study of sequences of Schottky groups
accumulating on Bers' boundary of Teichm\"uller space
\cite{ito:schottky}.

\subsection{Holonomy in fibers}

In contrast to the complicated global properties of the holonomy map,
its restriction to a fiber is very well-behaved:
\begin{theorem}
For each $X \in \T(S)$, the
  restriction $\left.\hol\right|_{P(X)}$ is a proper holomorphic
  embedding, whose image $\hol(P(X))$ is a complex-analytic subvariety
  of $\X(S)$.
\end{theorem}

As stated, this theorem incorporates several related but separate
results: Working in the context of systems of linear ODE on a fixed
Riemann surface, Poincar\'e showed that the holonomy map is injective
\cite[p.~310]{appell-goursat} (see also \cite{kra:generalization}
\cite[Thm.~15]{hejhal:monodromy}).  Gallo, Kapovich, and Marden showed
that the image is a complex-analytic subvariety \cite{gkm}, following
an outline given by Kapovich \cite{kapovich:monodromy}; when combined
with injectivity, this implies properness.  Tanigawa gave a more
geometric argument establishing properness of
$\left.\hol\right|_{P(X)}$ when considered as a map into the space
$\X'(S)$ of non-elementary characters \cite{tanigawa:divergence}.
Tanigawa's argument relies on the existence of loxodromic pants
decompositions\index{pants decomposition} (Theorem \ref{thm:pants}),
which was announced in \cite{kapovich:monodromy} and proved in
\cite{gkm}.

\boldpoint{Fuchsian and quasi-Fuchsian holonomy in fibers.}
\label{sec:fuchsian-fibers}
For any $X \in \T(S)$, let $P_\D(X) = P(X) \cap \P_\D(S)$ denote the
set of projective structures with discrete holonomy and with
underlying complex structure $X$.  Similarly, we define $P_\QF(X)$ and
$P_\F(X)$ as the subsets of $P(X)$ having quasi-Fuchsian and Fuchsian
holonomy, respectively.

We have already seen (in \S\ref{sec:fuchsian-holonomy}) that the
$P_\F(X)$ consists of the countable discrete set of projective
structures $\{ \sigma_{2 \pi \lambda}(X) \: | \: \lambda \in \ML_{2
\pi \Z}(S) \}$.  Since the holonomy map is continuous, and $\QF(S)$ is
an open neighborhood of $\F(S)$ in $\X(S)$, each of these Fuchsian
points has a neighborhood in $P(X)$ consisting of quasi-Fuchsian
projective structures with the same wrapping invariant.  Elements of
$P_\F(X)$ are sometimes called \emph{Fuchsian centers}\index{Fuchsian center} (or
\emph{centers of grafting}\index{center of grafting} \cite{anderson:projective}), because they
provide distinguished center points within these ``islands'' of
quasi-Fuchsian holonomy (see \cite[\S13]{dumas:schwarzian}
\cite[Thm.~6.6.10]{marden:outer-circles}).

\begin{figure}
\begin{center}
\setlength{\subfigbottomskip}{0cm} \subfigure[The Bers embedding of
the square punctured
torus.]{\fbox{\includegraphics[height=3.5cm]{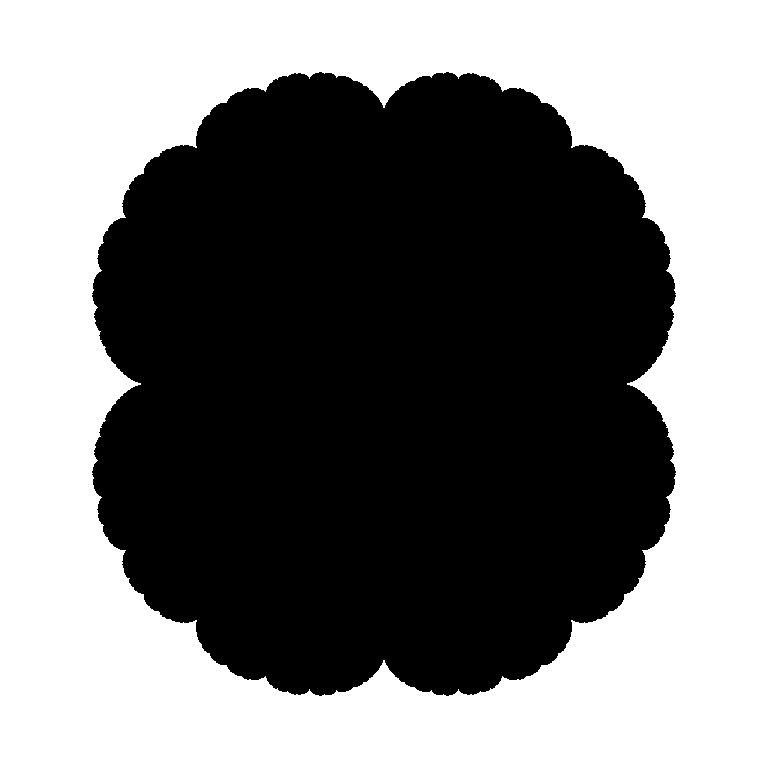}}}
\hspace{1mm}
\subfigure[In this larger view, the Bers embedding of a punctured
torus with a short geodesic appears as a small dot (center) surrounded
by many islands of exotic quasi-Fuchsian projective
structures.]{\fbox{\includegraphics[height=3.5cm]{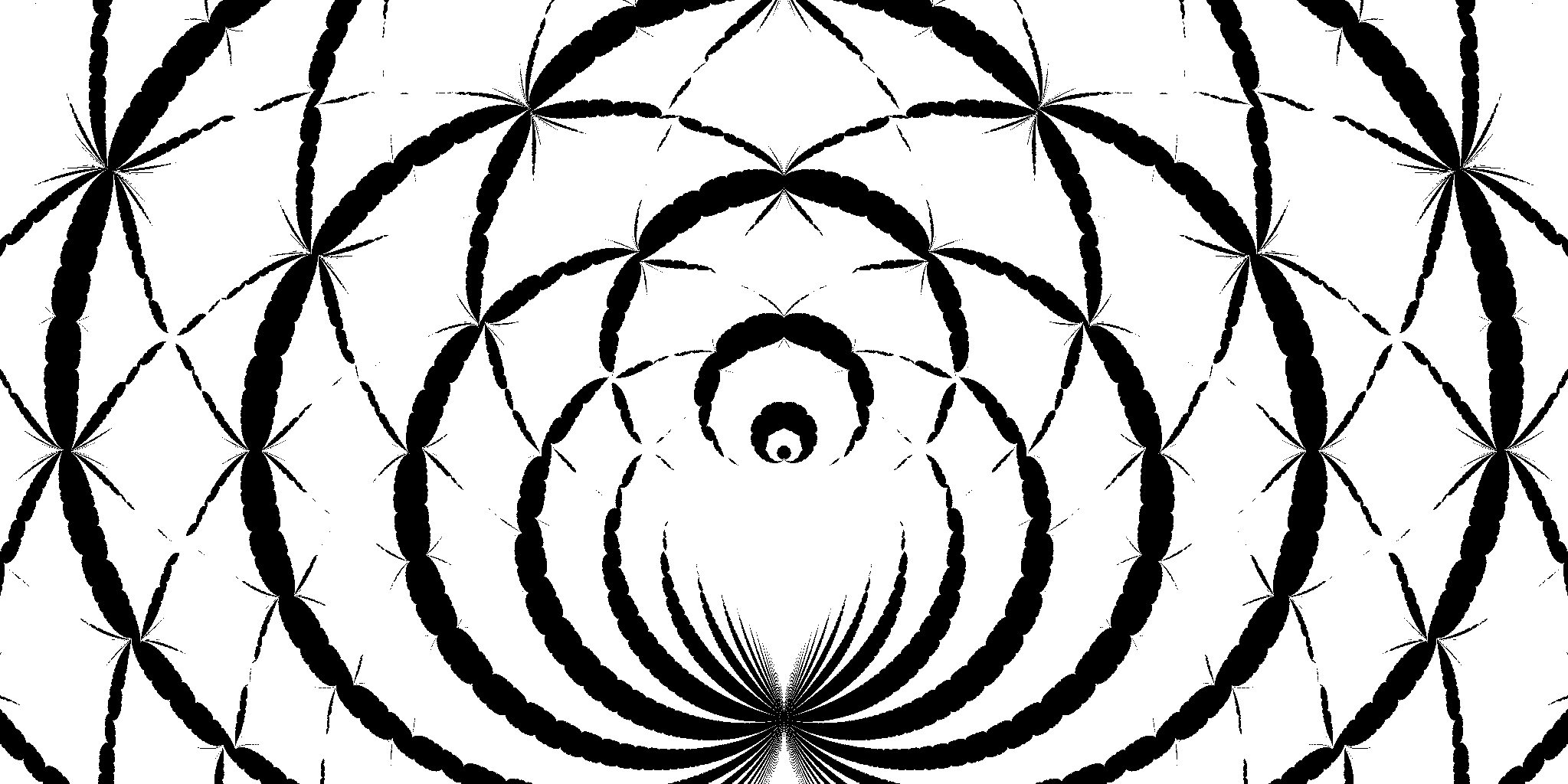}}}
\end{center}
\caption{Islands of quasi-Fuchsian holonomy\index{islands of
    quasi-Fuchsian holonomy} in $P(X) \simeq \C$ (where $X$ is a
  punctured torus) exhibit complicated structure at small and large
  scales.  These images were created using the software package
  \emph{Bear} \cite{dumas:bear}.\label{fig:bers} }
\end{figure}

Using the Schwarzian parameterization, the intersection $\P_0(S) \cap
  P(X)$, consisting of the standard quasi-Fuchsian projective
  structures on $X$, can be considered as an open set $B_X \subset
  Q(X) \simeq \C^{3g-3}$.  This set is the image of the holomorphic
  \emph{Bers embedding}\index{Bers embedding} of Teichm\"uller space \cite{shiga}, and in
  particular it is connected and contractible.  We also have $B(1/2)
  \subset B_X \subset B(3/2)$, where $B(r) = \{ \phi \in Q(X)
  \:|\:\|\phi\|_\infty < r \}$, as a consequence of Nehari's theorem
  \cite{nehari:schwarzian}.  See Figure \ref{fig:bers} for examples of
  Bers embeddings of the Teichm\"uller space of punctured tori.

\begin{figure}
\begin{center}
\renewcommand{\thesubfigure}{}
\subfigure[$X_1$]{\fbox{\includegraphics[width=0.3\textwidth]{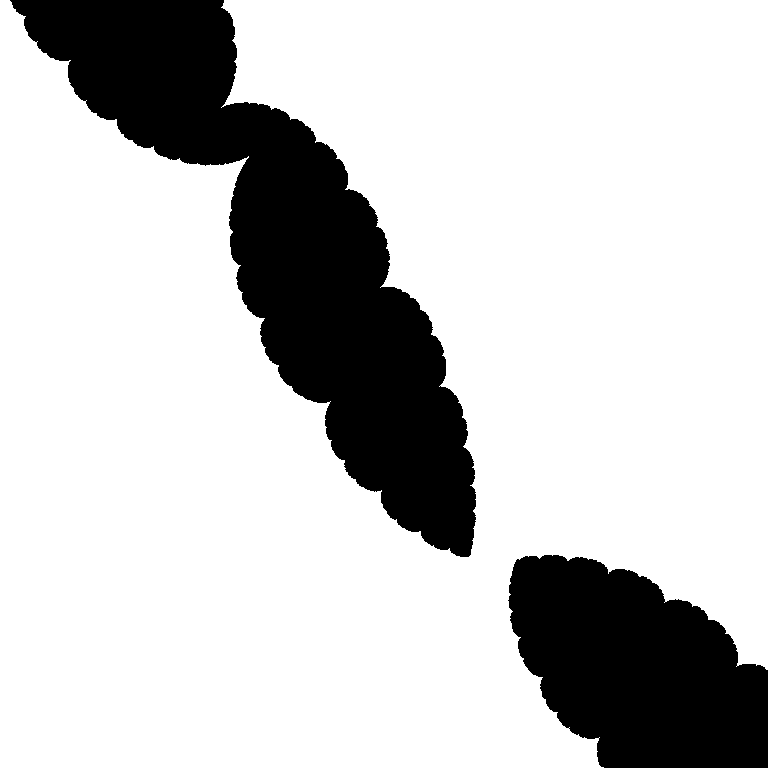}}}
\subfigure[$X_2$]{\fbox{\includegraphics[width=0.3\textwidth]{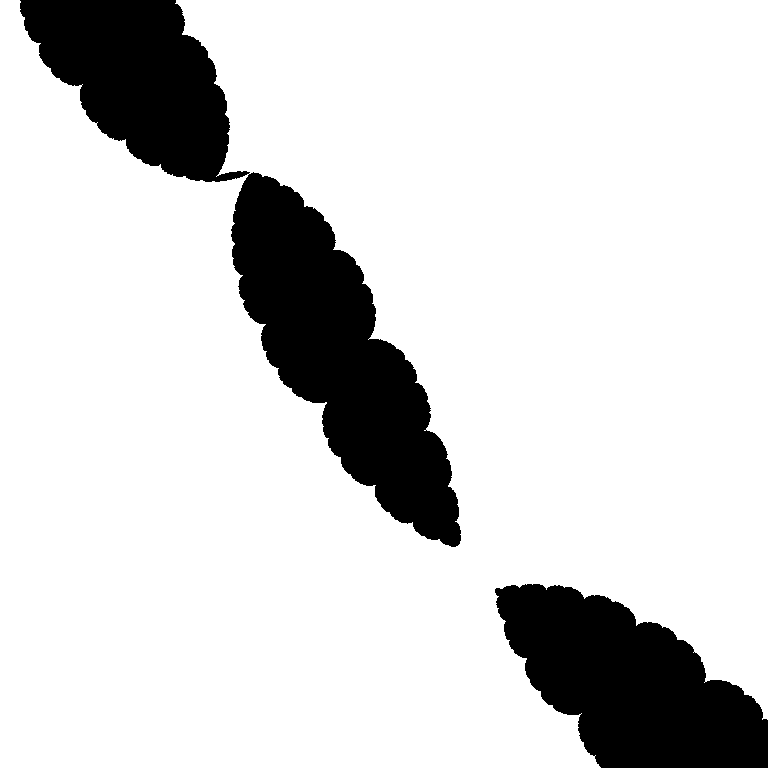}}}
\subfigure[$X_3$]{\fbox{\includegraphics[width=0.3\textwidth]{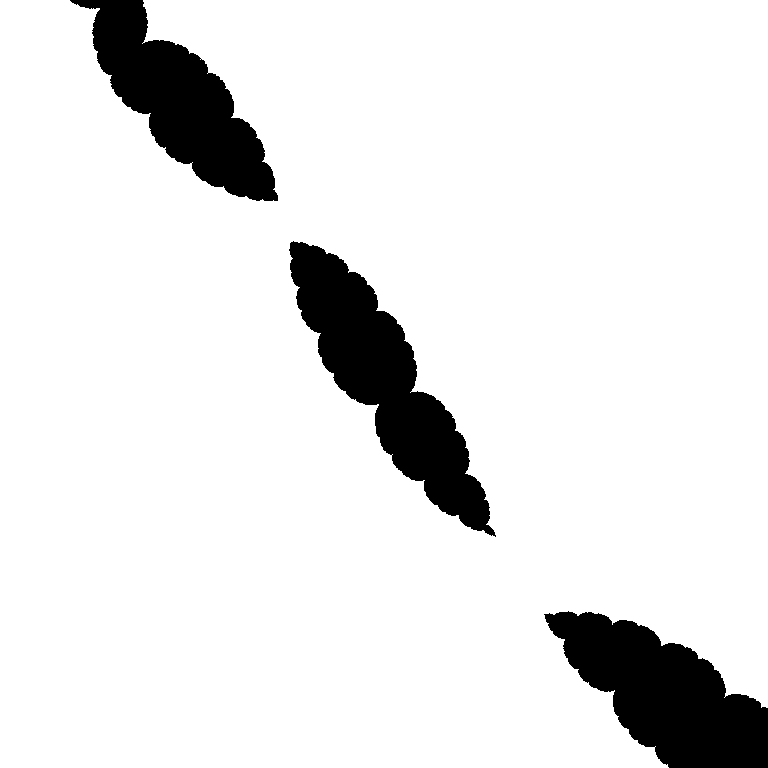}}}
\end{center}
\caption{Islands of quasi-Fuchsian holonomy in $P(X)$
  appear to break apart as the complex structure $X$ is changed,
  suggesting that some islands do not contain Fuchsian
  centers\index{Fuchsian center}.  Each
  image shows a small square in $P(X_i) \simeq \C$, where
  $\{X_1,X_2,X_3\}$ are closely-spaced points in the Teichm\"uller
  space of the punctured torus.\label{fig:disconnected}}
\end{figure}

For $\lambda \neq 0$, it is not known whether the set $\P_\lambda(S)
  \cap P(X)$ is connected (or bounded), though experimental evidence
  in the punctured case suggests that it often has many connected
  components, and that the structure of the connected components
  changes with $X$ (see Figure \ref{fig:disconnected}). Of course,
  only one component contains the Fuchsian structure
  $\sigma_{2 \pi \lambda}(X)$.

\boldpoint{Quasi-Fuchsian versus discrete in a fiber.}
In the space of all projective structures, the quasi-Fuchsian
structures form the interior of the set with discrete
holonomy\index{discrete holonomy}.  The same relationship holds for
$P_\QF(X)$ and $P_\D(X)$.

\begin{theorem}[Shiga and Tanigawa \cite{shiga-tanigawa}, Matsuzaki \cite{matsuzaki}]
  For any $X \in \T(S)$, we have $P_\QF(X) = \interior(P_\D(X))$.
\end{theorem}

In comparing these sets, one inclusion is immediate: Since
$\interior(\P_\D(S)) = \P_\QF(S)$, we have $\interior(P_\D(X)) \supset
P_\QF(X)$.  The opposite inclusion is more subtle.  Each component of
the interior of $P_\D(X)$ necessarily consists of quasiconformally
conjugate, discrete, faithful representations without accidental
parabolics.  However there exist $(3g-3)$-dimensional holomorphic
families of \emph{singly degenerate}\index{singly degenerate} surface
groups in $\X(S)$ which satisfy these conditions, but which are not
quasi-Fuchsian.  Such a family could account for an open subset of
$P_\D(X)$ (in either of two topologically distinct ways
\cite{matsuzaki}), and a key step in the proof of the theorem is to
exclude this possibility.

\section{Comparison of parameterizations}\label{sec:comparison}

\subsection{Compactifications}

\boldpoint{Compactification of $\ML(S)$.}
The space of measured laminations has the structure of a
cone: The group $\R^+$ acts by scaling the transverse measure
($\lambda \mapsto t \lambda$, $t \in \R^+$) and the empty lamination
$0 \in \ML(S)$ is the unique fixed point of this action.  The orbit of
a nonzero lamination is a \emph{ray} in $\ML(S)$.  The space of rays,
$$ \PML(S) = (\ML(S) - \{0\}) / \R^+,$$ or \emph{projective measured
  laminations}\index{projective measured lamination}\index{measured lamination!projective} forms a natural boundary for $\ML(S)$.  We say that a
sequence $\lambda_i \in \ML(S)$ converges to $[\lambda] = \R^+\cdot
\lambda \in \PML(S)$ if there exists a sequence of positive real
numbers $c_i$ such that $c_i \to 0$ and $ c_i \lambda_i \to \lambda$
in $\ML(S)$.  The induced compactification\index{compactification!of
  measured laminations}
$$ \bar{\ML(S)} = \ML(S) \cup \PML(S) $$ is homeomorphic to a
closed ball, with interior $\ML(S) \simeq \R^{6g-6}$ and boundary
$\PML(S) \simeq S^{6g-7}$.  See \cite[Ch.~3]{penner-harer} for further
discussion of the spaces $\ML(S)$ and $\PML(S)$, and \cite{flp} for
related discussion of the space of measured foliations, which is
naturally identified with $\ML(S)$ (as described in \cite{levitt}
\cite[\S11.8-11.9]{kapovich:book}).

\boldpoint{Compactification of $\T(S)$.}
Recall that $\scc$ denotes the set of isotopy classes of simple closed
curves on $S$, or equivalently, the simple closed geodesics of any
hyperbolic structure on $S$.  Thurston defined a compactification of
$\T(S)$ using the hyperbolic length map
\begin{equation*} \begin{split}
\L : \T(S) &\to \R^{\scc}\\
X &\mapsto ( \ell(\gamma,X) )_{\gamma \in \scc}. 
  \end{split}
\end{equation*}
This map is an embedding, as is its projectivization 
$$ \PL : \T(S) \to \mathbb{P}^+\R^\scc = (\R^\scc - \{0\}) / \R^+,$$ and in
each case, a suitable finite subset of $\scc$ suffices to determine
the image of a point.  The boundary $\partial \PL(\T(S))$ coincides with the image of $\PML(S)$ under the
projectivization of the embedding
\begin{equation*}
\begin{split}
\ML(S) &\to \R^\scc\\
\lambda &\mapsto ( i(\gamma,\lambda) )_{\gamma \in \scc}
\end{split}
\end{equation*}
where $i(\lambda,\gamma)$ denotes the total mass of $\gamma$ with
respect to the transverse measure of $\lambda$.  This gives the
\emph{Thurston compactification}\index{Thurston
  compactification}\index{compactification!of Teichm\"uller space}
$$ \bar{\T(S)} = \T(S) \cup \PML(S)$$ which has the topology of a
closed $(6g-6)$-ball.  Concretely, a sequence $X_n \to \infty$ in
Teichm\"uller space converges to $[\lambda] \in \PML(S)$ if for every
pair of simple closed curves $\alpha,\beta \in \scc$ we have
$$ \frac{\ell(\alpha,X_i)}{\ell(\beta,X_i)} \to
\frac{i(\alpha,\lambda)}{i(\beta,\lambda)} $$ whenever the right hand
side is well-defined (i.e. $i(\beta,\lambda) \neq 0$).  A detailed
discussion of the Thurston compactification can be found in
\cite[Exp.~7-8]{flp} (see also \cite{thurston:hyp2}
\cite{thurston:bulletin} \cite{bonahon:geodesic-currents}
\cite[Ch.~11]{kapovich:book} \cite[\S5.9]{marden:outer-circles}).

\boldpoint{Compactification of $Q(X)$.}  Since the vector space $Q(X)$
has an action of $\R^+$ by scalar multiplication, it supports a
natural compactification\index{compactification!of quadratic
  differentials} analogous to that of $\ML(S)$; in this case, the
boundary is the space of rays
$$ \PQ(X) = (Q(X) - \{0\}) / \R^+ $$ and we obtain $\bar{Q(X)} =
Q(X) \cup \PQ(X)$ which is homeomorphic to a closed ball.

\subsection{Quadratic differentials and measured laminations}

\boldpoint{The Hubbard-Masur theorem.}
\label{sec:hubbard-masur}
For any $X \in \T(S)$, there is a natural map
$$ \Lambda : Q(X) \to \ML(S)$$ which is defined by a two-step
procedure: First, a quadratic differential $\phi$ has an associated
\emph{horizontal foliation}\index{horizontal foliation}\index{foliation!horizontal} $\F(\phi)$, a singular foliation on $X$
which integrates the distribution of vectors $v \in TX$ such that
$\phi(v) \geq 0$.  This foliation is equipped with a transverse
measure, induced by integration of $|\Im \sqrt{\phi}|$.  In a local
coordinate where $\phi = dz^2$, the foliation is induced by the
horizontal lines in $\C$, with transverse measure $|dy|$.  Zeros of
$\phi$ correspond to singularities of the foliation, where three or
more half-leaves emanate from a point.  See e.g.~\cite[\S5.3,
\S11.3]{kapovich:book} \cite[\S2.2,Ch.~11]{gardiner:book} for a
discussion of quadratic differentials and their measured foliations.

Now lift the horizontal foliation of $\phi$ to the universal cover
$\Tilde{X} \simeq \H^2$.  Each non-singular leaf of the lifted
foliation is a uniform quasi-geodesic, so it is a bounded distance
from unique hyperbolic geodesic.  The hyperbolic geodesics obtained in
this way---the \emph{straightening}\index{straightening} of
$\F$---form the lift of a geodesic lamination on $X$, and the
transverse measure of the foliation induces a transverse measure on
this lamination in a natural way \cite{levitt}.  The result is a
measured lamination $\Lambda(\phi) \in \ML(S)$, which we call the
\emph{horizontal lamination}\index{horizontal
lamination}\index{measured lamination!horizontal} of $\phi$.

The same constructions can be applied to the distribution of vectors
satisfying $\phi(v) \leq 0$, which gives the \emph{vertical
  foliation}\index{vertical foliation}\index{foliation!vertical}
and \emph{vertical lamination}\index{vertical
  lamination}\index{measured lamination!vertical} of $\phi$.  The former is induced by
the foliation of $\C$ by vertical lines in local coordinates such that
$\phi = dz^2$.  Note that multiplication by $-1$ in $Q(X)$ exchanges
vertical and horizontal: for example, the horizontal lamination of
$-\phi$ is the vertical lamination of $\phi$.

The strong connection between quadratic differentials and measured
laminations is apparent in:

\begin{theorem}[Hubbard and Masur \cite{hubbard-masur}]
\label{thm:hubbard-masur}
For each $X \in \T(S)$, the map $\Lambda : Q(X) \to \ML(S)$ is a
homeomorphism.  In particular, every measured lamination is realized
by a unique quadratic differential on $X$.
\end{theorem}
\index{Hubbard-Masur theorem}

Note that Hubbard and Masur work with measured foliations rather than
measured laminations; the statement above incorporates the
aforementioned straightening procedure to identify the two notions.

We call the inverse of $\Lambda$ the \emph{foliation
  map}\index{foliation map}, denoted
$\phi_F : \ML(S) \to Q(X)$.  Note that the definition of both
$\Lambda$ and $\phi_F$ depend on the choice of a fixed conformal
structure $X$, but we suppress this dependence in the notation.

Since the transverse measure of $\Lambda(\phi)$ is obtained by
integrating $|\Im \sqrt{\phi}|$, these maps have the following homogeneity
properties:
\begin{equation*}
\begin{split}
\Lambda(c \phi) &= c^{\frac{1}{2}} \Lambda(\phi)\\
\phi_F(c \lambda) &= c^2 \phi_F(\lambda)\\ 
\end{split}
\end{equation*}
for all $c \in \R^+$.  Therefore $\Lambda$ and $\phi_F$ descend to
mutually inverse homeomorphisms between the spaces of rays $\PML(S)$
and $\PQ(X)$, and we also use $\Lambda$ and $\phi_F$ to denote these
induced maps.

\begin{figure}
\begin{center}
\includegraphics[width=0.9\textwidth]{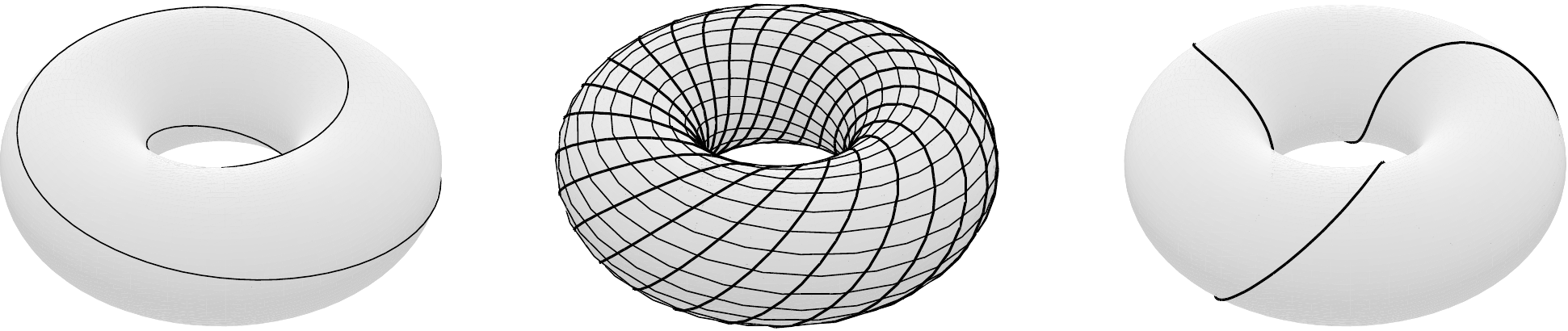}
\end{center}
\caption{A pair of closed curves on a compact Riemann surface of genus
 $1$ (as seen here on the far left and right) are ``orthogonal'' if
 they are isotopic to leaves of a pair of orthogonal geodesic
 foliations of the Euclidean metric (center).  This situation is
 non-generic; typically, at least one of the two foliations will have
 dense leaves.  For surfaces of higher genus, orthogonality of
 measured laminations is defined similarly, however there are many
 distinct singular Euclidean metrics.\label{fig:orthog}}
\end{figure}

\boldpoint{Orthogonality and the antipodal map.}
\label{sec:antipodal}
Given $X \in \T(S)$, a pair of measured laminations $\lambda, \mu \in
\ML(S)$ is \emph{orthogonal with respect to $X$}\index{orthogonality} if there exists $\phi
\in Q(X)$ such that
\begin{equation*}
\begin{split}
\Lambda(\phi) = \lambda\\
\Lambda(-\phi) = \mu
\end{split}
\end{equation*}
That is, $\lambda$ and $\mu$ appear as the horizontal and vertical
laminations of a single holomorphic quadratic differential on $X$.
(Compare the torus case shown in Figure \ref{fig:orthog}.)

By Theorem \ref{thm:hubbard-masur}, two laminations $\lambda$ and
$\mu$ are orthogonal with respect to $X$ if and only if 
$$ \phi_F(\lambda) = - \phi_F(\mu) \in Q(X).$$ Thus the
homeomorphism $\phi_F : \ML(S) \to Q(X)$ turns orthogonal pairs into
opposite quadratic differentials, and the set of $X$-orthogonal pairs is
the graph of the \emph{antipodal involution}\index{antipodal involution} $i_X : \ML(S) \to \ML(S)$
defined by
$$ i_X(\lambda) = \Lambda( - \phi_F(\lambda) ).$$

By homogeneity of $\Lambda$ and $\phi_F$, the antipodal map descends
to $i_X : \PML(S) \to \PML(S)$.  We say $[\lambda],[\mu] \in \PML(S)$
are orthogonal with respect to $X$ if $i_X([\lambda]) = [\mu]$.  See
\cite{dumas:antipodal} for further discussion of the antipodal map and
orthogonality.

\subsection{Limits of fibers}

Using the projective grafting\index{grafting!projective}\index{grafting} homeomorphism
$\Gr: \ML(S) \times \T(S) \to \P(S)$, we can regard $\bar{\ML(S)}
\times \bar{\T(S)}$ as a compactification of $\P(S)$.  This is the
\emph{grafting compactification}%
\index{grafting compactification}\index{compactification!grafting}.

Given $X \in \T(S)$, the fiber $P(X) \subset \P(S)$ corresponds to a
set of pairs $\Gr^{-1}(X) = \{ (\lambda,Y) \: | \: \gr_\lambda Y = X
\}$ in the grafting coordinates.  Since $P(X)$ is a distinguished
subset of the Schwarzian parameterization of $\P(S)$, studying its
behavior in the grafting parameterization is one way to study the
relationship between these two coordinate systems.  The asymptotic
behavior of $P(X)$ can be described in terms of orthogonality:

\begin{theorem}[Dumas \cite{dumas:antipodal}]
\label{thm:antipodal}
Let $(\lambda_n,Y_n) \in \ML(S) \times \T(S)$ be a divergent sequence
such that $\Gr_{\lambda_n} Y_n \in P(X)$ for all $n$.  Then
$$ \lim_{n \to \infty} \lambda_n = [\lambda] \;\; \text{ if and only
  if } \;\;
  \lim_{n \to \infty} Y_n = i_X([\lambda]), $$ where these limits are
  taken in $\bar{\ML(S)}$ and $\bar{\T(S)}$,
  respectively.\\ In particular, the boundary of $P(X)$ in the
  grafting compactification of $\P(S)$ is the graph of the antipodal
  involution $i_X : \PML(S) \to \PML(S)$.
\end{theorem}

This theorem can be considered as evidence of compatibility between
the grafting coordinates for $\P(S)$ and the foliation of $\P(S)$ by
fibers of $\pi$.  For example, we have:

\begin{corollary}
The closure of $P(X)$ in $\bar{\ML(S)} \times \bar{\T(S)}$
is homeomorphic to a closed ball of dimension $6g-6$.
\end{corollary}

The proof of Theorem \ref{thm:antipodal} in \cite{dumas:antipodal} is
essentially a study of the collapsing and
co-collapsing\index{co-collapsing map}\index{collapsing map} maps of a
complex projective structure, and their relation to the harmonic maps
variational problem.  We now describe this variational technique,
and then outline the main steps in the proof.

\boldpoint{Harmonic maps.}
\label{sec:harmonic-maps}
Let $(M,g)$ and $(N,h)$ be complete
Riemannian manifolds, and assume that $M$ is compact.  If $f : M \to
N$ is a smooth map, the \emph{energy}\index{energy (of a map)} of $f$ is defined by
$$ \E(f) = \frac{1}{2} \int_M \| df(x) \|_2^2 dg(x).$$ The map $f$ is
\emph{harmonic}\index{harmonic map} if it is a critical point of the
energy functional.  If $N$ is also compact and $h$ has negative
sectional curvature, then any nontrivial homotopy class of maps $M \to
N$ contains a harmonic map, and this map is an absolute minimum of the
energy functional in the homotopy class
\cite{eells-sampson:existence}.  Furthermore, the harmonic map is
unique in its homotopy class, unless the image of $M$ is a closed
geodesic in $N$, in which case there is a $1$-parameter family of
harmonic maps obtained by rotation.  General references for the theory
of harmonic maps include \cite{eells-lemaire:report}
\cite{eells-lemaire:report2} \cite{schoen:analytic-aspects}, with
particular applications to Teichm\"uller theory surveyed in
\cite{daskalopoulos-wentworth}.

\boldpoint{Equivariant harmonic maps.}  If $\pi_1(M)$ acts by
isometries on a Riemannian manifold $\hat{N}$, then we can define the
energy of an equivariant map $\Tilde{M} \to \hat{N}$ by integration of
$\|df\|_2^2$ over a fundamental domain for the action of $\pi_1M$ by
deck transformations.  This generalizes the energy of smooth maps $M
\to N$, because the action of $\pi_1M$ on $\hat{N}$ need not have a
Hausdorff quotient.  Existence of harmonic maps%
\index{equivariant harmonic map} is more delicate in this case, but
can sometimes be recovered under additional restrictions on the group
action.  For example if $M$ is a surface and $N = \H^3$ is equipped
with the isometric action coming from a non-elementary representation
$\rho : \pi_1(S) \to \PSL_2(\C)$, then there is a unique equivariant
harmonic map $h : \tilde{S} \to \H^3$ \cite{donaldson:harmonic}.

\boldpoint{Singular targets.}
Korevaar and Schoen developed a deep generalization of the theory of
harmonic maps\index{harmonic map} in which the Riemannian manifold $N$ is replaced by a
nonpositively curved (\npc, also
known as locally \catzero\index{locally CAT-0@locally $CAT(0)$}) metric
space \cite{korevaar-schoen:1} \cite{korevaar-schoen:2}
\cite{korevaar-schoen:3}.  Here the energy functional is approximated
by the average squared distance between the image of a point $x \in M$
and the image of a small sphere centered at $x$.  Inequalities
comparing distances in \npc\ spaces to those in Euclidean space have an
essential role in the development of this theory.

Generalizing the Riemannian case, we have the following equivariant
existence and uniqueness results: If $\hat{N}$ is a locally compact
\npc\ space on which $\pi_1(M)$ acts by isometries without fixing any
equivalence class of rays, then there is an equivariant
harmonic\index{equivariant harmonic map} map $h : \tilde{M} \to
\hat{N}$, which is Lipschitz and energy-minimizing
\cite{korevaar-schoen:2}.  If furthermore $\hat{N}$ is negatively
curved (locally \catk, for some $\kappa < 0$), then the
harmonic map\index{harmonic map} is unique unless its image is a geodesic
\cite{mese:uniqueness}.

\boldpoint{Harmonic maps from surfaces.}
When $M$ is $2$-dimensional, the energy functional depends only on the
conformal class of the metric $g$, so it makes sense to consider
harmonic maps from Riemann surfaces to Riemannian manifolds and
nonpositively curved metric spaces.  An important invariant of a
harmonic map $f : X \to (N,h)$ from a Riemann surface is its \emph{Hopf
differential}\index{Hopf differential}
\begin{equation}
\label{eqn:hopf}
\Phi(f) = [f^*(h)]^{2,0}
\end{equation}
which is a holomorphic quadratic differential.  In the Riemannian
case, the holomorphicity of $\Phi(f)$ is a consequence of the
Euler-Lagrange equation of the energy functional
\cite[\S10]{eells-lemaire:report}.  With a suitable generalization of
the pullback metric (see \cite[\S2.3]{korevaar-schoen:1}), a
holomorphic Hopf differential is also obtained from a harmonic
map\index{harmonic map} to a \npc\index{nonpositively curved (NPC)}\ 
metric space (compare \cite[\S5]{mese:harmonic-maps-into-singular}).

We can use the same formula \eqref{eqn:hopf} to define a Hopf
differential for any smooth map $X \to (N,h)$, which can be further
generalized to maps with $L^2$ distributional derivatives, and to
finite-energy maps to \npc\  metric spaces
\cite[Thm.~2.3.1]{korevaar-schoen:1}.  The result is a $L^1$
measurable quadratic differential that is not necessarily holomorphic.

\boldpoint{Harmonic maps and dual trees.}  Recall from
\S\ref{sec:dual-tree} that for each $\lambda \in \ML(S)$ we have a
dual $\R$-tree $T_\lambda$.  This tree is a \npc\  metric space (even
\catk\ for all $\kappa<0$) equipped with an isometric action of
$\pi_1(S)$.  The Hubbard-Masur construction of a quadratic
differential on $X \in \T(S)$ with lamination $\lambda$ can be
described in terms of an equivariant harmonic
\index{equivariant harmonic map} map $X \to T_\lambda$\index{harmonic map}.

\begin{theorem}[Wolf \cite{wolf:realizing-measured-foliations},
    Daskalopoulos-Dostoglou-Wentworth \cite{ddw:r-trees}]
\label{thm:wolf-tree}
Let $h : \Tilde{X} \to T_\lambda$ be an equivariant harmonic map to
the dual $\R$-tree of $\lambda \in \ML(S)$.  Then $\phi_F(\lambda) = -4
\Phi(h)$.
\end{theorem}

\boldpoint{Harmonic maps and the Thurston compactification.}  The
Thurston compactification of Teichm\"uller space can also be
characterized in terms of Hopf differentials\index{Hopf differential}
of harmonic maps\index{harmonic map} from a fixed Riemann surface as follows:
\begin{theorem}[Wolf \cite{wolf:thesis}]
\label{thm:wolf-thurston}
Fix $X \in \T(S)$ and let $Y_n \to \infty$ be a divergent sequence in
$\T(S)$.  Let $\Phi_n = \Phi(h_n)$ be the Hopf differential of the
harmonic map $h_n : X \to Y_n$ compatible with the markings.  Then
$$ \Lambda(-\Phi_n) \to [\lambda] \in \PML(S) \;\;\text{if and only
  if} \;\; Y_n \to [\lambda] \in \PML(S).$$
\end{theorem}

\boldpoint{Collapsing, co-collapsing, and harmonic
  maps.}\index{collapsing map}\index{co-collapsing map} Using the
harmonic maps\index{harmonic map} results presented above, we now
describe the main steps of the proof of Theorem \ref{thm:antipodal} in
\cite{dumas:antipodal}.  For simplicity, we will suppose that
$\Gr_{\lambda_n} Y_n \in P(X)$ and that both grafting coordinates have
limits in $\PML(S)$, i.e.
$$ \lim_{n \to \infty} \lambda_n = [\lambda] \;\; \lim_{n \to \infty}
Y_n = [\mu], $$ and we outline a proof that $i_X([\lambda]) = [\mu]$.
The stronger statement of the theorem is derived from the same set of
ideas.

\begin{proof}[Outline of proof of Theorem \ref{thm:antipodal}.]
\mbox{}
\begin{enumerate}
\item Both the collapsing maps $\kappa_n : X \to Y_n$ and the
  co-collapsing maps $\hat{\kappa}_n : \Tilde{X} \to T_{\lambda_n}$ are
  \emph{$C$-almost harmonic}\index{almost harmonic
    map}, meaning that their energies exceed the
  minimum energies in their homotopy classes by at most $C$.  Here $C$
  is a constant that depends only on the topology of $S$.  (Compare \cite{tanigawa:grafting}.)

\item The maps $\kappa_n$ and $\hat{\kappa}_n$ have an orthogonality
  relationship: their derivatives have rank $1$ in the same subset of
  $X$ (the Euclidean part of the Thurston metric), and in this set,
  the collapsed directions of $\kappa_n$ and $\hat{\kappa}_n$ are
  orthogonal.  This orthogonality relationship is expressed in terms
  of their Hopf differentials\index{Hopf differential} as 
\begin{equation}
\label{eqn:opposite}
\Phi(\kappa_n) + \Phi(\hat{\kappa}_n)  = 0.
\end{equation}

\item Let $h_n : X \to Y_n$ and $\hat{h}_n : \Tilde{X} \to
  T_{\lambda_n}$ denote the harmonic maps homotopic to $\kappa_n$ and
  $\hat{\kappa}_n$, respectively.  Then the projective limit of Hopf
  differentials $[\Phi] = \lim_{n \to \infty} \Phi(h_n)$ satisfies
  $[\Lambda(-\Phi)] = [\mu]$ by Theorem \ref{thm:wolf-thurston}.
  Similarly, by Theorem \ref{thm:wolf-tree}, the projective limit
  $[\hat{\Phi}] = \lim_{n \to \infty} \Phi(\hat{h}_n) = \lim_{n \to
    \infty} (- \phi_F(\lambda_n) / 4)$ satisfies
  $[\Lambda(-\hat{\Phi})] = [\lambda].$
 
\item Since the pair of almost harmonic maps $\kappa_n$ and
  $\hat{\kappa}_n$ have opposite Hopf differentials, one might expect
  that the associated harmonic maps\index{harmonic map} $h_n$ and $\hat{h}_n$ have
  ``almost opposite'' Hopf differentials.  Suppose that this is true in
  the sense of projective limits, i.e.~that
\begin{equation}
\label{eqn:opposite-phi}
 [\Phi] = [-\hat{\Phi}] \in \PQ(X).
\end{equation}
Then we would have $[\Lambda(\Phi)] = [\lambda]$ and $[\Lambda(-\Phi)]
= [\mu]$, or equivalently, that $i_X([\lambda]) = [\mu]$, completing
the proof.  Thus we need only derive \eqref{eqn:opposite-phi}.

\item The norm of the difference between the pullback metric of a
  $C$-almost harmonic map $f$ to an \npc\  space and that of its
  homotopic harmonic map $h$ is $O(C^{1/2} \E(h)^{1/2})$ as $\E(h) \to
  \infty$ (by an estimate of Korevaar and Schoen, see \cite[\S2.6]{korevaar-schoen:1}).  Phrasing this in terms of
  Hopf differentials, which are the $(2,0)$ parts of the pullback
  metrics, and using that $|\E(h) - 2 \|\Phi(h)\|| = O(1)$, we have
$$ \| \Phi(f) - \Phi(h) \|_1 \leq C' (1 + \|\Phi(h)\|_1^{\frac{1}{2}}).$$ In
  particular the norm of the difference is much smaller than either
  term as $\|\Phi(h)\| \to \infty$, and so the Hopf differentials of
  any sequence of $C$-almost harmonic maps with energy tending to
  infinity has the same projective limit as the Hopf differentials of
  the harmonic maps.  Applying this to the collapsing and
  co-collapsing\index{collapsing map}\index{co-collapsing map} maps, and using \eqref{eqn:opposite}, we have
$$ [\Phi] = \lim_{n \to \infty} \Phi(\kappa_n) = \lim_{n \to \infty}
  (- \Phi(\hat{\kappa}_n)) = [-\hat{\Phi}],$$
and \eqref{eqn:opposite-phi} follows.
\end{enumerate}
\end{proof}

\subsection{Limits of the Schwarzian}
We now connect the previous discussion of asymptotics of grafting
coordinates for $P(X)$ with the complex-analytic parameterization of
$\P(S)$.  Let $\bar{P(X)}$ denote the \emph{Schwarzian
  compactification}\index{Schwarzian compactification}\index{compactification!Schwarzian} of $P(X)$ obtained by attaching $\PQ(X)$ using the
limiting behavior of the Schwarzian derivative, i.e.~a sequence $Z_n
\in P(X)$ converges to $[\phi]$ if $(Z_n - Z_0) \to [\phi]$ in the
topology of $\bar{Q(X)}$.  Here $Z_0$ denotes an arbitrary
basepoint, which is used to identify $P(X)$ with $Q(X)$; the limit of
a sequence in $\PQ(X)$ does not depend on this choice.  Note that this
construction only compactifies the individual fibers of $\P(S)$, but
does not compactify $\P(S)$ itself.

There is a natural guess for the relationship between the Schwarzian
compactification and the closure of $P(X)$ in the grafting
compactification:  The boundary of the latter is the set of
$X$-antipodal pairs in $\PML(S) \times \PML(S)$, and each
$X$-antipodal pair arises from a ray in the space of quadratic
differentials, so one might expect a boundary point $[\phi] \in
\PQ(X)$ to correspond to the pair consisting of its vertical and
horizontal laminations.  The following makes this intuition precise:

\begin{theorem}[Dumas \cite{dumas:schwarzian}]
\label{thm:schwarzian-compactification}
The grafting and Schwarzian compactifications of $P(X)$ are naturally
homeomorphic, and the boundary map $\PQ(X) \to \PML(S) \times \PML(S)$
is given by
$$ [\phi] \mapsto ([\Lambda(-\phi)], [\Lambda(\phi)]).$$ That is, for
a divergent sequence in $P(X)$, the limit of the vertical
(resp.~horizontal) laminations of Schwarzian differentials is equal to
the limit of the measured laminations (resp.~hyperbolic structures) in
the grafting coordinates.
\end{theorem}

This result about compactifications involves a comparison between two
homeomorphisms $\ML(S) \to Q(X)$.  One of these we have already
seen---the foliation map\index{foliation map} $\phi_F$ which sends
$\lambda \in \ML(S)$ to a quadratic differential whose horizontal
foliation has straightening $\lambda$ (\S\ref{sec:hubbard-masur}).
The other homeomorphism is derived from the Schwarzian
parameterization of projective structures as follows.  Recall (from
\S\ref{sec:conformal-grafting}) that there is a homeomorphism
$\sigma_\param(X) : \ML(S) \to P(X)$ with the property that
$\sigma_\lambda(X) \in P(X)$ is a projective structure with grafting
lamination $\lambda$.  Using $\sigma_0(X)$ as a basepoint, we compose
with the Schwarzian parameterization $P(X) \simeq Q(X)$ to obtain the
\emph{Thurston map}\index{Thurston map}:
\begin{equation*}
\begin{split}
\phi_T : \ML(S) &\to Q(X)\\
\lambda &\mapsto \left(\sigma_\lambda(X) - \sigma_0(X)\right)
\end{split}
\end{equation*}
The Thurston map is a homeomorphism, and it satisfies $\phi_T(0) = 0$,
but unlike the foliation map there is no \emph{a priori} reason for
$\phi_T$ to map rays in $\ML(S)$ to rays in $Q(X)$.  However, the
Thurston map does preserve rays in an asymptotic sense:
\begin{theorem}[\cite{dumas:schwarzian}]
\label{thm:schwarzian}
For any $X \in \T(S)$, the foliation and Thurston maps are
asymptotically proportional.  Specifically, there exists a constant
$C(X)$ such that 
$$\| \phi_F(\lambda) + 2 \phi_T(\lambda) \|_1 \leq C(X) \left( 1 +
\|\phi_F(\lambda)\|_1^{\frac{1}{2}} \right )$$ for all $\lambda \in \ML(S)$.
\end{theorem}

Before discussing the proof of Theorem \ref{thm:schwarzian}, we
explain the connection with compactifications.  In terms of the
Thurston map, Theorem \ref{thm:schwarzian-compactification} asserts
that if $\phi_T(\lambda_n) = \Gr_{\lambda_n} Y_n$ is a divergent
sequence in $P(X)$, then we have
\begin{equation}
\label{eqn:phi-t}
\begin{split}
\lim_{n \to \infty} \lambda_n &= \lim_{n \to \infty}
\Lambda(-\phi_T(\lambda_n)) \in \bar{\ML(S)} \;\;\text{ and }\\
\lim_{n \to \infty} Y_n &= \lim_{n \to \infty}
\Lambda(\phi_T(\lambda_n)) \in \bar{\T(S)}.
\end{split}
\end{equation}
Theorem \ref{thm:antipodal} has already given a similar
characterization in terms of the map $\phi_F$; we have
\begin{equation}
\label{eqn:phi-f}
\begin{split}
\lim_{n \to \infty} \lambda_n &= \lim_{n \to \infty}
\Lambda(\phi_F(\lambda_n)) \in \bar{\ML(S)}\;\;\text{ and }\\
\lim_{n \to \infty} Y_n &= \lim_{n \to \infty}
\Lambda(-\phi_F(\lambda_n)) = \lim_{n \to \infty} i_X(\lambda_n) \in \bar{\T(S)},
\end{split}
\end{equation}
where the first line is trivial since $\Lambda \circ \phi_F =
\mathrm{Id}$, and the second line follows from the definition of the
antipodal map (\S\ref{sec:antipodal}).  However, since $\phi_F$ and
$\phi_T$ are asymptotically proportional by a negative constant
(Theorem \ref{thm:schwarzian}), the limit characterizations
\eqref{eqn:phi-t} and \eqref{eqn:phi-f} are equivalent, and Theorem
\ref{thm:schwarzian-compactification} follows.  See
\cite[\S14]{dumas:schwarzian} for details.

\boldpoint{Thurston metrics and the Schwarzian.}
\label{sec:schwarzian-proof}
We now sketch the main ideas involved in the proof of Theorem
\ref{thm:schwarzian}.  The proof is essentially a study of the
Thurston metric\index{Thurston metric} on a complex projective surface (see
\S\ref{sec:thurston-metric}).  Recall that the goal is to show that
$\| \phi_F(\lambda) + 2 \phi_T(\lambda) \|_1 \leq C(X)
\varepsilon(\lambda)$ where $\varepsilon(\lambda)$ is defined by
$$\varepsilon(\lambda) = 1 + \|\phi_F(\lambda)\|_1^{\frac{1}{2}}.$$

\begin{proof}[Outline of proof of Theorem \ref{thm:schwarzian}.]
\mbox{}
\begin{enumerate}
\item The functions $\varepsilon(\lambda)$ and $\lambda \mapsto \|
  \phi_F(\lambda) + 2\phi_T(\lambda) \|_1$ are continuous on $\ML(S)$.
  Since weighted simple closed geodesics are dense in $\ML(S)$, it
  suffices establish an inequality relating these functions for such
  weighted geodesics, and the general case follows by continuity.
  Thus we will assume $\lambda$ is a weighted simple closed geodesic
  for the rest of the proof.

\item Associated to such $\lambda$ we have the following objects:
\begin{itemize}
\item   The Thurston metric $\rho_\lambda$ of the projective structure
  $\sigma_\lambda(X) \in P(X)$
\item The decomposition $X = X_{0} \sqcup X_{-1}$ of $X$ into
  Euclidean and hyperbolic parts of $\rho_\lambda$.   Here $X_{0}$ is
  an open cylinder, the union of the $1$-dimensional strata\index{strata} in the
  canonical stratification\index{canonical stratification}.

\item The collapsing map $\kappa : X \to Y_\lambda =
  \gr_{\lambda}^{-1}(X)$ and its Hopf differential\index{Hopf differential} $\Phi(\kappa)$, which
    is a measurable (non-holomorphic) quadratic differential supported
    on $X_0$.  

\item The ratio of conformally equivalent metrics $\rho_\lambda /
  \rho_0$, a well-defined positive function on $X$.  Here $\rho_0$ is
  the hyperbolic metric.
\end{itemize}

\item The Schwarzian derivative $\phi_T(\lambda)$ of the projective
  structure $\sigma_\lambda(X)$ decomposes as a sum of two terms,
\begin{equation}
\label{eqn:decomposition}
\phi_T(\lambda) = -2 \Phi(\kappa) + 2 \B(\log(\rho_\lambda / \rho_0)),
\end{equation}
where the second-order differential operator $\B$ is defined by 
$$ \B(\eta) = \left [ \Hess(\eta) - d \eta \tensor d \eta
\right]^{2,0}. $$ In this expression, the Hessian is computed using
the hyperbolic metric $\rho_0$.  This decomposition follows from the
cocycle property for a generalization of the Schwarzian
derivative\index{Schwarzian derivative}%
\index{Osgood-Stowe Schwarzian derivative} introduced by Osgood and
Stowe \cite{osgood-stowe}.

\item The harmonic map\index{harmonic map} estimate from the proof of Theorem
  \ref{thm:antipodal} shows that the first term of the decomposition
  \eqref{eqn:decomposition} is approximately proportional to
  $\phi_F(\lambda)$.  Specifically, we have
\begin{equation}
\| \phi_F(\lambda) - 4 \Phi(\kappa) \|_1 \leq C \varepsilon(\lambda).
\label{eqn:harmonic-bound}
\end{equation}
Therefore it suffices to show that the $L^1$ norm of $\beta =
\B(\log(\rho_\lambda / \rho_0))$ is also bounded by a multiple of
$\varepsilon(\lambda)$.

\item By the definition of $\B$ and the Cauchy-Schwartz inequality,
  the $L^1$ norm of $\beta$ is bounded by the $L^2$ norms of the
  Hessian and gradient of $\log(\rho_\lambda / \rho_0)$ with respect
  to the hyperbolic metric.  By standard elliptic theory, these are in
  turn bounded by the $L^2$ norms of $\log(\rho_\lambda / \rho_0)$ and
  its Laplacian.

\item The Laplacian of $\log(\rho_\lambda / \rho_0)$ is essentially
  the difference of the curvature $2$-forms of $\rho_\lambda$ and
  $\rho_0$ (compare \eqref{eqn:curvature} above, also
  \cite{huber:zum}).  For large grafting, the surface $X$ is
  dominated by its Euclidean part, forcing most of the curvature of
  $\rho_\lambda$ to concentrate near a finite set of points.

\item This curvature concentration phenomenon provides a bound for the norm
  $\|\Delta \log(\rho_\lambda / \rho_0)\|_{L^2(D)}$ on a hyperbolic disk $D
  \subset X$ of definite size.  A bound on $\| \log(\rho_\lambda /
  \rho_0)\|_{L^2(D)}$ follows using a weak Harnack inequality, completing the
  local estimate $\|\beta\|_{L^1(D)} < C(X)$.

\item Finally, we make the local estimate global: If $\beta$ were
  holomorphic, then we would have $\|\beta\|_{L^1(X)} \leq C'(X)
  \|\beta\|_{L^1(D)}$ by compactness of the unit sphere in $Q(X)$.
  While $\beta$ is not holomorphic, the decomposition
  \eqref{eqn:decomposition} and the estimate
  \eqref{eqn:harmonic-bound} show that $\beta$ is close to a
  holomorphic quadratic differential, with difference of order
  $\varepsilon(\lambda)$.  Combining this with the holomorphic case, we
  obtain $\| \beta \|_{L^1(X)} \leq C(X) \varepsilon(\lambda)$, completing the proof.
\end{enumerate}
\end{proof}

\subsection{Infinitesimal compatibility}

In this final section we discuss infinitesimal aspects of the map between the
grafting and analytic coordinate systems for $\P(S)$.   

The forgetful projection $\pi : \P(S) \to \T(S)$ can be thought of as
a coordinate function in the Schwarzian parameterization of
$\P(S)$.  The other ``coordinate'' in this parameterization is an
element of the fiber $Q(X)$ of the bundle of quadratic differentials,
but lacking a canonical trivialization for this bundle, there is no
associated global coordinate map.\index{coordinate map}

On the other hand, in the grafting coordinate system, we have a pair
of well-defined coordinate maps $p_\ML : \P(S) \to \ML(S)$ and $p_\T :
\P(S) \to \T(S)$, which are defined by the property that the inverse of
projective grafting\index{grafting!projective}\index{grafting} is $\Gr^{-1}(Z) = (p_\ML(Z), p_\T(Z)) \in \ML(S)
\times \T(S)$.

The fiber of $p_\ML$ over $\lambda$ consists of the projective
structures $\{ \Gr_\lambda Y \: | \: Y \in \T(S) \}$.
Since $\Gr_\lambda : \T(S) \to \P(S)$ is a smooth map, these fibers
are smooth submanifolds of $\P(S)$.

The fiber of $p_\T$ over $Y$ consists of the projective structures $\{
\Gr_\lambda Y \: | \: \lambda \in \ML(S) \}$.  Bonahon showed that
$\lambda \mapsto \Gr_\lambda Y$ includes $\ML(S)$ into $\P(S)$
tangentiably (see Theorem \ref{thm:tangentiability}).  However, the
fibers of $p_\T$ have even more regularity than one might expect from
this tangentiable parameterization:
\begin{theorem}[{Bonahon \cite[Thm.~3, Lem.~13]{bonahon:variations}}]
\label{thm:regularity}
For each $Y \in \T(S)$, the set $p_\T^{-1}(Y)$ is a $C^1$ submanifold
of $\P(S)$.
\end{theorem}
Compare \cite[\S4]{dumas-wolf:grafting}.

Note that each of the three coordinate maps $\pi, p_\ML, p_\T$
projects $\P(S)$ onto a space of half its real dimension, i.e.~each
has both range and fibers of real dimension $6g-6$.  Thus one might
expect that for any two of these maps, the pair of fibers intersecting
at a generic point $Z \in \P(S)$ would have transverse tangent spaces
that span $T_Z \P(S)$.  In fact, this is true at every point, and
furthermore we have:
\begin{theorem}[{Dumas and Wolf \cite{dumas-wolf:grafting}}]\mbox{}
\label{thm:transversality}\index{transversality}
\begin{enumerate}
\item The maps $\pi, p_\ML, p_\T$ have pairwise transverse fibers.
\item The fiber of any one of them projects homeomorphically by each
  of the others.  Moreover, such a projection is a $C^1$
  diffeomorphism whenever its range is $\T(S)$, and is a
  bitangentiable\index{bitangentiable} homeomorphism when the range is $\ML(S)$.
\item The product of any two of these maps gives a homeomorphism from
  $\P(S)$ to a product of two spaces of real dimension $6g-6$.
\end{enumerate}
\end{theorem}

As before, we refer to Bonahon (see \cite[\S 2]{bonahon:variations})
for details about tangentiability, while limiting our focus to its
geometric consequences.  Also note that statement (1) of the theorem
does \emph{not} involve tangentiability, and only makes sense for
fibers of $p_\T$ due to Theorem \ref{thm:regularity}.

We sketch the proof of this theorem; the details we omit can be found
in \cite[Thms.~1.2, 4.1, 4.2, Cor.~4.3]{dumas-wolf:grafting}.

\begin{proof}[Sketch of proof of Theorem \ref{thm:transversality}]
Statement (3) follows because the inverse map for each pair of
coordinates can be written explicitly in terms of $\Gr$, $\gr_\lambda$,
and $\gr_\param X$ and their inverses (which exist by Theorems
\ref{thm:thurston}, \ref{thm:scannell-wolf}, and \ref{thm:dumas-wolf},
respectively).  For example, $p_{\T} \times \pi : \P(S) \to \T(S)
\times \T(S)$ is a homeomorphism with inverse
$$ (X, Y) \mapsto \Gr_{(\gr_\param X)^{-1}(Y)} X.$$ Similarly, the map
$(\lambda,X) \mapsto \sigma_\lambda(X)$ is inverse to $p_{\ML} \times
\pi$.

Statement (3) also shows that the restrictions of maps considered in
statement (2) are homeomorphisms.  To show that each case with target
$\T(S)$ is actually a diffeomorphism, it is enough to show that the
derivative of the restriction has no kernel (by the inverse
function theorem).  This kernel is the intersection of tangent spaces
to fibers of two coordinate maps, thus this case will follow from
statement (1).  Similar reasoning applies in cases with target
$\ML(S)$, where one deduces bitangentiability from transversality
using a criterion of Bonahon \cite[Lem.~4]{bonahon:variations}.

Thus the proof is reduced to the transversality statement (1), which
has one case for each pair of coordinate maps.  The pair $(p_\ML,
p_\T)$ follows easily from Thurston's theorem and the tangentiability
of grafting (Theorems \ref{thm:thurston} and
\ref{thm:tangentiability}).  For $(\pi,p_\ML)$ or $(\pi, p_\T)$, a
vector in the intersection of tangent spaces lies in the kernel of a
tangent map of either $\gr_\lambda$ or $\gr_\param X$, which must
therefore be zero, by Theorems \ref{thm:scannell-wolf} and
\ref{thm:dumas-wolf}.
\end{proof}

\nocite{} 
\newcommand{\removethis}[1]{}


\begin{thebibliography}{100}

\bibitem{ahlfors-bers}
L.~Ahlfors and L.~Bers.
\newblock Riemann's mapping theorem for variable metrics.
\newblock {\em Ann. of Math. (2)}, 72:385--404, 1960.

\bibitem{anderson:projective}
C.~Anderson.
\newblock {\em Projective Structures on Riemann Surfaces and Developing Maps to
  $\H^3$ and $\CP^n$}.
\newblock PhD thesis, University of California at Berkeley, 1998.

\bibitem{anderson-canary:bumping}
J.~Anderson and R.~Canary.
\newblock Algebraic limits of {K}leinian groups which rearrange the pages of a
  book.
\newblock {\em Invent. Math.}, 126(2):205--214, 1996.

\bibitem{anderson-canary-mccullough:bumping}
J.~Anderson, R.~Canary, and D.~McCullough.
\newblock The topology of deformation spaces of {K}leinian groups.
\newblock {\em Ann. of Math. (2)}, 152(3):693--741, 2000.

\bibitem{appell-goursat}
P.~Appell, {\'E}.~Goursat, and P.~Fatou.
\newblock {\em Th{\'e}orie des {F}onctions {A}lg{\'e}briques, {T}ome 2:
  {F}onctions {A}utomorphes}.
\newblock Gauthier-Villars, Paris, 1930.

\bibitem{bers:fiber-spaces}
L.~Bers.
\newblock Fiber spaces over {T}eichm\"uller spaces.
\newblock {\em Acta. Math.}, 130:89--126, 1973.

\bibitem{bers:holomorphic-families}
L.~Bers.
\newblock Holomorphic families of isomorphisms of {M}\"obius groups.
\newblock {\em J. Math. Kyoto Univ.}, 26(1):73--76, 1986.

\bibitem{bonahon:geodesic-currents}
F.~Bonahon.
\newblock The geometry of {T}eichm\"uller space via geodesic currents.
\newblock {\em Invent. Math.}, 92(1):139--162, 1988.

\bibitem{bonahon:shearing}
F.~Bonahon.
\newblock Shearing hyperbolic surfaces, bending pleated surfaces and
  {T}hurston's symplectic form.
\newblock {\em Ann. Fac. Sci. Toulouse Math. (6)}, 5(2):233--297, 1996.

\bibitem{bonahon:variations}
F.~Bonahon.
\newblock Variations of the boundary geometry of {$3$}-dimensional hyperbolic
  convex cores.
\newblock {\em J. Differential Geom.}, 50(1):1--24, 1998.

\bibitem{brock-bromberg}
J.~Brock and K.~Bromberg.
\newblock On the density of geometrically finite {K}leinian groups.
\newblock {\em Acta Math.}, 192(1):33--93, 2004.

\bibitem{bbcm:bumping}
J.~Brock, K.~Bromberg, R.~Canary, and Y.~Minsky.
\newblock In preparation.

\bibitem{bromberg:degenerate}
K.~Bromberg.
\newblock Projective structures with degenerate holonomy and the {B}ers density
  conjecture.
\newblock {\em Ann. of Math. (2)}, 166(1):77--93, 2007.

\bibitem{bromberg-holt:projective}
K.~Bromberg and J.~Holt.
\newblock Bumping of exotic projective structures.
\newblock Preprint.

\bibitem{bromberg-holt:bumping}
K.~Bromberg and J.~Holt.
\newblock Self-bumping of deformation spaces of hyperbolic 3-manifolds.
\newblock {\em J. Differential Geom.}, 57(1):47--65, 2001.

\bibitem{canary:bumponomics}
R.~Canary.
\newblock Introductory bumponomics: {T}he topology of deformation spaces of
  hyperbolic 3-manifolds.
\newblock Preprint, 2007.

\bibitem{ceg:notes-on-notes}
R.~Canary, D.~Epstein, and P.~Green.
\newblock Notes on notes of {T}hurston.
\newblock In {\em Analytical and geometric aspects of hyperbolic space
  (Coventry/Durham, 1984)}, volume 111 of {\em London Math. Soc. Lecture Note
  Ser.}, pages 3--92. Cambridge Univ. Press, Cambridge, 1987.

\bibitem{casson-bleiler}
A.~Casson and S.~Bleiler.
\newblock {\em Automorphisms of surfaces after {N}ielsen and {T}hurston},
  volume~9 of {\em London Mathematical Society Student Texts}.
\newblock Cambridge University Press, Cambridge, 1988.

\bibitem{chuckrow}
V.~Chuckrow.
\newblock On {S}chottky groups with applications to {K}leinian groups.
\newblock {\em Ann. of Math. (2)}, 88:47--61, 1968.

\bibitem{culler-shalen:varieties}
M.~Culler and P.. Shalen.
\newblock Varieties of group representations and splittings of {$3$}-manifolds.
\newblock {\em Ann. of Math. (2)}, 117(1):109--146, 1983.

\bibitem{ddw:r-trees}
G.~Daskalopoulos, S.~Dostoglou, and R.~Wentworth.
\newblock Character varieties and harmonic maps to {${\bf R}$}-trees.
\newblock {\em Math. Res. Lett.}, 5(4):523--533, 1998.

\bibitem{daskalopoulos-wentworth}
G.~Daskalopoulos and R.~Wentworth.
\newblock Harmonic maps and {T}eichm{\"u}ller theory.
\newblock In {\em Handbook of {T}eichm{\"u}ller {T}heory, ({A}. {P}apadopoulos,
  editor), {V}olume {I}}, pages 33--110. EMS Publishing House, Z\"urich, 2007.

\bibitem{donaldson:harmonic}
S.~Donaldson.
\newblock Twisted harmonic maps and the self-duality equations.
\newblock {\em Proc. London Math. Soc. (3)}, 55(1):127--131, 1987.

\bibitem{dumas:bear}
D.~Dumas.
\newblock \removethis{z}{B}ear:~{A} tool for studying {B}ers slices of
  punctured tori.
\newblock Free software, available for download from
  \texttt{http://bear.sourceforge.net/}.

\bibitem{dumas:antipodal}
D.~Dumas.
\newblock Grafting, pruning, and the antipodal map on measured laminations.
\newblock {\em J. Differential Geometry}, 74:93--118, 2006.
\newblock Erratum. 77:175--176, 2007.

\bibitem{dumas:schwarzian}
D.~Dumas.
\newblock The {S}chwarzian derivative and measured laminations on {R}iemann
  surfaces.
\newblock {\em Duke Math. J.}, 140(2):203--243, 2007.

\bibitem{dumas-wolf:grafting}
D.~Dumas and M.~Wolf.
\newblock Projective structures, grafting, and measured laminations.
\newblock {\em Geometry and Topology}, 12(1):351--386, 2008.

\bibitem{earle:variation}
C.~Earle.
\newblock On variation of projective structures.
\newblock In {\em Riemann surfaces and related topics: Proceedings of the 1978
  Stony Brook Conference (State Univ. New York, Stony Brook, N.Y., 1978)},
  volume~97 of {\em Ann. of Math. Stud.}, pages 87--99, Princeton, N.J., 1981.
  Princeton Univ. Press.

\bibitem{eells-lemaire:report}
J.~Eells and L.~Lemaire.
\newblock A report on harmonic maps.
\newblock {\em Bull. London Math. Soc.}, 10(1):1--68, 1978.

\bibitem{eells-lemaire:report2}
J.~Eells and L.~Lemaire.
\newblock Another report on harmonic maps.
\newblock {\em Bull. London Math. Soc.}, 20(5):385--524, 1988.

\bibitem{eells-sampson:existence}
J.~Eells and J.~Sampson.
\newblock Harmonic mappings of {R}iemannian manifolds.
\newblock {\em Amer. J. Math.}, 86:109--160, 1964.

\bibitem{epstein-marden:convex-hulls}
D.~Epstein and A.~Marden.
\newblock Convex hulls in hyperbolic space, a theorem of {S}ullivan, and
  measured pleated surfaces.
\newblock In {\em Analytical and geometric aspects of hyperbolic space
  (Coventry/Durham, 1984)}, volume 111 of {\em London Math. Soc. Lecture Note
  Ser.}, pages 113--253. Cambridge Univ. Press, Cambridge, 1987.

\bibitem{faltings:real-projective}
G.~Faltings.
\newblock Real projective structures on {R}iemann surfaces.
\newblock {\em Compositio Math.}, 48(2):223--269, 1983.

\bibitem{flp}
A.~Fathi, F.~Laudenbach, and V.~Poenaru.
\newblock {\em Travaux de {T}hurston sur les surfaces}, volume~66 of {\em
  Ast\'erisque}.
\newblock Soci\'et\'e Math\'ematique de France, Paris, 1979.
\newblock S\'eminaire Orsay, With an English summary.

\bibitem{frenkel-ben-zvi}
E.~Frenkel and D.~Ben-Zvi.
\newblock {\em Vertex algebras and algebraic curves}, volume~88 of {\em
  Mathematical Surveys and Monographs}.
\newblock American Mathematical Society, Providence, RI, second edition, 2004.

\bibitem{gkm}
D.~Gallo, M.~Kapovich, and A.~Marden.
\newblock The monodromy groups of {S}chwarzian equations on closed {R}iemann
  surfaces.
\newblock {\em Ann. of Math. (2)}, 151(2):625--704, 2000.

\bibitem{gardiner:book}
F.~Gardiner.
\newblock {\em Teichm\"uller theory and quadratic differentials}.
\newblock Pure and Applied Mathematics. John Wiley \& Sons Inc., New York,
  1987.

\bibitem{goldman:symplectic-nature}
W.~Goldman.
\newblock The symplectic nature of fundamental groups of surfaces.
\newblock {\em Adv. in Math.}, 54(2):200--225, 1984.

\bibitem{goldman:fuchsian-holonomy}
W.~Goldman.
\newblock Projective structures with {F}uchsian holonomy.
\newblock {\em J. Differential Geom.}, 25(3):297--326, 1987.

\bibitem{goldman:geometric-structures}
W.~Goldman.
\newblock Geometric structures on manifolds and varieties of representations.
\newblock In {\em Geometry of group representations (Boulder, CO, 1987)},
  volume~74 of {\em Contemp. Math.}, pages 169--198. Amer. Math. Soc.,
  Providence, RI, 1988.

\bibitem{goldman:topological}
W.~Goldman.
\newblock Topological components of spaces of representations.
\newblock {\em Invent. Math.}, 93(3):557--607, 1988.

\bibitem{goldman:ergodic-theory}
W.~Goldman.
\newblock Ergodic theory on moduli spaces.
\newblock {\em Ann. of Math. (2)}, 146(3):475--507, 1997.

\bibitem{gunning:riemann-surfaces}
R.~Gunning.
\newblock {\em Lectures on {R}iemann surfaces}.
\newblock Princeton Mathematical Notes. Princeton University Press, Princeton,
  N.J., 1966.

\bibitem{gunning:vector-bundles}
R.~Gunning.
\newblock {\em Lectures on vector bundles over {R}iemann surfaces}.
\newblock University of Tokyo Press, Tokyo, 1967.

\bibitem{gunning:special-coordinate}
R.~Gunning.
\newblock Special coordinate coverings of {R}iemann surfaces.
\newblock {\em Math. Ann.}, 170:67--86, 1967.

\bibitem{gunning:affine-projective}
R.~Gunning.
\newblock Affine and projective structures on {R}iemann surfaces.
\newblock In {\em Riemann surfaces and related topics: Proceedings of the 1978
  Stony Brook Conference (State Univ. New York, Stony Brook, N.Y., 1978)},
  volume~97 of {\em Ann. of Math. Stud.}, pages 225--244, Princeton, N.J.,
  1981. Princeton Univ. Press.

\bibitem{hejhal:monodromy}
D.~Hejhal.
\newblock Monodromy groups and linearly polymorphic functions.
\newblock {\em Acta Math.}, 135(1):1--55, 1975.

\bibitem{hejhal:monodromy-and-poincare}
D.~Hejhal.
\newblock Monodromy groups and {P}oincar\'e series.
\newblock {\em Bull. Amer. Math. Soc.}, 84(3):339--376, 1978.

\bibitem{heusener-porti}
M.~Heusener and J.~Porti.
\newblock The variety of characters in {${\rm PSL}\sb 2(\Bbb C)$}.
\newblock {\em Bol. Soc. Mat. Mexicana (3)}, 10(Special Issue):221--237, 2004.

\bibitem{holt:bumping}
J.~Holt.
\newblock Some new behaviour in the deformation theory of {K}leinian groups.
\newblock {\em Comm. Anal. Geom.}, 9(4):757--775, 2001.

\bibitem{hubbard:monodromy}
J.~Hubbard.
\newblock The monodromy of projective structures.
\newblock In {\em Riemann surfaces and related topics: Proceedings of the 1978
  Stony Brook Conference (State Univ. New York, Stony Brook, N.Y., 1978)},
  volume~97 of {\em Ann. of Math. Stud.}, pages 257--275, Princeton, N.J.,
  1981. Princeton Univ. Press.

\bibitem{hubbard:book}
J.~Hubbard.
\newblock {\em Teichm\"uller theory and applications to geometry, topology, and
  dynamics. {V}ol. 1}.
\newblock Matrix Editions, Ithaca, NY, 2006.

\bibitem{hubbard-masur}
J.~Hubbard and H.~Masur.
\newblock Quadratic differentials and foliations.
\newblock {\em Acta Math.}, 142(3-4):221--274, 1979.

\bibitem{huber:subharmonic}
A.~Huber.
\newblock On subharmonic functions and differential geometry in the large.
\newblock {\em Comment. Math. Helv.}, 32:13--72, 1957.

\bibitem{huber:zum}
A.~Huber.
\newblock Zum potentialtheoretischen {A}spekt der {A}lexandrowschen
  {F}l\"achentheorie.
\newblock {\em Comment. Math. Helv.}, 34:99--126, 1960.

\bibitem{imayoshi-taniguchi}
Y.~Imayoshi and M.~Taniguchi.
\newblock {\em An introduction to {T}eichm\"uller spaces}.
\newblock Springer-Verlag, Tokyo, 1992.
\newblock Translated and revised from the Japanese by the authors.

\bibitem{ito:exotic}
K.~Ito.
\newblock Exotic projective structures and quasi-{F}uchsian space.
\newblock {\em Duke Math. J.}, 105(2):185--209, 2000.

\bibitem{ito:schottky}
K.~Ito.
\newblock Schottky groups and {B}ers boundary of {T}eichm\"uller space.
\newblock {\em Osaka J. Math.}, 40(3):639--657, 2003.

\bibitem{ito:survey}
K.~Ito.
\newblock Grafting and components of quasi-{F}uchsian projective structures.
\newblock In {\em Spaces of {K}leinian groups}, volume 329 of {\em London Math.
  Soc. Lecture Note Ser.}, pages 355--373. Cambridge Univ. Press, Cambridge,
  2006.

\bibitem{ito:exotic2}
K.~Ito.
\newblock Exotic projective structures and quasi-{F}uchsian space, {II}.
\newblock {\em Duke Math. J.}, 140(1):85--109, 2007.

\bibitem{jorgensen:discrete}
T.~J{\o}rgensen.
\newblock On discrete groups of {M}\"obius transformations.
\newblock {\em Amer. J. Math.}, 98(3):739--749, 1976.

\bibitem{jost:compact-riemann-surfaces}
J.~Jost.
\newblock {\em Compact {R}iemann surfaces}.
\newblock Universitext. Springer-Verlag, Berlin, third edition, 2006.

\bibitem{kamishima:not-surjective}
Y.~Kamishima.
\newblock Conformally flat manifolds whose development maps are not surjective.
  {I}.
\newblock {\em Trans. Amer. Math. Soc.}, 294(2):607--623, 1986.

\bibitem{kamishima-tan:grafting}
Y.~Kamishima and S.~Tan.
\newblock Deformation spaces on geometric structures.
\newblock In {\em Aspects of low-dimensional manifolds}, volume~20 of {\em Adv.
  Stud. Pure Math.}, pages 263--299. Kinokuniya, Tokyo, 1992.

\bibitem{kapovich:monodromy}
M.~Kapovich.
\newblock On monodromy of complex projective structures.
\newblock {\em Invent. Math.}, 119(2):243--265, 1995.

\bibitem{kapovich:book}
M.~Kapovich.
\newblock {\em Hyperbolic manifolds and discrete groups}, volume 183 of {\em
  Progress in Mathematics}.
\newblock Birkh\"auser Boston Inc., Boston, MA, 2001.

\bibitem{klein:ausgewahlte}
F.~{K}lein.
\newblock {\em {A}usgew\"ahlte {K}apital aus der {T}heorie der linearen
  {D}ifferentialgleichungen zweiter {O}rdnung}, volume~1.
\newblock G\"ottingen, 1891.

\bibitem{klein}
F.~Klein.
\newblock {\em Vorlesungen {\"U}ber die {H}ypergeometrische {F}unktion}.
\newblock Springer-Verlag, Berlin, 1933.

\bibitem{kojima:survey}
S.~Kojima.
\newblock {C}ircle packing and {T}eichm{\"u}ller spaces.
\newblock To appear in \emph{{H}andbook of {T}eichm{\"u}ller {T}heory, ({A}.
  {P}apadopoulos, editor), {V}olume {II}}, EMS Publishing House, 2009.

\bibitem{kojima-mizushima-tan}
S.~Kojima, S.~Mizushima, and S.~Tan.
\newblock Circle packings on surfaces with projective structures: a survey.
\newblock In {\em Spaces of {K}leinian groups}, volume 329 of {\em London Math.
  Soc. Lecture Note Ser.}, pages 337--353. Cambridge Univ. Press, Cambridge,
  2006.

\bibitem{korevaar-schoen:3}
N.~Korevaar and R.~Schoen.
\newblock Global existence theorems for harmonic maps: finite rank spaces and
  an approach to rigidity for smooth actions.
\newblock Preprint.

\bibitem{korevaar-schoen:1}
N.~Korevaar and R.~Schoen.
\newblock Sobolev spaces and harmonic maps for metric space targets.
\newblock {\em Comm. Anal. Geom.}, 1(3-4):561--659, 1993.

\bibitem{korevaar-schoen:2}
N.~Korevaar and R.~Schoen.
\newblock Global existence theorems for harmonic maps to non-locally compact
  spaces.
\newblock {\em Comm. Anal. Geom.}, 5(2):333--387, 1997.

\bibitem{kra:deformations}
I.~Kra.
\newblock Deformations of {F}uchsian groups.
\newblock {\em Duke Math. J.}, 36:537--546, 1969.

\bibitem{kra:deformations2}
I.~Kra.
\newblock Deformations of {F}uchsian groups. {II}.
\newblock {\em Duke Math. J.}, 38:499--508, 1971.

\bibitem{kra:generalization}
I.~Kra.
\newblock A generalization of a theorem of {P}oincar\'e.
\newblock {\em Proc. Amer. Math. Soc.}, 27:299--302, 1971.

\bibitem{kra-maskit:remarks}
I.~Kra and B.~Maskit.
\newblock Remarks on projective structures.
\newblock In {\em Riemann surfaces and related topics: Proceedings of the 1978
  Stony Brook Conference (State Univ. New York, Stony Brook, N.Y., 1978)},
  volume~97 of {\em Ann. of Math. Stud.}, pages 343--359, Princeton, N.J.,
  1981. Princeton Univ. Press.

\bibitem{kulkarni-pinkall:metric}
R.~Kulkarni and U.~Pinkall.
\newblock A canonical metric for {M}\"obius structures and its applications.
\newblock {\em Math. Z.}, 216(1):89--129, 1994.

\bibitem{lehto:univalent}
O.~Lehto.
\newblock {\em Univalent functions and {T}eichm\"uller spaces}, volume 109 of
  {\em Graduate Texts in Mathematics}.
\newblock Springer-Verlag, New York, 1987.

\bibitem{levitt}
G.~Levitt.
\newblock Foliations and laminations on hyperbolic surfaces.
\newblock {\em Topology}, 22(2):119--135, 1983.

\bibitem{luo}
F.~Luo.
\newblock Monodromy groups of projective structures on punctured surfaces.
\newblock {\em Invent. Math.}, 111(3):541--555, 1993.

\bibitem{marden:outer-circles}
A.~Marden.
\newblock {\em Outer {C}ircles: {A}n {I}ntroduction to {H}yperbolic
  {$3$}-{M}anifolds}.
\newblock Cambridge University Press, Cambridge, 2007.

\bibitem{maskit:grafting}
B.~Maskit.
\newblock On a class of {K}leinian groups.
\newblock {\em Ann. Acad. Sci. Fenn. Ser. A I No.}, 442:8, 1969.

\bibitem{matsuzaki}
K.~Matsuzaki.
\newblock The interior of discrete projective structures in the {B}ers fiber.
\newblock {\em Ann. Acad. Sci. Fenn. Math.}, 32(1):3--12, 2007.

\bibitem{matsuzaki-taniguchi}
K.~Matsuzaki and M.~Taniguchi.
\newblock {\em Hyperbolic manifolds and {K}leinian groups}.
\newblock Oxford Mathematical Monographs. The Clarendon Press Oxford University
  Press, New York, 1998.
\newblock , Oxford Science Publications.

\bibitem{mcmullen:complex-earthquakes}
C.~McMullen.
\newblock Complex earthquakes and {T}eichm\"uller theory.
\newblock {\em J. Amer. Math. Soc.}, 11(2):283--320, 1998.

\bibitem{mese:between-surfaces}
C.~Mese.
\newblock Harmonic maps between surfaces and {T}eichm\"uller spaces.
\newblock {\em Amer. J. Math.}, 124(3):451--481, 2002.

\bibitem{mese:harmonic-maps-into-singular}
C.~Mese.
\newblock Harmonic maps into spaces with an upper curvature bound in the sense
  of {A}lexandrov.
\newblock {\em Math. Z.}, 242(4):633--661, 2002.

\bibitem{mese:uniqueness}
C.~Mese.
\newblock Uniqueness theorems for harmonic maps into metric spaces.
\newblock {\em Commun. Contemp. Math.}, 4(4):725--750, 2002.

\bibitem{morgan-shalen:valuations-trees}
J.~Morgan and P.~Shalen.
\newblock Valuations, trees, and degenerations of hyperbolic structures. {I}.
\newblock {\em Ann. of Math. (2)}, 120(3):401--476, 1984.

\bibitem{nehari:schwarzian}
Z.~Nehari.
\newblock The {S}chwarzian derivative and schlicht functions.
\newblock {\em Bull. Amer. Math. Soc.}, 55:545--551, 1949.

\bibitem{ohshika:bumping}
K.~Ohshika.
\newblock Divergence, exotic convergence, and self-bumping in quasi-{F}uchsian
  spaces.
\newblock In preparation.

\bibitem{osgood-stowe}
B.~Osgood and D.~Stowe.
\newblock The {S}chwarzian derivative and conformal mapping of {R}iemannian
  manifolds.
\newblock {\em Duke Math. J.}, 67(1):57--99, 1992.

\bibitem{penner-harer}
R.~Penner and J.~Harer.
\newblock {\em Combinatorics of train tracks}, volume 125 of {\em Annals of
  Mathematics Studies}.
\newblock Princeton University Press, Princeton, NJ, 1992.

\bibitem{poincare:groups-linear-equations}
H.~Poincar{\'e}.
\newblock Sur les groupes des \'equations lin\'eaires.
\newblock {\em Acta Math.}, 4(1):201--312, 1884.

\bibitem{pommerenke:univalent}
C.~Pommerenke.
\newblock {\em Univalent functions}.
\newblock Vandenhoeck \& Ruprecht, G\"ottingen, 1975.
\newblock With a chapter on quadratic differentials by Gerd Jensen, Studia
  Mathematica/Mathematische Lehrb\"ucher, Band XXV.

\bibitem{bcr}
A.~Rapinchuk, V.~Benyash-Krivetz, and V.~Chernousov.
\newblock Representation varieties of the fundamental groups of compact
  orientable surfaces.
\newblock {\em Israel J. Math.}, 93:29--71, 1996.

\bibitem{reshetnyak:survey}
Y.~Reshetnyak.
\newblock Two-dimensional manifolds of bounded curvature.
\newblock In {\em Geometry. {IV}: {N}onregular {R}iemannian geometry},
  Encyclopaedia of Mathematical Sciences, pages 3--163. Springer-Verlag,
  Berlin, 1993.
\newblock Translation of {\it Geometry, 4 (Russian)}, Akad.\ Nauk SSSR,
  Vsesoyuz.\ Inst.\ Nauchn.\ i Tekhn.\ Inform., Moscow, 1989.

\bibitem{reshetnyak:conformal}
Y.~Reshetnyak.
\newblock On the conformal representation of {A}lexandrov surfaces.
\newblock In {\em Papers on analysis}, volume~83 of {\em Rep. Univ.
  Jyv\"askyl\"a Dep. Math. Stat.}, pages 287--304. Univ. Jyv\"askyl\"a,
  Jyv\"askyl\"a, 2001.

\bibitem{riemann:hypergeometric}
B.~Riemann.
\newblock {V}orlesungen {\"u}ber die hypergeometrische {R}eihe.
\newblock In {\em Gesammelte mathematische {W}erke, wissenschaftlicher
  {N}achlass und {N}achtr\"age}, pages 667--692. Springer-Verlag, Berlin, 1990.

\bibitem{scannell:thesis}
K.~Scannell.
\newblock {\em Flat conformal structures and causality in de {S}itter
  manifolds}.
\newblock PhD thesis, University of California, Los Angeles, 1996.

\bibitem{scannell:desitter}
K.~Scannell.
\newblock Flat conformal structures and the classification of de {S}itter
  manifolds.
\newblock {\em Comm. Anal. Geom.}, 7(2):325--345, 1999.

\bibitem{scannell-wolf:grafting}
K.~Scannell and M.~Wolf.
\newblock The grafting map of {T}eichm\"uller space.
\newblock {\em J. Amer. Math. Soc.}, 15(4):893--927 (electronic), 2002.

\bibitem{schoen:analytic-aspects}
R.~Schoen.
\newblock Analytic aspects of the harmonic map problem.
\newblock In {\em Seminar on nonlinear partial differential equations
  (Berkeley, Calif., 1983)}, volume~2 of {\em Math. Sci. Res. Inst. Publ.},
  pages 321--358. Springer, New York, 1984.

\bibitem{sharpe}
R.~Sharpe.
\newblock {\em Differential geometry}, volume 166 of {\em Graduate Texts in
  Mathematics}.
\newblock Springer-Verlag, New York, 1997.
\newblock Cartan's generalization of Klein's Erlangen program, With a foreword
  by S. S. Chern.

\bibitem{shiga}
H.~Shiga.
\newblock Projective structures on {R}iemann surfaces and {K}leinian groups.
\newblock {\em J. Math. Kyoto Univ.}, 27(3):433--438, 1987.

\bibitem{shiga-tanigawa}
H.~Shiga and H.~Tanigawa.
\newblock Projective structures with discrete holonomy representations.
\newblock {\em Trans. Amer. Math. Soc.}, 351(2):813--823, 1999.

\bibitem{sullivan:qc2}
D.~Sullivan.
\newblock Quasiconformal homeomorphisms and dynamics. {II}. {S}tructural
  stability implies hyperbolicity for {K}leinian groups.
\newblock {\em Acta Math.}, 155(3-4):243--260, 1985.

\bibitem{sullivan-thurston}
D.~Sullivan and W.~Thurston.
\newblock Manifolds with canonical coordinate charts: some examples.
\newblock {\em Enseign. Math. (2)}, 29(1-2):15--25, 1983.

\bibitem{tanigawa:grafting}
H.~Tanigawa.
\newblock Grafting, harmonic maps and projective structures on surfaces.
\newblock {\em J. Differential Geom.}, 47(3):399--419, 1997.

\bibitem{tanigawa:divergence}
H.~Tanigawa.
\newblock Divergence of projective structures and lengths of measured
  laminations.
\newblock {\em Duke Math. J.}, 98(2):209--215, 1999.

\bibitem{thurston:hyp2}
W.~Thurston.
\newblock {H}yperbolic {S}tructures on 3-manifolds, {II}: {S}urface groups and
  3-manifolds which fiber over the circle.
\newblock Preprint.

\bibitem{thurston:notes}
W.~Thurston.
\newblock Geometry and topology of three-manifolds.
\newblock Princeton lecture notes, 1979.

\bibitem{thurston:earthquakes}
W.~Thurston.
\newblock Earthquakes in two-dimensional hyperbolic geometry.
\newblock In {\em Low-dimensional topology and {K}leinian groups
  (Coventry/Durham, 1984)}, volume 112 of {\em London Math. Soc. Lecture Note
  Ser.}, pages 91--112. Cambridge Univ. Press, Cambridge, 1986.

\bibitem{thurston:minimal-stretch}
W.~Thurston.
\newblock Minimal stretch maps between hyperbolic surfaces.
\newblock Unpublished preprint, 1986.

\bibitem{thurston:zippers}
W.~Thurston.
\newblock Zippers and univalent functions.
\newblock In {\em The Bieberbach conjecture (West Lafayette, Ind., 1985)},
  volume~21 of {\em Math. Surveys Monogr.}, pages 185--197. Amer. Math. Soc.,
  Providence, RI, 1986.

\bibitem{thurston:bulletin}
W.~Thurston.
\newblock On the geometry and dynamics of diffeomorphisms of surfaces.
\newblock {\em Bull. Amer. Math. Soc. (N.S.)}, 19(2):417--431, 1988.

\bibitem{wolf:thesis}
M.~Wolf.
\newblock The {T}eichm\"uller theory of harmonic maps.
\newblock {\em J. Differential Geom.}, 29(2):449--479, 1989.

\bibitem{wolf:realizing-measured-foliations}
M.~Wolf.
\newblock On realizing measured foliations via quadratic differentials of
  harmonic maps to {$\bold R$}-trees.
\newblock {\em J. Anal. Math.}, 68:107--120, 1996.

\bibitem{wright:maskit-boundary}
D.~Wright.
\newblock The shape of the boundary of the {T}eichm{\"u}ller space of
  once-punctured tori in {M}askit's embedding.
\newblock Unpublished preprint, 1987.

\end{thebibliography}
\end{document}